\documentclass[12pt]{article}

\usepackage{amssymb,amsmath,amsthm,bbm,enumitem,xcolor}
\usepackage[a4paper, total={6.3in, 9.2in}]{geometry}
\usepackage{setspace}
\usepackage{hyperref}
\usepackage{authblk}
\usepackage{graphicx,caption,subcaption}
\usepackage[title]{appendix}
 \usepackage{mathtools}
\usepackage{algorithm,algpseudocode}
\usepackage{tikz}
\usetikzlibrary{shapes.geometric, arrows, decorations.markings,shapes.multipart,cd,positioning,patterns,snakes}
\tikzstyle{startstop} = [rectangle, rounded corners, minimum width=2cm, minimum height=1cm,text centered, draw=black]
\tikzstyle{vecArrow} = [thick, decoration={markings,mark=at position
   1 with {\arrow[semithick]{open triangle 60}}},
   double distance=1.4pt, shorten >= 5.5pt,
   preaction = {decorate},
   postaction = {draw,line width=1.4pt, white,shorten >= 4.5pt}]
\tikzstyle{innerWhite} = [semithick, white,line width=1.4pt, shorten >= 4.5pt]
   \tikzset{|/.tip={Bar[width=.8ex,round]}}

\newtheorem{theorem}{Theorem}[section]
\theoremstyle{definition}
\newtheorem{definition}{Definition}[section]
\theoremstyle{definition}

\theoremstyle{definition}
\newtheorem{assumption}{Assumption}[section]
\theoremstyle{definition}
\newtheorem{proposition}{Proposition}[section]
\theoremstyle{definition}

\theoremstyle{remark}
\newtheorem{remark}{Remark}[section]
\theoremstyle{definition}
\newtheorem{lemma}{Lemma}[section]

\newtheorem{assumex}{Assumption}
\newenvironment{assume}[1]
 {\renewcommand\theassumex{#1}\assumex}
 {\endassumex}

\theoremstyle{definition}
\newtheorem{problem}{Problem}
\newtheorem{modelx}{Model}
\newenvironment{model}[1]
 {\renewcommand\themodelx{#1}\modelx}
 {\endmodelx}

\DeclareMathOperator*{\argmin}{arg\,min}

\usepackage[]{authblk}

\setlist[itemize]{leftmargin=0.4cm,labelindent=\parindent}

\title{Mean Field Analysis of Two-Party Governance: Competition versus Cooperation among Leaders}
\date{\today}

\author[$\dagger$]{Dantong Chu}
\author[$\star$]{Kenneth Tsz Hin Ng} 
\author[$\dagger$]{Sheung Chi Phillip Yam} 
\author[$\ddag$]{Harry Zheng} 

 \affil[$\dagger$]{Department of Statistics, The Chinese University of Hong Kong, Shatin, N.T., Hong Kong. Email: dtchu@link.cuhk.edu.hk, scpyam@sta.cuhk.edu.hk }
 \affil[$\star$]{Department of Mathematics, University of Illinois at Urbana-Champaign, Illinois, US. \protect \\   Email: tszhinn2@illinois.edu}
 \affil[$\ddag$]{Department of Mathematics, Imperial College, London, UK. \protect \\  Email: h.zheng@imperial.ac.uk}

\begin{document}
\maketitle

\begin{abstract}                          
     This article studies linear-quadratic Stackelberg games between two dominating players (or equivalently, leaders) and a large group of followers, each of {\color{black}whom interacts} under a mean field game (MFG) framework. Unlike the conventional major-minor player game, the mean field term herein is endogenously {\color{black}affected} by the two leaders simultaneously. These homogeneous followers are {\color{black}non-}cooperative, whereas the two leaders can either  compete or cooperate with each other, which are respectively formulated as a Nash and a Pareto game. 
    The complete solutions of the leader-follower game can be expressed in terms of the solutions of some non-symmetric Riccati equations. {\color{black}Notably, our analysis suggests that both} modes of interactions between leaders has their own merits and neither {\color{black}of them} is always more favourable to the community {\color{black}of} followers. {\color{black}In our knowledge,}  a comparative {\color{black}study} of the effect of different modes of governance on the society  is relatively rare in the {\color{black}existing} literature,  we {\color{black}here} provide {\color{black}its first} preliminary  quantitative analysis; under a broad class of practically relevant models, we provide sufficient conditions to decide whether cooperation or competition between leaders is more favourable to the followers. {\color{black}Being in common with modern folklore}, the relative merits of the two {\color{black}Stackelberg} games depend on whether the interests between the two leaders and the followers align {\color{black}among themselves}. {\color{black}Representative} numerical examples are also {\color{black}supplemented}. 	
\end{abstract}

\section{Introduction}

    A mean field game is a macroscopic framework in studying the optimal decisions of players in a stochastic differential game (SDG) when the number of players involved is too large to be handled explicitly, due to the intricate interactions among them. For a detailed introduction to the subject, consult the monographs by  Bensoussan et al. \cite{BFY}, Carmona and Delarue \cite{carmona:2018}, and the references therein.  The simplicity and competence of mean field theory in modeling  multi-agent systems have prompted  fruitful interdisciplinary applications, for instance, in finance \cite{cardaliaguet2018mean,casgrain2019algorithmic,firoozi2017optimal,HAN2022}; economics \cite{achdou2014partial,achdou2022income}; industrial engineering \cite{alasseur2020extended,de2019mean,kizilkale2019integral}; and more recently in machine learning \cite{ruthotto2020machine} and studies on the pandemic \cite{elie2020contact,hubert2020incentives}.
  
\nocite{BSYY}\nocite{CAR_PROB}

In this article, we consider a linear-quadratic mean field game between two heterogeneous dominating players (or equivalently, leaders) and a group of homogeneous followers. This multi-leader multi-follower environment is ubiquitous in many interdisciplinary studies, for instance, a two-party {\color{black}campaign} voting system in political science, where the leaders and followers are respectively the major political parties and the voters \cite{alesina:1987,alesina:1988}; the interactions between {\color{black}oligopolistic} price setters (large suppliers) and price takers (smaller firms or countries) in resources and energy markets \cite{ALROUSAN201828,alasseur2020extended}; the resource management problem between {\color{black}oligopolistic} service providers and individual miners in block-chain mining \cite{jiang:block:chain:2022}, etc. The notion of dominating players of a MFG was first introduced in Bensoussan et al. \cite{BCY}.  In contrast to the notion of major players in \cite{Buckdahn:2014,carmona:major,LASRY2018886,csen2014mean}, the mean field term herein is endogenously determined by the leaders. 
This manifests the governing nature of the leaders who publicly {\color{black}impact} the society. Also, followers need to comply with the initiatives taken by the leaders. On the other hand, the mean field term (see the definition \eqref{meanfieldfixpoint}) also serves as a feedback mechanism that influences leaders' decisions, whose importance in maintaining the stability of a financial system is highlighted in Cooper \cite{cooper2008origin}. Under a single-leader game, Bensoussan et al. \cite{2020:parametrized,BCLY}  considered a model where followers are subject to delayed effects from the state dynamics of themselves and the leader;  Fu and Horst \cite{fu:horst:2020} studied a mean field type leader-follower control problem and its associated McKean-Vlasov forward backward stochastic differential equations (FBSDEs).



\nocite{bensoussan:sethi:2019}

 {\color{black} Here, w}e consider two different modes of interactions between the two heterogeneous leaders, namely, competition and cooperation,  which are formulated as a Nash and a Pareto game, respectively. A Pareto equilibrium is an optimal solution in  allocations of resources such that no player could be better off without hurting other players' benefits if he moves away from the equilibrium. Under the SDG framework, the Pareto equilibrium is achieved by minimizing the sum of  individuals' cost functionals; see, e.g.,  Elliot \cite{elliott:game}  for details. In the context of mean field games or mean field type control, Pareto equilibria are often considered in social optimization problems, see, e.g.,  Wang et al. \cite{WANG:social}. On the other hand, a Nash equilibrium under a non-cooperative environment is the situation where each player is well aware of others' decision, and no one will be better off by only changing his own strategy from the equilibrium; see, e.g., Li and Sethi \cite{tao:sethi} for a general review of non-cooperative Stackelberg games.

 {We decompose our {\color{black}overall problem} into three {\color{black}nested} sub-problems {\color{black}as follows}. In Problem \ref{P1}, given the state dynamics of the two leaders and the mean field term, we solve for the representative follower's optimal control and the $\epsilon$-Nash equilibrium for the community.
In Problem \ref{P2}, we find the mean field term by solving a fixed point problem, which is set to equal the conditional expectation of the state dynamics of the follower given the filtration generated by the Wiener processes of the two leaders. In Problem \ref{P3}, we define and solve  the corresponding competitive Nash game and the cooperative Pareto game between the two leaders, where the solution can be characterized by systems of FBSDEs using the stochastic maximum principle; {\color{black}as a result}, we express explicitly the complete solutions of the leader-follower games by some non-symmetric Riccati differential equations and address their well-posedness.}


The relative merits of competition and cooperation among individuals are often discussed in game theory, see {\color{black}for example,} Cohen \cite{Cohen9724}. {\color{black}In our knowledge, due to the mathematical complexity, it is still rare in the literature that compares different interactions among two leaders, competition versus cooperation, about their effects on another third party, let alone  a group of  followers.} {\color{black}Here}, we  provide {\color{black} the ever first} {\color{black}systematic}  quantitative discussion {\color{black}on this matter.}
{\color{black}According to our later analysis},  neither cooperation or competition between leaders is always ideal for the society. While leaders ought to cooperate to share the limited resources in a society, sole cooperation would however lessen their incentives to innovate; this is supported by the findings in political science,  where electoral competitions between major parties can provide incentives to bolster educational and economic development, see {\color{black} for example}, \cite{de2005logic,gamm_kousser_2021}.  To illustrate this {\color{black}quantitatively}, in our Theorem \ref{thm:compare}{\color{black},} under the conditions therein,  we compare the cost functionals of the followers under the two games, which can be expressed by the solutions of systems of Lyapunov {\color{black} and} Riccati equations.   Since there are no general comparison theorems for non-symmetric Riccati equations, direct comparison of the {\color{black}same} cost functionals {\color{black}but} under general setting is highly intractable, and the discussion of this topic is almost absent in the literature.  Herein, we provide a first symmetric analysis where followers and leaders depend on each other in a fair manner. Under {\color{black}the} model (Model \ref{degeneratedmodel}), we obtain decisive conditions in Theorem  \ref{thm:compare} under which cooperation or competition among leaders is favourable to the followers. Intutively, the conditions {\color{black} are grounded} on whether the interests between the two leaders and the followers are aligned {\color{black}among themselves}, which are quantified by the signs of the coefficients in the state dynamics and the cost functionals, along with the {\color{black}magnitude} of the initial conditions. In general, we can extend this comparison result to a broader class of models via a perturbation approach; indeed,  we show  in Theorem \ref{compare_extend} that the decisive conditions in Theorem \ref{thm:compare} also hold for models that are reasonably close to Model \ref{degeneratedmodel}, which is measured by the deviations of the coefficients from the latter.


 



The article is organized as follows. In Section \ref{sec:problem:formulation}, we provide the problem settings and present the $\epsilon$-equilibrium of the follower-game in the presence of the two leaders. The Nash and the Pareto game between the leaders, as well as Problems \ref{P1}-\ref{P3}, are also formalized.  The three sub-problems are  solved in Section \ref{sec:sol}. In Section \ref{sec:compare}, we compare the effect of the two mode of interactions of leaders on the follower and present Theorems \ref{thm:compare} and \ref{compare_extend}. {\color{black}Supportive n}umerical examples  are provided in Section \ref{sec:numerical}. The article is concluded in Section \ref{sec:conclusion}.  The {\color{black}appendices} {\color{black}provide more detailed} proofs of the main results. 



\section{Problem Setting}\label{sec:problem:formulation}
Throughout this article, we consider a complete probability space $(\Omega,\mathcal{F},\mathbb{P})$. We shall call the two heterogeneous leaders as the $\alpha$-leader and the $\beta$-leader, respectively, and use a superscript or a subscript $\alpha$ and $\beta$ to represent {\color{black}their corresponding mathematical} objects.  Let $T>0$ be a fixed terminal time, also let $d_\alpha$, $d_\beta$, $d_1$, $n_\alpha$, $n_\beta$, $n_1$, $m_\alpha$, $m_\beta$ and  $m_1 \in \mathbb{N}^+$. Assume that $W^{\alpha},W^{\beta}$ and $\{{\color{black}W^{1,i}}\}_{i\in\{1,...,N\}}$ are independent Wiener processes taking values in $\mathbb{R}^{d_\alpha}$, $\mathbb{R}^{d_\beta}$ and $\mathbb{R}^{d_1}$, respectively. Let  $\eta_{\alpha}$ and $\eta_{\beta}$ be independent random variables  representing the {\color{black}respective} initial states of the two leaders, which are square-integrable and are independent of the aforementioned Wiener processes. Also let the random variables $\{\eta_1^i\}_{i\in\{1,...,N\}}$, representing the initial states of {\color{black}followers}, be square integrable, independent and identically distributed, which are also independent of  $\eta_\alpha$, $\eta_\beta$, and the aforementioned Wiener processes. Define the filtrations $\mathcal{F}^\alpha=\{\mathcal{F}^\alpha_t\}_{t\geq 0}$, $\mathcal{F}^\beta=\{\mathcal{F}^\beta_t\}_{t\geq 0}$ and $\mathcal{F}^{1,i}=\{{\color{black}\mathcal{F}^{1,i}_t}\}_{t\geq 0}$, $i=1,2\dots,{\color{black}N,}$ by, for any $t\geq 0$, $\mathcal{F}^\alpha_t:=\sigma(\eta_\alpha,W^{\alpha}(s):s\leq t)$, $\mathcal{F}^\beta_t:=\sigma(\eta_\beta,{\color{black}W^{\beta}}(s):s\leq t)$ and  $\mathcal{F}^{1,i}_t:=\sigma(\eta_1^i,{\color{black}W^{1,i}}(s):s\leq t)$.
We now introduce the dynamical systems of the two leaders and the $N$ followers. For $\gamma = \alpha,\beta$, the empirical state dynamics for the $\gamma$-leader is  described by the following stochastic differential equations on $\mathbb{R}^{n_{\gamma}}$:  %
  \begin{equation} \label{eq:leader:state:empirical}
    \begin{cases}
dy_{\gamma}(t)=\Big(A_{\gamma}y_{\alpha}(t)+B_{\gamma}y_{\beta}(t)+C_{\gamma}\frac{\sum_{j=1	
	}^{N}{\color{black}y_{1,j}}(t)}{N}+D_{\gamma}v_{\gamma}(t)\Big)dt+\sigma_{\gamma}dW^{\gamma}(t);\\
y_{\gamma}(0)=\eta_{\gamma},
    \end{cases}
\end{equation}
where $\{{\color{black}y_{1,i}}(t)\}_{t\geq 0}$, $i=1,2,\dots,N$, are the state dynamics of the homogeneous followers taking values in $\mathbb{R}^{n_1}$, which are given by: 
\begin{equation}
\label{y1SDE}
    \begin{cases}
    d{\color{black}y_{1,i}}(t)=\Big(A_1{\color{black}y_{1,i}}(t)+B_1\frac{\sum_{j=1,	
    		j\neq i}^{N}{\color{black}y_{1,j}}(t)}{N-1}+C_1y_{\alpha}(t) +D_1y_{\beta}(t)+E_1{\color{black}v_{1,i}}(t)\Big)dt+\sigma_1d{\color{black}W^{1,i}}(t);\\
    {\color{black}y_{1,i}}(0)=\eta_1^i.
    \end{cases}
\end{equation}
Here, {\color{black}the controls} $v_{\alpha},v_{\beta}$ and ${\color{black}v_{1,i}}$  are respectively $\mathbb{R}^{m_\alpha}$, $\mathbb{R}^{m_\beta}$ and $\mathbb{R}^{m_1}$-valued {\color{black}square-integrable} process  for two leaders and the $i$-th follower. The matrices $A_\alpha,B_\alpha,C_\alpha,D_\alpha, \sigma_\alpha$, $A_\beta,B_\beta,C_\beta,D_\beta,\sigma_\beta,A_1,B_1,C_1$, $D_1$, $E_1,\sigma_1$ are assumed to be constant with appropriate dimensions, and $D_\alpha,D_\beta,E_1$ are assumed to have full column ranks. We also assume that $v_{\alpha},v_{\beta},v_{1,i}\in\mathcal{F}^{\alpha}\vee\mathcal{F}^{\beta}{\color{black}\bigvee_{j=1}^N\mathcal{F}^{1,j}}$. {\color{black} Let ${\bf v}_1 := (v_{1,1},\dots,v_{1,N})$ and $\mathbf{v}_{1,-i}:=(v_{1,1},$ $\dots,v_{1,i-1},v_{1,i+1},\dots,$ $v_{1,N})$}. For $\gamma = \alpha,\beta$, the $\gamma$-leader and the $i$-th follower aim to minimize the following cost functionals, $\mathcal{J}_\gamma=\mathcal{J}_\gamma(v_{\gamma};v_{\alpha+\beta-\gamma}{\color{black}, {\bf v}_1})$ 
and {\color{black}$\mathcal{J}_{1,i}=\mathcal{J}_{1,i}({\color{black}v_{1,i}; v_\alpha,v_\beta, {\bf v}_{1,-i}})$}, respectively:  
  \begin{align*}
     \mathcal{J}_{\gamma} 
   &= \mathbb{E}\bigg[\int_0^T \bigg(\bigg|y_{\gamma}(t)-F_{\gamma}\frac{\sum_{j=1}^{N}{\color{black}y_{1,j}}(t)}{N}-G_\gamma y_{\alpha+\beta-\gamma}(t)   -M_{\gamma}\bigg|_{Q_{\gamma}}^2+|v_{\gamma}(t)|_{R_{\gamma}}^2\bigg)dt \\&\qquad  +\bigg|y_{\gamma}(T)   -\bar{F}_{\gamma}\frac{\sum_{j=1}^{N}{\color{black}y_{1,j}}(T)}{N}   -\bar{G}_{\gamma}y_{\alpha+\beta-\gamma}(T)-\bar{M}_{\gamma}\bigg|_{\bar{Q}_{\gamma}}^2\bigg];\\
 {\color{black}\mathcal{J}_{1,i}} &= \mathbb{E}\bigg[\int_0^T\bigg(\bigg|{\color{black}y_{1,i}}(t)-F_1\dfrac{\sum_{j=1,j\neq i}^{N}{\color{black}y_{1,j}}(t)}{N-1}  -G_1y_{\alpha}(t)  -H_1y_{\beta}(t)-M_1\bigg|_{Q_1}^2
 + |{\color{black}v_{1,i}}(t)|^2_{R_1}\bigg) dt \\ &\qquad  +\bigg|{\color{black}y_{1,i}}(T)   -\bar{F}_1\frac{\sum_{j=1,j\neq i}^{N}{\color{black}y_{1,j}}(T)}{N-1} -\bar{G}_1y_{\alpha}(T) -\bar{H}_1y_{\beta}(T)-\bar{M}_1\bigg|_{\bar{Q}_1}^2\bigg],
\end{align*}
where 
$|\cdot|_{\mathcal{Q}}:=\sqrt{\langle \cdot, \mathcal{Q} \cdot\rangle}$ for any positive definite matrix $\mathcal{Q}$ and $\langle \cdot, \star\rangle$ is the usual Euclidean inner product. The matrices $F_{\alpha}, G_{\alpha}, M_{\alpha}, Q_{\alpha},R_\alpha, F_\beta,G_\beta$, $M_{\beta}, Q_\beta,R_\beta,F_1,G_1,H_1,Q_1,R_1$, $\bar{F}_\alpha, \bar{G}_\alpha, \bar{M}_{\alpha}, \bar{Q}_\alpha, \bar{F}_\beta, \bar{G}_\beta$, $\bar{M}_{\beta}, \bar{Q}_\beta, \bar{F}_1, \bar{G}_1, \bar{H}_1, \bar{Q}_1$ are assumed to be constant with appropriate dimensions. In addition, $Q_\alpha, Q_\beta, Q_1, R_\alpha$, $R_\beta$ and $R_1$ are positive definite, while $\bar{Q}_\alpha, \bar{Q}_\beta$ and $\bar{Q}_1$ are positive semi-definite.

Solving the empirical SDGs, either in the sense of a Nash game or a Pareto game to be described in Definitions \ref{NG} and \ref{PG} below, is rather complicated as $N$ becomes very large. Instead of directly solving the empirical game, we consider the limiting case{\color{black}s} when the number of followers $N\to\infty$ and solve the resulting mean field game{\color{black}s}. This leads to the following mean field dynamics for the $\gamma$-leader, $\gamma=\alpha,\beta$, and the $i$-th follower, respectively:  
\begin{align}
&\begin{cases}\label{eq:dyn:leader}
dx_{\gamma}(t)&=\Big(A_{\gamma}x_{\alpha}(t)+B_{\gamma}x_{\beta}(t)+C_{\gamma}z(t) +D_{\gamma}v_{\gamma}(t)\Big)dt+\sigma_{\gamma}dW^{\gamma}(t);\\
x_{\gamma}(0)&=\eta_{\gamma},\\
\end{cases}\\
\label{X1SDE}&\begin{cases}
d{\color{black}x_{1,i}}(t)&=\Big(A_1{\color{black}x_{1,i}}(t)+B_1z(t)+C_1x_{\alpha}(t)    +D_1x_{\beta}(t)+E_1{\color{black}v_{1,i}}(t)\Big)dt+\sigma_1d{\color{black}W^{1,i}}(t);\\
{\color{black}x_{1,i}}(0)&={\color{black}\eta_{1,i}}, \ {\color{black}i=1,2\dots.}
\end{cases}
\end{align}
Their cost functionals, $J_\gamma=J_\gamma(v_\gamma;v_{\alpha+\beta-\gamma})$ and ${\color{black}J_{1,i}}={\color{black}J_{1,i}}({\color{black}v_{1,i}};x_\alpha,x_\beta,z)$, are respectively given by  
    \begin{align}
        J_{\gamma} =&\ \mathbb{E}\bigg[\int_0^T \big(\big|x_{\gamma}(t)-F_{\gamma}z(t)-G_{\gamma}x_{\alpha+\beta-\gamma}(t)  -M_\gamma   \big|_{Q_{\gamma}}^2  +|v_{\gamma}(t)|_{R_{\gamma}}^2\big)dt \nonumber \\ 
        &\qquad +\big|x_{\gamma}(T)
         -\bar{F}_{\gamma}z(T)  -\bar{G}_{\gamma}x_{\alpha+\beta-\gamma}(T)-\bar{M}_{\gamma}\big|_{\bar{Q}_{\gamma}}^2\bigg];\\
    {\color{black}J_{1,i}} =&\  \mathbb{E}\bigg[\int_0^T\big(\big|{\color{black}x_{1,i}}(t)-F_1z(t)-G_1x_{\alpha}(t)-H_1x_{\beta}(t)  - M_1\big|_{Q_1}^2+|{\color{black}v_{1,i}}(t)|_{R_1}^2\big)dt \nonumber \\ &\ +\big|{\color{black}x_{1,i}}(T)-\bar{F}_1z(T) -\bar{G}_1x_{\alpha}(T) -\bar{H}_1x_{\beta}(T)-\bar{M}_1\big|_{\bar{Q}_1}^2\bigg]. \label{X1CF}
 \end{align}
{\color{black}Here, $z$ is {\color{black} an {\color{black}$\mathbb{R}^{n_1}$}-valued process adapted to $\mathcal{F}^{\alpha}\vee\mathcal{F}^{\beta}$. Under the mean field equilibrium, $z$} 
represents the average behavior of followers under the influence of two leaders:
\begin{align}\label{meanfieldfixpoint}
	z(t)=\mathbb{E}^{\mathcal{F}^{\alpha}_t\vee\mathcal{F}^{\beta}_t}[x_{1,i}(t)],\ 
    {\color{black} i=1,2,\dots},
	\end{align}
{\color{black}  where $\mathbb{E}^{\mathcal{F}^{\alpha}_t\vee\mathcal{F}^{\beta}_t}[\cdot]$ denotes the conditional expectation given the $\sigma$-algebra $\mathcal{F}^{\alpha}_t\vee\mathcal{F}^{\beta}_t$.} We shall refer to $z$ as the mean field term.} {\color{black}Unlike the case without leaders, {\color{black}$z$ herein is a stochastic process that is a functional of information generating the states of leaders}. As a result, the mean field term is exogenous {\color{black} to} followers, but endogenous {\color{black} to} leaders. {\color{black} We thus omit writing the argument dependence of $z$ in $J_\alpha$ and $J_\beta$, as it is more of an implicit consequence of the leaders' decision than an input.} This also highlights the difference with the notion of major players in  \cite{Buckdahn:2014,carmona:major,LASRY2018886,csen2014mean} in which leaders play a more influential role and followers have to comply with the actions taken by the leaders.} Comparing with the single-leader game in \cite{BCY}, the mean field term herein is jointly determined by the two leaders, {\color{black}indeed} it is adapted to $\mathcal{F}^\alpha\vee\mathcal{F}^\beta$, and neither of the leaders alone {\color{black}is supposed to} dominate the community {\color{black}\textit{a priori}}.%

Note that given $x_\alpha$, $x_\beta$ and $z$, $\{x_{1,i}\}_{i=1}^\infty$ {\color{black}share identical probabilistic behaviour}. We therefore drop the index $i$ in (\ref{X1SDE}) and (\ref{X1CF}),  and simply call $x_1$ the (representative) follower. The corresponding Wiener process, filtration and {\color{black}cost functional} will be denoted by $W^1$, $\mathcal{F}^1$ and {\color{black} $J_1$}, respectively. Depending on the mode of interactions between leaders, either competitive or cooperative, we respectively formulate the leader game in terms of a Nash and a Pareto game as follows.  In the sequel, we shall use the superscripts $\mathcal{N}$ and $\mathcal{P}$ to denote the objects corresponding to the Nash and the Pareto game, respectively. 



\begin{definition}[Nash Game]\label{NG}
	The optimal control $u_{\alpha}^{\mathcal{N}}$ and $u_{\beta}^{\mathcal{N}}$ for the non-cooperative Nash game between the two leaders are defined 
    {\color{black} by the following fixed point relationship:}
    \begin{equation}
    \left\{\begin{aligned}
	u_{\alpha}^{\mathcal{N}}&:=\argmin_{v_{\alpha}}J_{\alpha}(v_{\alpha};u_{\beta}^{\mathcal{N}});\\
    u_{\beta}^{\mathcal{N}}&:=\argmin_{v_{\beta}}J_{\beta}({\color{black}v_{\beta};u_{\alpha}^{\mathcal{N}}}).
	\end{aligned}\right.
	\end{equation}
\end{definition}

\begin{definition}[Pareto Game]\label{PG}
	The optimal control $u_{\alpha}^{\mathcal{P}}$ and $u_{\beta}^{\mathcal{P}}$ for the cooperative Pareto game between the two leaders are defined as follows:
	\begin{equation}
(u_{\alpha}^{\mathcal{P}},u_{\beta}^{\mathcal{P}}):=\argmin_{v_{\alpha},v_{\beta}}\Big(J_{\alpha}(v_{\alpha};v_{\beta})+J_{\beta}({\color{black}v_{\beta};v_{\alpha}})\Big).
	\end{equation}
\end{definition} %
{\color{black}In particular, the two leaders in the Pareto game can be considered as a single leader with a higher-dimensional state 
{\color{black}from} the follower's perspective. This can be generalized to the case with more than two leaders; see  Appendix \ref{sec:app:pareto:one}.} We now introduce the complete leader-follower MFG  by {\color{black} the following nested sub-problems, which shall be solved in a sequential manner}:


\begin{problem}[{\color{black}Follower-Game}]\label{P1}
	Given $x_\alpha$, $x_\beta$ and $z$, find a control $u_1$ such that
	\begin{align*}
	u_1&=\argmin_{v_1}J_1(v_1;x_\alpha,x_\beta,z).
	\end{align*}
\end{problem}
\begin{problem}[{\color{black}Mean Field Equilibrium of Followers}]\label{P2}
{\color{black} Given $x_\alpha$, $x_\beta$, 	f}ind the {\color{black} $\mathcal{F}^{\alpha}\vee\mathcal{F}^{\beta}$-adapted} process $z$ such that the fixed point property {\color{black} \eqref{meanfieldfixpoint}} is satisfied.
\end{problem}
    
\begin{problem}[{\color{black}Leader-Game}]\label{P3} 
	Find the optimal control for the Nash and the Pareto game defined in Definition \ref{NG} and \ref{PG}, respectively,
	where $z$ is 
 {\color{black} now an endogenous process satisfying \eqref{meanfieldfixpoint}} in Problem \ref{P2}.
\end{problem}

{\color{black} The formulation of Problems \ref{P1}-\ref{P3} is explained below. Under the influence of leaders, Problems \ref{P1} and \ref{P2} form an MFG problem at the follower's level. 
At the leader's level in Problem \ref{P3}, decisions are made based on the relationship among leaders, i.e., competition or cooperation, and also,  the average behavior of followers under the impact of their decisions. {\color{black}The latter is described by} the mean field term $z$ solved in Problem \ref{P2}, 
which creates an endogenous feedback mechanism  such that leaders take the influence of their own actions {\color{black} on the community} into consideration, keeping the essence of a Stackelberg game.} 
{\color{black} In the sequel, we shall first assume Problems \ref{P1}-\ref{P3} are solvable and leave the discussion of their well-posedness in Sections \ref{sec:sol:follower} and \ref{sec:well-posed:FBSDE}.}

{\color{black}
\begin{remark}
    The model introduced herein could be potentially generalized to controlled diffusion coefficients by employing the lifting approach and continuation method in \cite{bensoussan2023degenerate}. To maintain the primary focus on the two regimes of governance and their impact on the follower, we do not {\color{black}pursue} 
    this extension herein, which shall be left to future research endeavors. 
\end{remark}}


{\color{black} 
\subsection{$\epsilon$-Nash Equilibrium of the Follower-Game} }

    {\color{black}This subsection is devoted to establishing the $\epsilon$-Nash equilibrium of the follower-game given the dynamics of leaders. 
    } {\color{black} Let ${\color{black}{\bf u}_1:=(u_{1,1},\dots,u_{1,N})}$, where $u_{1,i}$ is the optimal control achieved in Problem \ref{P1}, i.e.,  $u_{1,i}:=\argmin_{v_{1,i}}J_{1,i}(v_{1,i};x_\alpha,x_\beta,z)$, and $z$ satisfies the fixed point property \eqref{meanfieldfixpoint}. We shall denote by ${\color{black}x_{1,i}}={\color{black}x^{{\color{black}u_{1,i}}}_{1,i}}$ and ${\color{black}\hat{y}_{1,i}}={\color{black}y^{{\color{black}u_{1,i}}}_{1,i}}$ the corresponding  MF state dynamics \eqref{X1SDE} and empirical state dynamics \eqref{y1SDE} of the $i$-th follower under ${\color{black}u_{1,i}}$. Let also $\hat{y}_\alpha$ and $\hat{y}_\beta$ be the empirical state dynamics of the two leaders given in \eqref{eq:leader:state:empirical} under the followers' dynamics $\hat{{\bf y}}_1:=(\hat{y}_{1,1},\dots,\hat{y}_{i,N})$, {\color{black} and $x_\alpha$, $x_\beta$ be the MF state dynamics of the leaders.}} Theorem \ref{equilibrium} below shows that the optimal control ${\bf u}_1$ of the {\color{black} MFG} for followers {\color{black}achieving} an $\epsilon$-Nash equilibrium\footnote{We remark that the notion of $\epsilon$-Nash equilibrium among the group of followers introduced in Theorem \ref{equilibrium} should not be confused with the Nash (or Pareto) game in Definition \ref{NG} (or \ref{PG}) between leaders.  } for the original empirical problem, which can be proven by following the same argument in Section 2 of \cite{BCY2} and the proof of Theorem 4.1 of \cite{BSYY}, {\color{black}and the proof is omitted. }

\begin{theorem}\label{equilibrium}
    {\color{black} Let ${\bf u}_1$ be the solution of Problem \ref{P1} and $z$ be the solution of Problem \ref{P2} defined above.} Then, 
 \begin{equation}
	\mathbb{E}\Big[\sup_{0\leq t\leq T}|{\color{black}\hat{y}_{\alpha}}-x_{\alpha}|^2(t)+\sup_{0\leq t\leq T}|{\color{black}\hat{y}_{\beta}}-x_{\beta}|^2(t)+\sup_{i=1,2\dots,N}\sup_{0\leq t\leq T}|{\color{black}\hat{y}_{1,i}}-{\color{black}x_{1,i}}|^2(t)\Big]=O\bigg(\frac{1}{N}\bigg).
	\end{equation} 
	Moreover, the control ${\bf u}_1$ achieves an $\epsilon$-Nash equilibrium of order $O\left(\frac{1}{\sqrt{N}}\right)$ among the followers for the original empirical problem, that is, {\color{black} for $i=1,\dots,N$ and} any arbitrary control ${\color{black}v_{1,i}}$, we have \vspace{-0.6cm}%
 {\color{black}
    \begin{equation}
        \mathcal{J}_{1,i}(u_{1,i};{\bf u}_{1,-i}) \le \mathcal{J}_{1,i}(v_{1,i};{\bf u}_{1,-i})+O\bigg(\frac{1}{\sqrt{N}}\bigg),
    \end{equation}}%
   where ${\bf u}_{1,-i} := (u_{1,1},\dots,u_{1,i-1},u_{1,i+1},\dots,u_{1,N})$. 
\end{theorem}%

\section{Constructions of Solutions}\label{sec:sol}

For any symmetric matrix $M_1$, we denote by $\lambda(M_1)$ its {\color{black}generic} eigenvalues with $\lambda_{\min}(M)$ being the smallest one, and $\lambda_{\max}(M)$ being the largest one. For any matrix $M_2$, we denote by $\rho(M_2):=\sqrt{\lambda_{\max}(M_2^{\top}M_2)}$ the spectral norm of $M_2$. We also write $I$ as the identity matrix of compatible dimensions.  Let $C([0,T]; \mathbb{R}^{n\times k})$ be the space of $\mathbb{R}^{n\times k}$-valued continuous functions equipped with the norm $\rho_T({\bf \Phi}) := \sup_{0\leq t\leq T}\rho({\bf \Phi}(t))$, for any $\Phi\in C([0,T]; \mathbb{R}^{n\times k})$.

\subsection{Optimal Control for the Follower}
\label{sec:sol:follower}
Theorem \ref{C3_follower} below provides the optimal control of the follower given the dynamics of the leaders $x_\alpha$, $x_\beta$ and the mean field term $z$.
\begin{theorem}\label{C3_follower}
	Given $x_\alpha$, $x_\beta$ and $z$, the following FBSDE is uniquely solvable  
   \begin{equation}\label{x1andn}
\left\{\begin{aligned}
dx_1(t)&=\Big(A_1x_1(t)+B_1z(t)+C_1x_{\alpha}(t)+D_1x_{\beta}(t)-E_1R_1^{-1}E_1^{\top}\xi(t)\Big)dt+\sigma_1dW^1(t);\\
-d\xi(t)&=\Big(A_1^{\top}\xi(t)+Q_1\Big(x_1(t)-F_1z(t)-G_1x_{\alpha}(t)  -H_1x_{\beta}(t)-M_1\Big)\Big)dt\\ &\qquad-Z_{\xi,\alpha}(t)dW^{\alpha}(t) -Z_{\xi,\beta}(t)dW^{\beta}(t)-Z_{\xi,1}(t)dW^1(t);\\	
x_1(0)&=\eta_1, \ \xi(T)=\bar{Q}_1\Big(x_1(T)-\bar{F}_1z(T)   -\bar{G}_1x_{\alpha}(T)
 -\bar{H}_1x_{\beta}(T)-\bar{M}_1\Big).
	\end{aligned}\right.
\end{equation}
    {\color{black}Moreover,} the optimal control for Problem \ref{P1} is given by $u_1(t) = -R_1^{-1}E_1^{\top} \xi(t)$.
\end{theorem}

{\color{black}Given $x_\alpha, x_\beta$ and $z$, the FBSDE \eqref{x1andn} and thus the optimal control can be deduced by following the lines of argument of Lemma 3.4 of \cite{BCY}. In addition, it is easy to verify that \eqref{x1andn} satisfies Assumptions 2.1 and {\color{black}3.1} 
of \cite{Peng}, which are respectively the standard Lipschitz and monotonicity conditions. The unique existence of the \eqref{x1andn} is then an immediate consequence of Theorem {\color{black}3.2} 
of \cite{Peng}.  }

To obtain the mean field equilibrium stated in Problem \ref{P2}, we take conditional expectation {\color{black}on} $\mathcal{F}_t^{\alpha}\vee \mathcal{F}_t^{\beta}$ {\color{black}of} both sides of \eqref{x1andn}, which yields the following FBSDE:  
\begin{equation} \label{ZM}
\left\{\begin{aligned}
dz(t)&=\Big((A_1+B_1)z(t)+C_1x_{\alpha}(t)+D_1x_{\beta}(t) -E_1R_1^{-1}E_1^{\top} \zeta(t)\Big)dt;\\
-d\zeta(t)&=\Big(A_1^{\top} \zeta(t)+Q_1(I-F_1)z(t)-Q_1G_1x_{\alpha}(t)-Q_1H_1x_{\beta}(t)-Q_1M_1\Big)dt\\
&\qquad -Z_{\zeta,\alpha}(t)dW^{\alpha}(t)-Z_{\zeta,\beta}(t)dW^{\beta}(t);\\
z(0)&=\mathbb{E}[\eta_1],\ \zeta(T)=\bar{Q}_1\Big((I-\bar{F}_1)z(T)  -\bar{G}_1x_{\alpha}(T) -\bar{H}_1x_{\beta}(T)-\bar{M}_1\Big).
\end{aligned}\right.
\end{equation}
Hence, given $x_\alpha$ and $x_\beta$, Problem \ref{P2} is solvable if the FBSDE \eqref{ZM} admits a unique solution; {\color{black}t}o this end, we {\color{black}provide with} following {\color{black} sufficient} condition:
\begin{assumption}\label{fixedptcond}
    Assume that the following hold:  
\begin{equation}\label{fixedptcond1}
    \begin{aligned}
  \rho^2(B_1) <&\  4 \lambda_{\min}\left(\frac{Q_1(I-F_1)+(I-F_1)^{\top}Q_1}2\right) \lambda_{\min}(E_1R_1^{-1}E_1^{\top}),
  \end{aligned}
\end{equation}
and $\lambda_{\min}\left( \bar{Q}_1(I-\bar{F}_1)+(I-\bar{F}_1)^{\top}\bar{Q}_1\right)\geq 0$.
\end{assumption}

In particular, Assumption \ref{fixedptcond} holds when the mean field effect on the follower is not too large, i.e., when $B_1$, $F_1$ and $\bar{F}_1$ have a small magnitude.  The unique existence of solution of \eqref{ZM} is asserted by the following theorem, whose proof is relegated to  Appendix \ref{app:bsde}.
\begin{theorem}\label{ZMEX}
 Under Assumption \ref{fixedptcond}, the system \eqref{ZM} admits a unique solution for any given $x_\alpha$ and $x_\beta$.
\end{theorem}



\subsection{Optimal Control for the Leaders: Nash Game}

{\color{black}To introduce the optimal control of the leaders under the Nash game, we define the following matrices:}
 \begin{equation}\label{bigmatrixA}
 {\bf  A}:=
\begin{pmatrix}
 A_{\alpha} & B_{\alpha} & C_{\alpha} & 0 & 0\\
 A_{\beta} & B_{\beta} & C_{\beta} & 0 & 0\\
 C_1 & D_1 & A_1+B_1 & 0 & 0\\
 0 & 0 & 0 & A_1 & 0\\
 0 & 0 & 0 & 0 & A_1
\end{pmatrix} , 
 \end{equation}
\begin{equation}\label{bigmatrixB}
\begin{aligned}
& {\bf  B}{\color{black}^\mathcal{N}}:= 
\left(\begin{matrix}
 D_{\alpha}R_{\alpha}^{-1}D_{\alpha}^{\top} & 0 & 0 & 0 & 0 & {\color{black} 0} & {\color{black} 0} \\
 0 & D_{\beta}R_{\beta}^{-1}D_{\beta}^{\top} & 0 & 0 & 0 & {\color{black} 0} & {\color{black} 0} \\
 0 & 0 & E_1R_1^{-1}E_1^{\top} & 0 & 0 & {\color{black} 0} & {\color{black} 0} \\
 0 & 0 & 0 & E_1R_1^{-1}E_1^{\top} & 0 & {\color{black} 0} & {\color{black} 0} \\
 0 & 0 & 0 & 0 & E_1R_1^{-1}E_1^{\top} & {\color{black} 0} & {\color{black} 0} \\
\end{matrix}\right), 
\end{aligned}
 \end{equation}
\begin{equation*}
 {\bf  C}^{\mathcal{N}}:=
\left(\begin{matrix}
 A_{\alpha}^{\top} & 0 & 0 & C_1^{\top} & 0 & {\color{black} A_\beta^\top} & {\color{black} 0} \\
 0 & B_{\beta}^{\top} & 0 & 0 & D_1^{\top} & {\color{black} 0} & {\color{black} B_\alpha^\top}\\
 0 & 0 & A_1^{\top} & 0 & 0 & {\color{black} 0} & {\color{black} 0} \\
 C_{\alpha}^\top & 0 & 0 & (A_1+B_1)^{\top} & 0 & {\color{black} C_\beta^\top} & {\color{black} 0}\\
 0 & C_{\beta}^\top & 0 & 0 & (A_1+B_1)^{\top} & {\color{black} 0} & {\color{black} C_\alpha^\top}\\
 {\color{black} B_\alpha^\top} & {\color{black} 0} & {\color{black} 0} & {\color{black} D_1^\top} & {\color{black} 0} & {\color{black} B_\beta^\top} & {\color{black} 0} \\ 
 {\color{black} 0} & {\color{black} A_\beta^\top} & {\color{black} 0} & {\color{black} 0} & {\color{black} C_1^\top} & {\color{black} 0} & {\color{black} A_\alpha^\top}   \\ 
\end{matrix}\right),
 \end{equation*}
 \begin{equation}
\label{eq:D:nash}
 \begin{aligned}
&{\bf  D}^{\mathcal{N}}:= 
\left(\begin{matrix}
 Q_{\alpha} & -Q_{\alpha}G_{\alpha} & -Q_{\alpha}F_{\alpha} & -(Q_1G_1)^{\top} & 0 
 \\
-Q_{\beta}G_{\beta} & Q_{\beta} & -Q_{\beta}F_{\beta} & 0 & -(Q_1H_1)^{\top} 
\\
-Q_1G_1 & -Q_1H_1 & Q_1(I-F_1) & 0 &0 
\\
-(Q_{\alpha}F_{\alpha})^{\top} & (Q_{\alpha}F_{\alpha})^{\top}G_{\alpha} & (Q_{\alpha}F_{\alpha})^{\top}F_{\alpha} & (Q_1(I-F_1))^{\top} & 0 
\\
(Q_{\beta}F_{\beta})^{\top}G_{\beta} & -(Q_{\beta}F_{\beta})^{\top} & (Q_{\beta}F_{\beta})^{\top}F_{\beta} & 0 & (Q_1(I-F_1))^{\top} 
\\
{\color{black}-(Q_\alpha G_\alpha)^\top} & {\color{black}(Q_\alpha G_\alpha)^\top G_\alpha} & {\color{black}(Q_\alpha G_\alpha)^\top F_\alpha }& {\color{black}-(Q_1H_1)^\top} & {\color{black}0} 
\\  
{\color{black}(Q_\alpha G_\beta)^\top G_\beta} & {\color{black}-(Q_\alpha G_\beta)^\top  } & {\color{black}(Q_\alpha G_\beta)^\top F_\beta }& {\color{black}0} & {\color{black}-(Q_1G_1)^\top}   
\end{matrix}\right),
\end{aligned}
\end{equation}
 \begin{equation}\label{bigmatrixVol}
{\bf  \Sigma}:=
\left(\begin{matrix}
 \sigma_{\alpha} & 0\\
 0 & \sigma_{\beta}\\
 0 & 0\\
 0 & 0\\
 0 & 0\\
\end{matrix}\right),
\ {\bf  W}(t):=
\left(\begin{matrix}
W^{\alpha}(t)\\W^{\beta}(t)
\end{matrix}\right),
\ {{\bf  x}_0}:=
\left(\begin{matrix}
\eta_{\alpha}\\
\eta_{\beta}\\
\mathbb{E}[\eta_1]\\
0\\
0\\
\end{matrix}\right),
\end{equation}
 \begin{equation}
\begin{aligned}
&{\bf  E}^{\mathcal{N}}:= 
\left(\begin{matrix}
 \bar{Q}_{\alpha} & -\bar{Q}_{\alpha}\bar{G}_{\alpha} & -\bar{Q}_{\alpha}\bar{F}_{\alpha} & -(\bar{Q}_1\bar{G}_1)^{\top} & 0 
 \\
-\bar{Q}_{\beta}\bar{G}_{\beta} & \bar{Q}_{\beta} & -\bar{Q}_{\beta}\bar{F}_{\beta} & 0 & -(\bar{Q}_1\bar{H}_1)^{\top} 
\\
-\bar{Q}_1\bar{G}_1 & -\bar{Q}_1\bar{H}_1 & \bar{Q}_1(I-\bar{F}_1) & 0 &0 
\\
-(\bar{Q}_{\alpha}\bar{F}_{\alpha})^{\top} & (\bar{Q}_{\alpha}\bar{F}_{\alpha})^{\top}\bar{G}_{\alpha} & (\bar{Q}_{\alpha}\bar{F}_{\alpha})^{\top}\bar{F}_{\alpha} & (\bar{Q}_1(I-\bar{F}_1))^{\top} & 0 
\\
(\bar{Q}_{\beta}\bar{F}_{\beta})^{\top}\bar{G}_{\beta} & -(\bar{Q}_{\beta}\bar{F}_{\beta})^{\top} & (\bar{Q}_{\beta}\bar{F}_{\beta})^{\top}\bar{F}_{\beta} & 0 & (\bar{Q}_1(I-\bar{F}_1))^{\top} 
\\
{\color{black}-(\bar{Q}_\alpha \bar{G}_\alpha)^\top} & {\color{black}(\bar{Q}_\alpha \bar{G}_\alpha)^\top \bar{G}_\alpha} & {\color{black}(\bar{Q}_\alpha \bar{G}_\alpha)^\top \bar{F}_\alpha }& {\color{black}-(\bar{Q}_1\bar{H}_1)^\top} & {\color{black}0} 
\\  
{\color{black}(\bar{Q}_\alpha \bar{G}_\beta)^\top \bar{G}_\beta} & {\color{black}-(\bar{Q}_\alpha \bar{G}_\beta)^\top  } & {\color{black}(\bar{Q}_\alpha \bar{G}_\beta)^\top \bar{F}_\beta }& {\color{black}0} & {\color{black}-(\bar{Q}_1\bar{G}_1)^\top} 
\end{matrix}\right),
\end{aligned}
\label{eq:E:nash}
\end{equation}
and
 \begin{equation*}
 {\bf  M}^{\mathcal{N}}:=
\left(\begin{matrix}
-Q_{\alpha}M_{\alpha}\\
-Q_{\beta}M_{\beta}\\
-Q_1M_1\\
(Q_{\alpha}F_{\alpha})^{\top}M_{\alpha}\\	(Q_{\beta}F_{\beta})^{\top}M_{\beta}\\
{\color{black}(Q_\alpha G_\alpha)^\top M_\alpha} \\
{\color{black}(Q_\beta G_\beta)^\top M_\beta}
\end{matrix}\right),\ {\bf  h}^{\mathcal{N}}:=
\left(\begin{matrix}
-\bar{Q}_{\alpha}\bar{M}_{\alpha}\\
-\bar{Q}_{\beta}\bar{M}_{\beta}\\
-\bar{Q}_1\bar{M}_1\\
(\bar{Q}_{\alpha}\bar{F}_{\alpha})^{\top}\bar{M}_{\alpha}\\	(\bar{Q}_{\beta}\bar{F}_{\beta})^{\top}\bar{M}_{\beta}\\
{\color{black}(\bar{Q}_\alpha \bar{G}_\alpha)^\top \bar{M}_\alpha} \\
{\color{black}(\bar{Q}_\beta \bar{G}_\beta)^\top \bar{M}_\beta}
\end{matrix}\right).
 \end{equation*}
 {\color{black} The singularity of the matrix ${\bf B}^\mathcal{N}$ shall be addressed in Section \ref{sec:well-posed:FBSDE}, see also Assumption \ref{Fullcond2:Nash} therein, whereas the singularity of ${\bf \Sigma}$ is not relevant in the subsequent discussion.}  By a standard first-order method, the solution of the Nash game introduced in Definition \ref{NG} is given by the following theorem, whose proof is relegated to  Appendix \ref{app:bsde}. 

\begin{theorem}\label{thm_Nash}
{\color{black}Suppose that the following FBSDE admits a unique solution:}
 \begin{equation}\label{NumericalNG}
    \left\{\begin{aligned}
    d{\bf  x}^{\mathcal{N}}(t)&=\Big({\bf  A}{\bf  x}^{\mathcal{N}}(t)-{\color{black}{\bf  B}^\mathcal{N}}{\bf  p}^{\mathcal{N}}(t)\Big)dt+{\bf  \Sigma} d{\bf  W}(t);\\
    -d{\bf  p}^{\mathcal{N}}(t)&=\Big({\bf  C}^{\mathcal{N}}{\bf  p}^{\mathcal{N}}(t)+{\bf  D}^{\mathcal{N}}{\bf  x}^{\mathcal{N}}(t)+{\bf  M}^{\mathcal{N}}\Big)dt -{\bf  Z}^{\mathcal{N}}(t)d{\bf  W}(t);\\
        {\bf  x}^{\mathcal{N}}(0)&={\bf  x}_0
,\quad {\bf  p}^{\mathcal{N}}(T)={\bf  E}^{\mathcal{N}}{\bf  x}^{\mathcal{N}}(T)+{\bf  h}^{\mathcal{N}}.
    \end{aligned}\right.
    \end{equation} 
The solution components are given by   ${\bf x}^\mathcal{N} := (x_{\alpha}^{\mathcal{N}}, x_{\beta}^{\mathcal{N}},$  $z^{\mathcal{N}}, s^{\mathcal{N}}_{\alpha}, s^{\mathcal{N}}_{\beta})$,  ${\bf p}^\mathcal{N}:=(p_{\alpha}^{\mathcal{N}},  p_{\beta}^{\mathcal{N}}, \zeta^{\mathcal{N}}$, $r^{\mathcal{N}}_{\alpha} , r^{\mathcal{N}}_{\beta}, {\color{black} q_\alpha^{\mathcal{N}},q_\beta^{\mathcal{N}}}$). 
{\color{black}Then, the optimal controls of the Nash game in Definition \ref{NG} are given by, for $0\leq t\leq T$,}
	\begin{equation}\label{Nashoptimalcontrol}
	(u_{\alpha}^{\mathcal{N}}(t), u_\beta^\mathcal{N}(t))=(-R_{\alpha}^{-1}D_{\alpha}^{\top}p_{\alpha}^{\mathcal{N}}(t), -R_{\beta}^{-1}D_{\beta}^{\top}p_{\beta}^{\mathcal{N}}(t)).
	\end{equation}
\end{theorem}%
{\color{black} The solution of the leader-game is described by an augmented system of FBSDE \eqref{NumericalNG}. In contrast to a major-minor player game, the additional variables $s_\alpha^\mathcal{N}$ and $s_\alpha^\mathcal{N}$ describe the variation {\color{black} on the mean field term as a result of respective leaders' action}
, and $r_\alpha^\mathcal{N}$ and $r_\beta^\mathcal{N}$ are their adjoint processes, respectively.} {\color{black} The process $q_\alpha^{\mathcal{N}}$ (resp.~$q_\beta^{\mathcal{N}}$) describes the variations of the $\alpha$-leader's (resp.~$\beta$-leader's) control due to its impact on the counterparty $x_\beta^{\mathcal{N}}$ (resp.~$x_\alpha^{\mathcal{N}}$) when he makes his own decision.}  



\subsection{Optimal Control for the Leaders: Pareto Game}
{\color{black} To introduce the optimal control for the leaders under the Pareto game, we define} {\color{black} ${\bf  B}^{\mathcal{P}}:={\bf  B}^{\mathcal{N}}{\bf I}^{\mathcal{P},\mathcal{N}}_2$, ${\bf  C}^{\mathcal{P}}:={\bf I}^{\mathcal{P},\mathcal{N}}_1{\bf  C}^{\mathcal{N}}{\bf I}^{\mathcal{P},\mathcal{N}}_2$, }
 {\color{black}${\bf D}^{\mathcal{P}} := {\bf I}^{\mathcal{P},\mathcal{N}}_1{\bf D}^{\mathcal{N}}  $,  ${\bf E}^{\mathcal{P}} := {\bf I}^{\mathcal{P},\mathcal{N}}_1{\bf E}^{\mathcal{N}}  $,  ${\bf M}^{\mathcal{P}} := {\bf I}^{\mathcal{P},\mathcal{N}}_1{\bf M}^{\mathcal{N}}  $ and  ${\bf h}^{\mathcal{P}} := {\bf I}^{\mathcal{P},\mathcal{N}}_1{\bf h}^{\mathcal{N}}  $, where
\begin{equation*}
       {\bf I}^{\mathcal{P},\mathcal{N}}_1 := 
      \left(  \begin{matrix}
            I & 0 & 0 &0 & 0 & 0 & I \\
            0 & I & 0 &0 & 0 & I & 0 \\
            0 & 0 & I &0 & 0 & 0 & 0 \\
            0 & 0 & 0 &I & 0 & 0 & 0 \\
         0 & 0 & 0 &0 & I & 0 & 0 \\
      \end{matrix}\right),\
      {\color{black}{\bf I}^{\mathcal{P},\mathcal{N}}_{2} := 
      \left(  \begin{matrix}
            I & 0 & 0 &0 & 0  \\
            0 & I & 0 &0 & 0  \\
            0 & 0 & I &0 & 0  \\
            0 & 0 & 0 &I & 0  \\
         0 & 0 & 0 &0 & I  \\
         0 & 0 & 0 &0 & 0  \\
         0 & 0 & 0 &0 & 0  \\
      \end{matrix}\right).}
\end{equation*}
In particular, the matrix ${\bf I}_2^{\mathcal{P},\mathcal{N}}$ removes the last two columns of ${\bf B}^\mathcal{N}$ that renders a non-singular matrix ${\bf B}^\mathcal{P}$. 

}%

{\color{black}Theorem \ref{thm_Pareto} below provides the  solution of the Pareto game between leaders, which can again be deduced by the first-order method; see  Appendix  \ref{app:bsde} for a detailed proof. }

 
\begin{theorem}\label{thm_Pareto}
{\color{black}Suppose that the following FBSDE admits a unique solution:}
\begin{equation}\label{NumericalPG}
\left\{\begin{aligned}
d{\bf  x}^{\mathcal{P}}(t)&=\Big({\bf  A}{\bf  x}^{\mathcal{P}}(t)-{\bf  B}{\color{black}^\mathcal{P}}{\bf  p}^{\mathcal{P}}(t)\Big)dt+{\bf  \Sigma} d{\bf  W}(t);\\
-d{\bf  p}^{\mathcal{P}}(t)&=\Big({\bf  C}^{\mathcal{P}}{\bf  p}^{\mathcal{P}}(t)+{\bf  D}^{\mathcal{P}}{\bf  x}^{\mathcal{P}}(t)+{\bf  M^{\mathcal{P}}}\Big)dt  -{\bf  Z}^{\mathcal{P}}(t)d{\bf  W}(t);\\
{\bf  x}^{\mathcal{P}}(0)&={\bf  x}_0,\quad {\bf  p}^{\mathcal{P}}(T)={\bf  E}^{\mathcal{P}}{\bf  x}^{\mathcal{P}}(T)+{\bf  h}^{\mathcal{P}},
\end{aligned}\right.
\end{equation}
 The solution components are given by ${\bf x}^\mathcal{P} := (x_{\alpha}^{\mathcal{P}}, x_{\beta}^{\mathcal{P}}, z^{\mathcal{P}},$ $s^{\mathcal{P}}_{\alpha}, s^{\mathcal{P}}_{\beta})$, ${\bf p}^\mathcal{P}:=(p_{\alpha}^{\mathcal{P}}, p_{\beta}^{\mathcal{P}} , \zeta^{\mathcal{P}} ,r^{\mathcal{P}}_{\alpha} , r^{\mathcal{P}}_{\beta})$. 
{\color{black}Then,	the pair of optimal controls of the Pareto game in Definition \ref{PG} is given by, for $0\leq t\leq T$,}
\begin{equation}
	\label{Paretooptimalcontrol}
(u_{\alpha}^{\mathcal{P}}(t),u_{\beta}^{\mathcal{P}}(t))=(-R_{\alpha}^{-1}D_{\alpha}^{\top}p^{\mathcal{P}}_{\alpha}(t),-R_{\beta}^{-1}D_{\beta}^{\top}p^{\mathcal{P}}_{\beta}(t)).
\end{equation}
\end{theorem}


{\color{black} Unlike the Nash game, the impact of the variation of a leader's state process on the control of the opposite leader has been incorporated in the processes $p_\alpha^\mathcal{P}$ and $p_\beta^\mathcal{P}$  under the cooperative environment. Hence, the auxiliary processes $q_\alpha$ and $q_\beta$ for the Nash game can be dropped herein, and the dimension of System \eqref{NumericalPG} is slightly smaller than that of \eqref{NumericalNG}. }

The well-posedness of System \eqref{NumericalPG}, and thus that of the Pareto game between the two leaders, will be discussed in Theorem \ref{wellpose:dominate} in the next subsection. 


\subsection{The Well-posedness of FBSDEs \eqref{NumericalNG} and \eqref{NumericalPG}}
\label{sec:well-posed:FBSDE}

This subsection is devoted to provide sufficient conditions for the {\color{black}unique existence of equilibrium controls for both the Nash and Pareto game, which are respectively characterized by} Systems   \eqref{NumericalNG} and \eqref{NumericalPG}. In the sequel, we shall write ${\bf M}^\mathcal{I}$, $\mathcal{I}=\mathcal{N}$ or $\mathcal{P}$, to denote the associated matrix ${\bf M}$ for the Nash game or the Pareto game, respectively. {\color{black} We also define $\mathbbm{1}_{\mathcal{I}=\mathcal{N}}$ by  $\mathbbm{1}_{\mathcal{I}=\mathcal{N}} =1$ if $\mathcal{I}=\mathcal{N}$, or zero otherwise. The same definition applies to  $\mathbbm{1}_{\mathcal{I}=\mathcal{P}}$.} 

{\color{black} For $\mathcal{I}=\mathcal{N}$ and $\mathcal{P}$, we consider an \textit{ansatz} of the adjoint process taking the affine form ${\bf p}^\mathcal{I}(t) ={\bf \Gamma}^\mathcal{I}(t) {\bf x}^\mathcal{I}(t) + {\bf g}^\mathcal{I}(t)$, for some deterministic functions ${\bf \Gamma}^\mathcal{I}$ and ${\bf g}^\mathcal{I}$. An application of It\^o's lemma yields that   ${\bf \Gamma}^\mathcal{I}$ satisfies the following non-symmetric Riccati equation:
 
	\begin{equation}\left\{
	\begin{aligned}
	-\frac{d {\bf \Gamma}^{\mathcal{I}}(t)}{dt} =&\ {\bf \Gamma}^{\mathcal{I}}(t){\bf  A}+{\bf  C}^{\mathcal{I}}{\bf \Gamma}^{\mathcal{I}}(t)-{\bf \Gamma}^{\mathcal{I}}(t){\bf  B}^{\mathcal{I}}{\bf \Gamma}^{\mathcal{I}}(t)  + {\bf  D}^{\mathcal{I}},\\
  {\bf \Gamma}^{\mathcal{I}}(T)=&\ {\bf E}^{\mathcal{I}};
	\end{aligned}\right.
 \label{dynGamma}
	\end{equation}}%
while ${\bf g}^\mathcal{I}$ satisfies the following equation: 
 
	\begin{equation}\label{NumericalAll2g}
    \left\{\begin{aligned}
      -d{\bf  g}^{\mathcal{I}}(t)=&\ \left(({\bf  C}^{\mathcal{I}}-{\bf \Gamma}^{\mathcal{I}}(t){\bf  B}^{\mathcal{I}}){\bf  g}^{\mathcal{I}}(t)+{\bf  M}^{\mathcal{I}}\right)dt,\\
    {\bf  g}^{\mathcal{I}}(T)=&\ {\bf h}^{\mathcal{I}}.
    \end{aligned}\right.
    \end{equation}
{\color{black} Since System \eqref{NumericalNG} (resp. \eqref{NumericalPG}) is linear, it admits a unique solution 
if and only if \eqref{dynGamma} does
for $\mathcal{I}=\mathcal{N}$ (resp. $\mathcal{I}=\mathcal{P}$), regardless of the singularity of the coefficient matrices such as  ${\bf B}^\mathcal{I}$ and ${\bf D}^\mathcal{I}$  (see Lemma \ref{unique:FBSDE} of  Appendix \ref{sec:app:FBSDE:riccati}). 
In the following, we introduce a generic global condition for the well-posedness of \eqref{dynGamma}; see \cite{matrix:riccati}.   

}

{\color{black} 
\begin{assume}{2A}
  \label{Fullcond2}
  There exists a symmetric matrix $\Pi^\mathcal{I}_1 \in \mathbb{R}^{{(3n_1+n_{\alpha}+n_{\beta})\times (3n_1+n_{\alpha}+n_{\beta})}}$, and a matrix   $\Pi^{\mathcal{I}}_2 \in \mathbb{R}^{{(3n_1+n_{\alpha}+n_{\beta})\times (3n_1+(1+\mathbbm{1}_{\mathcal{I}=\mathcal{N}})(n_{\alpha}+n_{\beta}))}}$, with
    \begin{equation*}
        \mathbf{\Lambda}^\mathcal{I}:= \begin{pmatrix}
            -\Pi^\mathcal{I}_1{\bf A}+\Pi^\mathcal{I}_2{\bf D}^\mathcal{I} & \Pi^\mathcal{I}_1{\bf B}^\mathcal{I}-{\bf A}^{\top}\Pi^\mathcal{I}_2+\Pi^\mathcal{I}_2{\bf C}^\mathcal{I}\\ {\bf 0} & ({\bf B}^\mathcal{I})^{\top}\Pi^\mathcal{I}_2
        \end{pmatrix},
    \end{equation*}
    such that
            \begin{equation}
            \label{eq:new:condition}
                \begin{aligned}
                    \lambda_{\min}\left( \mathbf{\Lambda}^\mathcal{I}+(\mathbf{\Lambda}^\mathcal{I})^{\top}\right) &\geq 0,\\  \text{and } \lambda_{\min}\left(\Pi^\mathcal{I}_1+\Pi^\mathcal{I}_2{\bf E}^\mathcal{I}+({\bf E}^\mathcal{I})^{\top}(\Pi^\mathcal{I}_2)^{\top} \right)&>0. 
                \end{aligned}
            \end{equation}
     
\end{assume}}%

{\color{black} Assumption \ref{Fullcond2} is practically relevant for the Pareto game. For example,  when the mean field effect and the interactions between leaders are small (in other words, a limited degree of cooperation between leaders), one natural choice of $\Pi^\mathcal{P}_1$ and $\Pi^\mathcal{P}_2$ would be the zero matrix and the identity matrix of appropriate dimensions, respectively. 
Another natural choice is to take both  $\Pi^\mathcal{P}_1$ and $\Pi^\mathcal{P}_2$ to be the identity matrix. 



For the Nash game, however, the singularity of ${\bf B}^\mathcal{N}$ and the highly non-symmetric nature of the system could make it difficult to choose appropriate matrices $\Pi_1^\mathcal{N}$ and $\Pi_2^\mathcal{N}$ such that  \eqref{eq:new:condition} holds, even for the case with a small mean field effect. 
To prepare a more easily verifiable condition, we augment the state process as $\bar{{\bf x}}^\mathcal{N} := ({\bf x}^\mathcal{N}, \hat{{\bf x}}^\mathcal{N})$, where $\hat{{\bf x}}^\mathcal{N}$ is a $\mathbb{R}^{n_\alpha+n_\beta}$-valued process with  $\hat{{\bf x}}^\mathcal{N}(0) = {\bf 0}$ 
, and 
    \begin{equation*}
        d\hat{{\bf x}}^\mathcal{N}(t) = \left(\hat{{\bf A}}^\mathcal{N}\hat{{\bf x}}^\mathcal{N}(t) - \hat{\bf B}^\mathcal{N}{\bf p}^\mathcal{N}(t)\right)dt, 
    \end{equation*}
where the matrices $\hat{{\bf A}}^\mathcal{N}$ and $\hat{{\bf B}}^\mathcal{N}$ are arbitrarily chosen matrices of compatible dimensions. By doing so, the processes 
$\bar{{\bf x} }^\mathcal{N}$ and ${\bf p}^\mathcal{N}$ share the same dimension, which satisfy the following FBSDE:
   
 \begin{equation}\label{eq:Nash:FBSDE:augment}
    \left\{\begin{aligned}
    d\bar{{\bf x}}^{\mathcal{N}}(t)&=\Big(\bar{{\bf  A}}^\mathcal{N}\bar{{\bf x}}^{\mathcal{N}}(t)-\bar{{\bf  B}}^\mathcal{N}{\bf  p}^{\mathcal{N}}(t)\Big)dt+{\bf  \Sigma} d{\bf  W}(t);\\
    -d{\bf  p}^{\mathcal{N}}(t)&=\Big({\bf  C}^{\mathcal{N}}{\bf  p}^{\mathcal{N}}(t)+\bar{{\bf  D}}^{\mathcal{N}}\bar{{\bf  x}}^{\mathcal{N}}(t)+{\bf  M}^{\mathcal{N}}\Big)dt  -{\bf  Z}^{\mathcal{N}}(t)d{\bf  W}(t);\\
    \bar{{\bf  x}}^{\mathcal{N}}(0)& =\begin{pmatrix}
            {\bf  x}_0 \\
            {\bf 0}
        \end{pmatrix}
,\quad {\bf  p}^{\mathcal{N}}(T)=\bar{{\bf  E}}^{\mathcal{N}}\bar{{\bf  x}}^{\mathcal{N}}(T)+{\bf  h}^{\mathcal{N}},
    \end{aligned}\right.
    \end{equation}
where 
    \begin{equation}\label{matrix:Nash:augment}
    \begin{aligned}
         & \bar{{\bf  A}}^\mathcal{N} = \begin{pmatrix}
            {\bf A} & {\bf 0} \\
            {\bf 0} & \hat{{\bf A}}^\mathcal{N}
        \end{pmatrix},\ \bar{{\bf B}}^\mathcal{N} =\begin{pmatrix}
            {\bf B}^\mathcal{N} \\
            \hat{{\bf B}}^\mathcal{N}
        \end{pmatrix}, \\
        & \bar{{\bf D}}^\mathcal{N}=  \begin{pmatrix}   {\bf D}^\mathcal{N} & {\bf 0}
        \end{pmatrix},\ \bar{{\bf E}}^\mathcal{N}=  \begin{pmatrix}
             {\bf E}^\mathcal{N} & {\bf 0}
        \end{pmatrix}.
    \end{aligned}
    \end{equation}
Notice that the processes ${\bf x}^\mathcal{N}$ and ${\bf p}^\mathcal{N}$ are unaffected by   $\hat{{\bf x}}^\mathcal{N}$, indeed, by the construction of the augmented matrices \eqref{matrix:Nash:augment}, the first $(3n_1+n_{\alpha}+n_{\beta})$ entries of $\bar{{\bf x}}^\mathcal{N}$, along with ${\bf p}^\mathcal{N}$, always constitute the  solution of \eqref{NumericalNG}. Hence, the well-posedness of \eqref{NumericalNG} is equivalent to that of \eqref{eq:Nash:FBSDE:augment}, whose well-posedness is shown as below. This augmentation also allows us to avoid directly considering the  singularity of ${\bf B}^\mathcal{N}$ in $\boldsymbol{\Lambda}^\mathcal{N}$ that makes Condition \eqref{eq:new:condition} difficult to be fulfilled.


Following the same argument,  ${\bf p}^\mathcal{N}$ takes an affine form ${\bf p}^\mathcal{N}(t)=\bar{\Gamma}^\mathcal{N}(t)\bar{{\bf x}}^\mathcal{N}(t) + \bar{{\bf g}}^\mathcal{N}(t)$ for some suitable function $\bar{{\bf g}}^\mathcal{N}$, and $\bar{\bf \Gamma}^\mathcal{N}$ is characterized by the augmented non-symmetric Riccati equation: 
	\begin{equation}\left\{
	\begin{aligned}
	-\frac{d \bar{{\bf \Gamma}}^{\mathcal{N}}(t)}{dt} =&\ \bar{{\bf \Gamma}}^{\mathcal{N}}(t)\bar{{\bf  A}}^\mathcal{N}+{\bf  C}^{\mathcal{N}}\bar{{\bf \Gamma}}^{\mathcal{N}}(t)-\bar{{\bf \Gamma}}^{\mathcal{N}}(t)\bar{{\bf  B}}^\mathcal{N}\bar{{\bf \Gamma}}^{\mathcal{N}}(t)   +\bar{{\bf  D}}^{\mathcal{N}}, \\
  \bar{{\bf \Gamma}}^{\mathcal{N}}(T)=&\ \bar{{\bf E}}^{\mathcal{N}}. 
	\end{aligned}\right.
 \label{dynGamma:agumented}
	\end{equation}
When \eqref{dynGamma:agumented} is well-posed, it is easy to verify that the solution is given by $\bar{{\bf \Gamma}}^\mathcal{N}(t) = ({\bf \Gamma}^\mathcal{N}(t) \ {\bf 0})$, where ${\bf \Gamma}^\mathcal{N}$ is 
a solution of \eqref{dynGamma}.  On the other hand, when \eqref{dynGamma} admits a solution ${\bf \Gamma}^\mathcal{N}(t)$, we can easily show that the matrix-valued function $\bar{{\bf \Gamma}}^\mathcal{N}(t) = ({\bf \Gamma}^\mathcal{N}(t) \ {\bf 0})$ solves \eqref{dynGamma:agumented}. We can thus establish the well-posedness of \eqref{dynGamma} for $\mathcal{I}=\mathcal{N}$ by that of \eqref{dynGamma:agumented}. This motivates us to  introduce the following alternative condition for the Nash game that is customized for the augmented Riccati equation $\bar{{\bf \Gamma}}^\mathcal{N}$ : 
\begin{assume}{2B}
  \label{Fullcond2:Nash}
  {\color{black}There exist a symmetric matrix $\bar{\Pi}^\mathcal{N}_1 \in \mathbb{R}^{{(3n_1+2(n_{\alpha}+n_{\beta}))\times (3n_1+2(n_{\alpha}+n_{\beta}))}}$, a matrix   $\bar{\Pi}^{\mathcal{I}}_2 \in \mathbb{R}^{{(3n_1+2(n_{\alpha}+n_{\beta})})\times (3n_1+2(n_{\alpha}+n_{\beta}))}$, a matrix $\hat{{\bf A}}^\mathcal{N} \in \mathbb{R}^{{(n_{\alpha}+n_{\beta})\times (n_{\alpha}+n_{\beta})}}$, and a matrix $\hat{{\bf B}}^\mathcal{N} \in \mathbb{R}^{{(n_{\alpha}+n_{\beta})\times (3n_1+2(n_{\alpha}+n_{\beta}))}}$, with}
    \begin{equation*}
       \bar{\mathbf{\Lambda}}^\mathcal{N}:= \left(\begin{matrix}
            -\bar{\Pi}^\mathcal{N}_1\bar{{\bf A}}^{\mathcal{N}}+\bar{\Pi}^\mathcal{N}_2\bar{{\bf D}}^\mathcal{N} & \bar{\Pi}^\mathcal{N}_1\bar{{\bf B}}^\mathcal{N}-(\bar{{\bf A}}^\mathcal{N})^{\top}\bar{\Pi}^\mathcal{N}_2+\bar{\Pi}^\mathcal{N}_2{\bf C}^\mathcal{N}\\ {\bf 0} & (\bar{{\bf B}}^\mathcal{N})^{\top}\bar{\Pi}^\mathcal{N}_2
        \end{matrix}\right),
    \end{equation*}
    such that
            \begin{equation}
            \label{eq:new:condition:Nash}
                \begin{aligned}
                    \lambda_{\min}\left( \bar{\mathbf{\Lambda}}^\mathcal{N}+(\bar{\mathbf{\Lambda}}^\mathcal{N})^{\top}\right) &\geq 0,\\  \text{and } \lambda_{\min}\left(\bar{\Pi}^\mathcal{N}_1+\bar{\Pi}^\mathcal{N}_2\bar{{\bf E}}^\mathcal{N}+(\bar{{\bf E}}^\mathcal{N})^{\top}(\bar{\Pi}^\mathcal{N}_2)^{\top} \right)&>0. 
                \end{aligned}
            \end{equation}
     
\end{assume}

}

{\color{black} The flexibility of choosing $\hat{{\bf A}}^\mathcal{N}$ and $\hat{{\bf B}}^\mathcal{N}$ has made it  easier to find simple matrices $\bar{\Pi}^\mathcal{N}_1$ and $\bar{\Pi}^\mathcal{N}_2$, such as the identify matrix, that fulfills Assumption \ref{Fullcond2:Nash}. This allows us to consider a broad spectrum of parameters that is relevant to real-life applications. Combining Assumptions \ref{Fullcond2} and \ref{Fullcond2:Nash}, we introduce the following assumption:  

\begin{assume}{2}
 \label{ass:2:combined}
 At least one of the following holds:
    \begin{enumerate}
        \item Assumption \ref{Fullcond2} is satisfied for both $\mathcal{I}=\mathcal{N}$ and $\mathcal{P}$; 
        \item Assumption \ref{Fullcond2} is satisfied for $\mathcal{I}=\mathcal{P}$, and Assumption \ref{Fullcond2:Nash} is satisfied for $\mathcal{I}=\mathcal{N}$. 
    \end{enumerate}
\end{assume}
}

The following is the main result in this subsection, {\color{black}which is a direct consequence of Theorem 3.6.6 of \cite{matrix:riccati}, and the locally Lipschitz property of Riccati equations; see also Lemma \ref{unique:Riccati} in  Appendix \ref{sec:app:FBSDE:riccati}.} 
{\color{black} 
\begin{theorem}
\label{wellpose:dominate}
    Under Assumption \ref{ass:2:combined},  the Riccati equation \eqref{dynGamma} admits a unique solution for $\mathcal{I}=\mathcal{N}$ and $\mathcal{P}$. As a result, the FBSDEs \eqref{NumericalNG} and \eqref{NumericalPG}  are uniquely solvable,  where the adjoint process can be written as ${\bf p}^\mathcal{I}(t)={\bf \Gamma}^\mathcal{I}(t) {\bf x}^\mathcal{I}(t) + {\bf g}^\mathcal{I}(t)$, where ${\bf g}^\mathcal{I}$ is the solution of \eqref{NumericalAll2g}. 
\end{theorem}



}

\section{Cooperation versus Competition}
\label{sec:compare}
{\color{black} In this section, we provide the first systematic analysis which compares the effect of leaders' interactions on the community of followers. Our study indicates that} 
neither mode of the interactions between leaders will always be more favourable to the follower in terms of providing a smaller cost functional. {\color{black}Specifically, we conduct}  a comprehensive comparison of the cost functionals of the follower under the two leaders' games. 
To this end, we first represent the cost functional of the follower in terms of the solution of {\color{black}the respective} Lyapunov {\color{black}and} Riccati equations {\color{black}(see \eqref{dynGamma}, and \eqref{eq:K} below)}. The {\color{black} minimal} {\color{black} costs} of the follower under the two games can thus be compared by the solutions of the associated equations. However, there {\color{black}is} no general comparison theorem available for non-symmetric Riccati equations {\color{black}in the literature}. {\color{black} We examine one-dimensional cases, and} first consider a  broad class of simplified models (Model \ref{degeneratedmodel}), where the {\color{black} values} of the optimal cost functionals {\color{black}under the two modes of governance} can be {\color{black} compared} by a component-wise comparison of the {\color{black}corresponding} simplified Lyapunov {\color{black} and} Riccati equations. Then,  {\color{black}in Section \ref{sec:asymptotic:analysis},} we  extend the comparison to more sophisticated models which are fairly close to Model \ref{degeneratedmodel}  by a perturbation analysis.

\subsection{Cost Functional of the Follower}

{\color{black}We here} represent the optimal cost functional of followers in terms of the solution of {\color{black}the following} Lyapunov equation {\color{black} for $\mathcal{I}=\mathcal{N}$ and $\mathcal{P}$:}
     
    \begin{equation}
    \label{eq:K}
    \left\{
    \begin{aligned}
      - \frac{d{\bf K}^\mathcal{I}(t)}{dt} =&\ {\bf Q} + {\bf K}^\mathcal{I}(t)({\bf A}-{\bf B}{\color{black}^\mathcal{I}}{\bf \Gamma}^\mathcal{I}(t))  +({\bf A}-{\bf B}{\color{black}^\mathcal{I}}{\bf \Gamma}^\mathcal{I}(t))^{\top}{\bf K}^\mathcal{I}(t) + {\bf R}^\mathcal{I}(t),\\
      {\bf K}^\mathcal{I}(T)=&\ {\bf \bar{Q}}, 
        \end{aligned}\right.
    \end{equation}
{\color{black}and the following equation system:}  
    \begin{equation}
    \label{eq:k} \left\{
    \begin{aligned}
      - \frac{d{\bf k}^\mathcal{I}(t)}{dt} =&\ -{\bf q}Q_1M_1 - {\bf K}^\mathcal{I}(t){\bf B}{\color{black}^\mathcal{I}}{\bf g}^\mathcal{I}(t)   +({\bf A}-{\bf B}{\color{black}^\mathcal{I}}{\bf \Gamma}^\mathcal{I}(t))^{\top}{\bf k}^\mathcal{I}(t) \\ & \  + ({\bf \Gamma}^{\mathcal{I}}(t))^{\top} {\bf i}^{ {\color{black}\mathcal{I}}}_3   E_1R_1^{-1}E_1^\top ({\bf i}^{ {\color{black}\mathcal{I}}}_3 )^{\top}{\bf g}^{\mathcal{I}}(t) ,\\
      {\bf k}^\mathcal{I}(T)=&\ -\bar{\mathbf{q}} \bar{Q}_1 \bar{M}_1;
        \end{aligned}\right.
    \end{equation}
    \begin{equation}
    \label{eq:l}\left\{
    \begin{aligned}
      - \frac{dl^\mathcal{I}(t)}{dt}=&\   -2({\bf k}^\mathcal{I}(t))^{\top}{\bf B}{\color{black}^\mathcal{I}}{\bf g}^\mathcal{I}(t) +M_1^{\top}Q_1M_1 + ({\bf g}^{\mathcal{I}}(t))^{\top}{\bf i}^{ {\color{black}\mathcal{I}}}_3   E_1R_1^{-1}E_1^\top ({\bf i}^{ {\color{black}\mathcal{I}}}_3)^\top  {\bf g}^{\mathcal{I}}(t) ,\\
      l^\mathcal{I}(T)=&\  \bar{M}_1^{\top}\bar{Q}_1\bar{M}_1,
        \end{aligned}\right.
    \end{equation}
    {\color{black}where  ${\bf q} :=(-G_1^\top,-H_1^\top,(1-F_1)^\top,0,0)$, ${\bf \bar{q}} :=(-\bar{G}_1^\top,-\bar{H}_1^\top,1-\bar{F}_1)^\top,0,0)$,  ${\bf i}^\mathcal{N}_3 := (0,0,I,0,0$, $0,0)$, ${\bf i}^{{\color{black}\mathcal{P}}}_3 := (0,0,I,0,0)$, 
            ${\bf Q}:= \mathbf{q} Q_1 \mathbf{q}^{\top}$, ${\bf \bar{Q}}:=\bar{\mathbf{q}} \bar{Q}_1 \bar{\mathbf{q}}^{\top}$,  ${\bf R}^{\mathcal{I}}(t):=({\bf \Gamma}^{\mathcal{I}}(t))^{\top}  {\bf i}^{ {\color{black}\mathcal{I}}}_3  E_1R_1^{-1}E_1^\top ({\bf i}^{ {\color{black}\mathcal{I}}}_3 )^{\top} {\bf \Gamma}^{\mathcal{I}}(t)$ for $t\in[0,T]$.} The proof {\color{black}of Theorem \ref{thm:cost} below} is relegated to Appendix \ref{app:cost}.  
\begin{theorem}\label{thm:cost}
Let $\Gamma_{0}(t)$, $t\in[0,T]$ be the solution of the following Riccati equation  
{\color{black}\begin{equation}
\label{Gamma0}\left\{
\begin{aligned}
   -\frac{d  \Gamma_0(t)}{dt}=&\ A^{\top}_1  \Gamma_0(t)+ \Gamma_0(t)A_1- \Gamma_0(t)E_1R^{-1}_1E_1^{\top} \Gamma_0(t)   + Q_1, \\
  \Gamma_0(T)=&\ \bar{Q}_1.
\end{aligned}\right.
\end{equation}}
Then we have:
\begin{enumerate}
    \item  The optimal control $u^\mathcal{I}_1$ of the representative follower is given by, for $t\in[0,T]$, 
       \begin{equation}
       \label{optimalfollowerform}
       \begin{aligned}
        u^{\mathcal{I}}_1(t) =&\ -R^{-1}_1 E^{\top}_1\bigg( \Gamma_{0}(t)(x^{\mathcal{I}}_1(t)-z^{\mathcal{I}}(t)) + ({\bf i}^{ {\color{black}\mathcal{I}}}_3 )^{\top}{\bf \Gamma}^{\mathcal{I}}(t) {\bf x}^{\mathcal{I}}(t)+({\bf i}^{ {\color{black}\mathcal{I}}}_3 )^{\top}{\bf g}^{\mathcal{I}}(t) \bigg) .
        \end{aligned}
    \end{equation}%

    \item  The {\color{black}optimal} cost functional of the representative follower {\color{black}$J_1^{\mathcal{I}}$ under the equilibrium solution of Problems \ref{P1}-\ref{P3}}
    can be written as  
        \begin{equation}
        \label{eq:lem:cost}
        \begin{aligned}
            J_1^{\mathcal{I}} =&\  \mathbb{E}[   {\bf x}^\top(0){\bf K}^{\mathcal{I}}(0){\bf x}(0)] +\int_0^T \textup{tr}(  {\bf \Sigma}^\top{\bf K}^{\mathcal{I}}(t)  {\bf \Sigma} )dt   +2\mathbb{E}[   {\bf x}^\top(0)  {\bf k}^{\mathcal{I}}(0)]  +l^\mathcal{I}(0)+c,
            \end{aligned}
        \end{equation}%
    where 
    \begin{equation*}
    \begin{aligned}
    c:=&\ \int_0^T \mathrm{tr}\Big(( \Gamma^{\top}_0(t)E_1R^{-1}_1E_1^{\top} \Gamma_0(t)+Q_1)  (\Sigma_0(t)+P_0(t)) \Big)dt \\
       &\     + \mathrm{tr}\left(\bar{Q}_1(\Sigma_0 (T)+P_0(T))\right),
            \end{aligned}
    \end{equation*}%
    and for $t\in[0,T]$,
    \begin{equation}
    \label{Gamma0Lambda0A}
    \begin{aligned}
    \Sigma_0(t)&:=\int_0^te^{\int_s^t \bar{A}_1(u) du}\sigma_1\sigma^{\top}_1e^{\int_s^t \bar{A}^{\top}_1(u) du}ds, \\
    P_0(t)&:=e^{\int_0^t \bar{A}_1(s) ds}  \textup{Cov}[\eta_1] e^{\int_0^t \bar{A}^{\top}_1(s) ds},\\
    \bar{A}_1(t)&:=A_1-E_1R^{-1}_1E_1^{\top} \Gamma_0(t).
\end{aligned}  
\end{equation}%
Here, $\textup{Cov}[\eta_1]$ is the variance-covariance matrix of $\eta_1$, and
    ${\bf K}^{\mathcal{I}}(t)$ 
    ${\bf k}^{\mathcal{I}}(t)$, $l^\mathcal{I}(t)$ are  the solutions of \eqref{eq:K}, \eqref{eq:k} and \eqref{eq:l}, respectively. 
  
\end{enumerate}
          
\end{theorem}

\subsection{Comparison in a Symmetric Model}
\label{sec:comparedegenerated}
In this subsection, we examine the following base model, chosen as a representative first analysis:

\begin{model}{${\bf b_0}$}
 We make the following assumptions {\color{black}on coefficients}:
\label{degeneratedmodel}
    \begin{enumerate}

        \item $d_1=d_\alpha=d_\beta=1$, $n_1=n_\alpha=n_\beta=1$, $m_1=m_\alpha=m_\beta=1$;
        \item $C_\alpha=C_\beta=F_\alpha=F_\beta=\bar{F}_\alpha=\bar{F}_\beta=0$;
        \item $M_1=M_\alpha = M_\beta =\bar{M}_1=\bar{M}_\alpha = \bar{M}_\beta = 0$;
        \item there exists $k>0$ such that
        \begin{align*}
        & A_{\alpha\beta}:=A_\alpha = B_\beta, \ B_\alpha = k^2 A_\beta ,\ D_{\alpha\beta}:=D_\alpha=D_\beta, \\ 
        & G_\alpha = k^2 G_\beta, \ Q_{\alpha} = k^{-2}Q_\beta, \  R_\alpha = k^{-2}R_\beta, \\  & C_1 = k^{-1} D_1,\ H_1=k G_1, \ \bar{G}_\alpha = k^2 \bar{G}_\beta, \\ 
        & \bar{Q}_{\alpha} = k^{-2}\bar{Q}_\beta, \ \bar{H}_1=k \bar{G}_1.
    \end{align*}
    
    \end{enumerate}
\end{model}

Under Model \ref{degeneratedmodel}, the dynamics for the two leaders and the representative follower read as follows:  
\begin{align}
&\begin{cases}
dx_{\alpha}(t)&=\Big(A_{\alpha\beta}x_{\alpha}(t)+k^2 A_{\beta}x_{\beta}(t) +D_{\alpha\beta}v_{\alpha}(t)\Big)dt+\sigma_{\alpha}dW^{\alpha}(t);\\
x_{\alpha}(0)&=\eta_{\alpha},\\
\end{cases}\\
&\begin{cases}\label{dynbeta}
dx_{\beta}(t)&=\Big(A_{\beta}x_{\alpha}(t)+A_{\alpha\beta}x_{\beta}(t) +D_{\alpha\beta}v_{\beta}(t)\Big)dt+\sigma_{\beta}dW^{\beta}(t);\\
x_{\beta}(0)&=\eta_{\beta},
 \end{cases}\\
&\begin{cases}
d{\color{black}x_{1}}(t)&=\Big(A_1{\color{black}x_{1}}(t)+B_1z(t)+C_1x_{\alpha}(t)    +kC_1x_{\beta}(t)+E_1{\color{black}v_{1}}(t)\Big)dt+\sigma_1d{\color{black}W^{1}}(t);\\
{\color{black}x_{1}}(0)&={\color{black}\eta_{1}},
\end{cases}
\end{align}
where $W^{\alpha}$, $W^{\beta}$, $W^{1}$, $x_{\alpha}$, $x_{\beta}$, $x_{1}$, $v_{\alpha}$, $v_{\beta}$, $v_{1}$ are all real-valued process.
Their cost functionals, $J_{\alpha}(v_{\alpha};v_{\beta})$, $J_{\beta}(v_{\beta};v_{\alpha})$ and ${\color{black}J_{1}}={\color{black}J_{1}}({\color{black}v_{1}};x_\alpha,x_\beta,z)$ are respectively
    \begin{align*}
        J_{\alpha}(v_{\alpha};v_{\beta})=&\ \mathbb{E}\bigg[\int_0^T \big(\big|x_{\alpha}(t)-k^2G_{\beta}x_{\beta}(t)\big|_{\frac{Q_{\beta}}{k^2}}^2 +|v_{\alpha}(t)|_{\frac{R_{\beta}}{k^2}}^2\big)dt +\big|x_{\alpha}(T) -k^2\bar{G}_{\beta}x_{\beta}(T)\big|_{\frac{\bar{Q}_{\beta}}{k^2}}^2\bigg], \\
         J_{\beta}(v_{\beta};v_{\alpha})=&\ \mathbb{E}\bigg[\int_0^T \big(\big|x_{\beta}(t)-G_{\beta}x_{\alpha}(t)\big|_{Q_{\beta}}^2 +|v_{\beta}(t)|_{R_{\beta}}^2\big)dt +\big|x_{\beta}(T) -\bar{G}_{\beta}x_{\alpha}(T)\big|_{\bar{Q}_{\beta}}^2\bigg], \\
    {\color{black}J_{1}} =&\  \mathbb{E}\bigg[\int_0^T\big(\big|{\color{black}x_{1}}(t)-F_1z(t)-G_1x_{\alpha}(t)-kG_1x_{\beta}(t)\big|_{Q_1}^2+|{\color{black}v_{1}}(t)|_{R_1}^2\big)dt \nonumber \\ &\ +\big|{\color{black}x_{1}}(T)-\bar{F}_1z(T) -\bar{G}_1x_{\alpha}(T) -k\bar{G}_1x_{\beta}(T)\big|_{\bar{Q}_1}^2\bigg]. 
 \end{align*}

         Model \ref{degeneratedmodel} is a simplified base model which is symmetric in the following sense.  First,  we assume that the behaviour of the two leaders are symmetric (up to the factor $k$) on how they depend on themselves and on each other. Second, the dependence of the leaders for the follower is the same, up to the factor $k$. {\color{black} Third, the leaders are unaffected by the followers, but we allow the followers to be influenced by the mean field effect and the leaders' states, which respectively signifies the indirect and direct effect of the leaders on the follower; the former effect is often found in reality: leaders normally affect the followers indirectly via the media. } 
         
         {\color{black} Despite its simplicity,} Model \ref{degeneratedmodel} is a very realistic one, in practice,   leaders often share reflexive goals and decisions; for instance,  the dependence of evolution of one leader's dynamics on the other should synchronize with that of the cost functionals; see the aligning relationships of $B_\alpha,A_\beta$ with $G_\alpha,G_\beta$ in Model \ref{degeneratedmodel}. 
        {\color{black} For example, two leading companies in a competitive market often share similar goals, at the most up to the companies' scale, where they  strive for market share {\color{black} expansion and profit} 
        maximization{\color{black},} and evaluate their performance relative to each other. Their decisions regarding development, pricing, and market entry may exhibit symmetric dynamics as they compete for the same customer base.}
        This model also allows us to compare the effect of the two games between the leaders on the follower with intuitive conditions, {\color{black}especially} on the coefficients {\color{black} representing the interactions between different leaders and followers, e.g., $A_\beta,B_\alpha,C_1,D_1$}, and shedding light on performing comparisons for more general models (see Section \ref{sec:asymptotic:analysis}). We remark that, as a consequence of the assumptions in Model \ref{degeneratedmodel}, the functions ${\bf k}^{\mathcal{I}}$, $l^{\mathcal{I}}$ in Theorem \ref{thm:cost} and ${\bf g}^{\mathcal{I}}$ in \eqref{NumericalAll2g} are all identically zero,  and the {\color{black}dimensions of} functions ${\bf \Gamma}^{\mathcal{I}}$, ${\bf K}^{\mathcal{I}}$ are reduced to  $3\times3$ matrices.

    Before proceeding to the main result of this subsection, we introduce the following technical condition. {\color{black}To this end, we define the following constants: $\tilde{A}:=A_1+A_{\alpha\beta}+B_{\alpha\beta}$,}
    
        \begin{align}
            \label{eq:U}
           U &:= \max\bigg\{\frac{R_1}{2E_1^2}\bigg(\sqrt{(2A_1+B_1)^2 + \frac{4E_1^2Q_1(1-F_1)}{R_1}}   + 2A_1+B_1 \bigg), (1-\bar{F}_1)\bar{Q}_1\bigg\},\\
            V&:=\frac{Q_1 G_1 E^2_1}{R_1},\ \bar{V}=\frac{\bar{Q}_1\bar{G}_1E_1^2}{R_1}. \nonumber
        \end{align}%
    
   \begin{assumption}
    \label{ass:compare}
     {\color{black}
     All of the followings hold:
     \begin{enumerate}
         \item  $\max\{ F_1, \bar{F}_1\}<1$;
         
         
     \item 
     $\min\{\vert G_{\alpha} \vert, \vert G_{\beta} \vert\}<1$ and  $\min\{\vert \bar{G}_{\alpha} \vert, \vert \bar{G}_{\beta} \vert\}<1$;

     \item $\max\{\vert B_{\alpha} \vert, \vert G_{\alpha} \vert,  \vert \bar{Q}_{\alpha}\bar{G}_{\alpha}  \vert\}>0$ and \\  $\max\{\vert C_1 \vert,  \vert G_1 \vert,  \vert \bar{Q}_1\bar{G}_1 \vert\}>0$;
     \end{enumerate}
    and at least one of the following holds:
     \begin{enumerate}
       \setcounter{enumi}{3}

      \item 
      $B_1< -\frac{2E_1^2 U}{R_1}$ and $\frac{1}{2}\leq F_1,\bar{F}_1 < 1$;
     
    \item  $C_1\neq 0$,  $\frac{\bar{V}}{C_1} \ge -1$, and   at least one of the following conditions holds:
       \begin{enumerate}
           \item 

    $\tilde{A} \leq \frac{V}{C_1}$;

           \item the terminal time $T$ satisfies
           
            \begin{equation}
           \label{eq:ass:2}
            T < \begin{cases}
             \frac{ |C_1+\bar{V}|   }{\left|C_1 U-V\right|}, \ \text{if } \ \tilde{A} = U;\\
            \frac{\left| \log   \left| C_1 U+\bar{V}(U-\tilde{A})-V\right|-\log \left| C_1 \tilde{A}-V\right| \right| }{\vert U-\tilde{A} \vert }, \ \text{otherwise.}
            \end{cases}
        \end{equation}

       \end{enumerate}
           \end{enumerate}
       }
     \end{assumption}
      Assumption \ref{ass:compare} is a technical condition used to show the component-wise properties of the Lyapunov equation \eqref{eq:K} in  Proposition \ref{poK},  which eventually allows us to compare the cost functionals of the followers under the two {\color{black}leaders'} games using the representation in \eqref{eq:lem:cost}.  {\color{black} Condition 1) is a standard assumption  that limits excessive deviations from the mean field term in the representative follower's cost functional. Condition 2) suggests that each leader's objective should prioritize their own state over the alternative leader's state, ensuring that the influence from the other leader does not dominate their own. Condition 3) requires that each leader will impose certain influences over both the alternative leader and the followers. 
      
      To allow meaningful comparison, we require that the representative follower's state to be at least affected by the mean field term (Condition 4)) or the leaders' states (Condition 5)). In the former case, the positive coefficients $F_1$ and $\bar{F}_1$  encourages the representative follower to align with the community via the mean field term, which is sketched by the leaders. In the special scenario where $C_1=0$, Condition 4) about $B_1$ preserves a certain degree of influence from the leaders on the representative 
      follower's state via the mean field term, which allows a clear contrast of the effect on the representative follower between the two types of  governance.} 
      {\color{black}In the latter case when $C_1\neq 0$, we require} that the effect of leaders on the terminal cost of the follower is relatively small compared to their effects on the follower's state dynamics when $C_1G_1 \le 0$ and $C_1\bar{G}_1\leq0$, {\color{black}which are the cases we shall study in Theorem \ref{thm:compare}}. {\color{black}Condition} 5a) means that the overall state dynamics of the {\color{black}leaders and the representative follower} are self-decaying. Otherwise, we require a relatively small terminal time to avoid blow-up of the associated equations {\color{black}\eqref{dynGamma}}. Examples of model parameters which satisfy Assumption \ref{ass:compare} will be given in Section \ref{sec:numerical}.

          The following comparison in Theorem \ref{thm:compare} is the main result of this subsection, {\color{black}whose proof} is mainly due to component-wise comparison theorems of the associated Riccati {\color{black}and} Lyapunov systems {\color{black}\eqref{dynGamma} and \eqref{eq:K}},  requires a series of elementary computations and is relegated to   Appendix \ref{sec:app:compare}. To summarize, due to the symmetry {\color{black}of} Model \ref{degeneratedmodel}, we can take the advantage of the reduced dimensions and the coupling of the underlying Lyapunov {\color{black}and} Riccati systems. In this case, we can apply comparison theorems for ODEs to determine the relative size of the corresponding components of the associated Lyapunov {\color{black}and} Riccati equations under the two {\color{black}leaders'} games, which mainly depend on the signs of the coefficients representing the interactions between two different players. Using the representation \eqref{eq:lem:cost} and by taking into account of the size of the initial conditions, we are able to arrive at Theorem \ref{thm:compare}. 

    \begin{theorem}
    \label{thm:compare}
        Consider Model \ref{degeneratedmodel}. Let $J_1^{\mathcal{N},{\bf b}_0}$ and $J_1^{\mathcal{P},{\bf b}_0}$ be the optimal cost functional in \eqref{eq:lem:cost} for the {\color{black} representative} follower {\color{black}under} the Nash and the Pareto game, respectively. Suppose that Assumptions {\color{black}\ref{fixedptcond}-\ref{ass:compare}} hold. Then, {\color{black} given the initial conditions $\eta_\alpha,\eta_\beta$,} we have the following conclusions {\color{black}under different scenarios}:
            \begin{enumerate}
                \item Suppose that $B_{\alpha},A_\beta \ge 0$, ${\color{black} -1} < G_{\alpha},G_\beta,   \bar{G}_\alpha,\bar{G}_\beta \leq  0$, $C_1 \ge 0$ {\color{black} and $G_1, \bar{G}_1\le 0$}. Then, there exists $L>0$ independent of $\eta_1$, such that
                    \begin{enumerate}
                        \item  $J_1^{\mathcal{N},{\bf b}_0}>J_1^{\mathcal{P},{\bf b}_0}$ if $\mathbb{E}[\eta_1]\mathbb{E}[\eta_\alpha+k\eta_\beta]>-L$;
                        \item  $J_1^{\mathcal{N},{\bf b}_0}=J_1^{\mathcal{P},{\bf b}_0}$ if $\mathbb{E}[\eta_1]\mathbb{E}[\eta_\alpha+k\eta_\beta]=-L$; 
                        \item $J_1^{\mathcal{N},{\bf b}_0}<J_1^{\mathcal{P},{\bf b}_0}$ if $\mathbb{E}[\eta_1]\mathbb{E}[\eta_\alpha+k\eta_\beta]<-L$.
                    \end{enumerate}
                        

                \item Suppose that $B_{\alpha},A_\beta \ge 0$, ${\color{black} -1 }< G_{\alpha},G_\beta, {\color{black} \bar{G}_\alpha,\bar{G}_\beta  }\leq  0$, $C_1  \le 0$ {\color{black}and $G_1, \bar{G}_1\ge 0$}. Then, there exists $L>0$ independent of $\eta_1$, such that 
                    \begin{enumerate}
                        \item $J_1^{\mathcal{N},{\bf b}_0}>J_1^{\mathcal{P},{\bf b}_0}$ if $\mathbb{E}[\eta_1]\mathbb{E}[\eta_\alpha+k\eta_\beta]<L$; 
                        \item   $J_1^{\mathcal{N},{\bf b}_0}=J_1^{\mathcal{P},{\bf b}_0}$ if $\mathbb{E}[\eta_1]\mathbb{E}[\eta_\alpha+k\eta_\beta]=L$;
                        \item $J_1^{\mathcal{N},{\bf b}_0}<J_1^{\mathcal{P},{\bf b}_0}$ if $\mathbb{E}[\eta_1]\mathbb{E}[\eta_\alpha+k\eta_\beta]>L$.
                    \end{enumerate}

                        
                    
                 \item Suppose that  $B_{\alpha},A_\beta \le 0$, ${\color{black} 0  }\leq G_{\alpha},G_\beta, {\color{black} \bar{G}_\alpha,\bar{G}_\beta}  \leq  1$,  $C_1 \ge0$ {\color{black} and $G_1, \bar{G}_1\le 0$}. Then, there exists  $L>0$ independent of $\eta_1$, such that 
                    \begin{enumerate}
                        \item  $J_1^{\mathcal{N},{\bf b}_0}<J_1^{\mathcal{P},{\bf b}_0}$ if $\mathbb{E}[\eta_1]\mathbb{E}[\eta_\alpha+k\eta_\beta]>-L$; 
                        \item  $J_1^{\mathcal{N},{\bf b}_0}=J_1^{\mathcal{P},{\bf b}_0}$ if $\mathbb{E}[\eta_1]\mathbb{E}[\eta_\alpha+k\eta_\beta]=-L$;
                        \item  $J_1^{\mathcal{N},{\bf b}_0}>J_1^{\mathcal{P},{\bf b}_0}$ if $\mathbb{E}[\eta_1]\mathbb{E}[\eta_\alpha+k\eta_\beta]<-L$.
                    \end{enumerate}

                        

                 \item Suppose that  $B_{\alpha},A_\beta \le 0$, ${\color{black} 0 } \leq  G_{\alpha},G_\beta, {\color{black} \bar{G}_\alpha,\bar{G}_\beta < }  1$, $C_1 \le0$   and $G_1, \bar{G}_1\ge 0$. Then, there exists $L>0$ independent of $\eta_1$, such that
                    \begin{enumerate}
                        \item   $J_1^{\mathcal{N},{\bf b}_0}<J_1^{\mathcal{P},{\bf b}_0}$ if $\mathbb{E}[\eta_1]\mathbb{E}[\eta_\alpha+k\eta_\beta]<L$;
                        \item   $J_1^{\mathcal{N},{\bf b}_0}=J_1^{\mathcal{P},{\bf b}_0}$ if $\mathbb{E}[\eta_1]\mathbb{E}[\eta_\alpha+k\eta_\beta]=L$;
                        \item $J_1^{\mathcal{N},{\bf b}_0}>J_1^{\mathcal{P},{\bf b}_0}$ if $\mathbb{E}[\eta_1]\mathbb{E}[\eta_\alpha+k\eta_\beta]>L$.
                    \end{enumerate}
                 
                        
                
            \end{enumerate}
    \end{theorem}

 By looking at the state dynamics $x_1$ of the follower under Model \ref{degeneratedmodel}, the conditions  $C_1 \ge0$, $G_1, \bar{G}_1\leq0$ and $\mathbb{E}[\eta_1]\mathbb{E}[(\eta_\alpha+k\eta_\beta)]>0$ indicate the magnitude of the state dynamics of the follower $x_1$ tends to increase with that of the sum of the leaders, while the cost functional of the {\color{black} representative} follower penalizes {\color{black}with} the {\color{black} square of the sum of his} state dynamics with that of the leaders.  Intutively, the follower's interest is in line with that of both leaders. Similar intuition can be given to the conditions $C_1 \le0$, $G_1, \bar{G}_1 \geq0, \mathbb{E}[\eta_1]\mathbb{E}[\eta_\alpha+k\eta_\beta] <0$ by considering the reflected state dynamics $-x_1$ of the {\color{black}representative} follower. On the other hand, under the conditions  $C_1 \le0$, $G_1, \bar{G}_1 \geq0, \mathbb{E}[\eta_1]\mathbb{E}[\eta_\alpha+k\eta_\beta] >0$, the state dynamics of the {\color{black}representative} follower $x_1$ tends to decrease with that of the sum of the leaders, and the cost functional {\color{black}penalizes with the square of} the difference between the {\color{black}representative} follower's state dynamics with that of the leaders. In this case, the {\color{black}representative} follower's interest is in conflict with that of the leaders. Similar interpretations can be given to the conditions $C_1 \ge0$, $G_1, \bar{G}_1\leq0$ and $\mathbb{E}[\eta_1]\mathbb{E}[(\eta_\alpha+k\eta_\beta)]<0$, again, by considering the reflected state dynamics $-x_1$ of the {\color{black}representative} follower. The constant $L>0$ in Theorem \ref{thm:compare} is the threshold which measures by how much the initial conditions can or should deviate from {\color{black}the above positivity or negativity requirements}\footnote{{\color{black}The exact value  of $L$ in the above four cases can be written down explicitly in terms of the solutions of the associated Lyapunov equations, which are provided in Appendix \ref{app:proof:compare}.}}. Similarly, we interpret the conditions $B_\alpha,A_\beta \ge0$, $G_\alpha,G_\beta, {\color{black}\bar{G}_\alpha,\bar{G}_\beta}  \le0$ as the scenario when the leaders share a common interest; while they have conflicting interests when $B_\alpha,A_\beta \le0$, and $G_\alpha,G_\beta,\bar{G}_\alpha, \bar{G}_\beta  \ge0$.

 With the above discussions,  Theorem \ref{thm:compare} can be interpreted in the following ways. When the {\color{black}representative} follower's interest is in line with that of the leaders, he will be better off from leaders' cooperation if their interests are also in line with each other (Statements 1a) and 2a)); otherwise the {\color{black}representative} follower will be worsen off from the cooperation if the leaders have conflicting interests (Statements 3a) and 4a)). The opposite is observed when the {\color{black}representative} follower's interest is conflicting with that of the leaders, where they will be worsen off from leaders' cooperation if their interests are in line with each other (Statements 1c) and 2c)); and {\color{black}being} better off from leaders' cooperation when the they have conflicting interests (Statements 3c) and 4c)). Mathematically, under Model \ref{degeneratedmodel}, the difference of the optimal cost functionals is essentially affine in $\mathbb{E}[\eta_1]$ {\color{black}for given $\eta_\alpha,\eta_\beta$}, i.e.,  $J^{\mathcal{N},{\bf b}_0}_1-J^{\mathcal{P},{\bf b}_0}_1 = a{\color{black}\mathbb{E}[\eta_\alpha+k\eta_\beta]}\mathbb{E}[\eta_1] + b$, {\color{black} where  $a,b$ are independent of $\eta_1$, with} $a{\color{black}>} 0$, $b>0$ under Statement 1); $a{\color{black}<} 0$, $b> 0$ under Statement 2); $a{\color{black}<} 0$, $b< 0$ under Statement 3); and $a{\color{black}>} 0$, $b<0$ under Statement 4).

\subsection{Perturbation Analysis}
\label{sec:asymptotic:analysis}
We then extend the result in Theorem \ref{thm:compare} by considering models which {\color{black}are close enough to} Model \ref{degeneratedmodel}. For simplicity, we again consider the case $d_1=d_\alpha=d_\beta=n_1=n_\alpha=n_\beta=m_1=m_\alpha=m_\beta=1$. Denote a model by $\mathfrak{M}:=(\Theta,\boldsymbol{\eta})$, where $\Theta:=(\theta_i)_{i=1}^{44}$ and $\boldsymbol{\eta}$ are the model parameters defined by
\begin{align*}
            \Theta :=  &\ \big(A_\alpha,B_\alpha,C_\alpha,D_\alpha,\sigma_\alpha,A_\beta,B_\beta,C_\beta,D_\beta,\sigma_\beta,A_1,B_1, \\ &\quad C_1,D_1,E_1,\sigma_1, F_{\alpha}, G_{\alpha},M_{\alpha}, Q_{\alpha},R_\alpha, F_\beta,G_\beta,M_{\beta}, \\ &\quad Q_\beta,R_\beta,  F_{1}, G_{1},H_1,M_{1}, Q_{1},R_1, \bar{F}_{\alpha}, \bar{G}_{\alpha},\bar{M}_{\alpha}, \bar{Q}_{\alpha}, \\ 
            &\quad \bar{F}_\beta,\bar{G}_\beta,\bar{M}_{\beta}, \bar{Q}_\beta,\bar{F}_{1}, \bar{G}_{1},\bar{H}_1,\bar{M}_{1}, \bar{Q}_{1} \big),\\
            \boldsymbol{\eta} :=&\ (\eta_\alpha,\eta_\beta,\eta_1).
        \end{align*}
 We define the following  which measures the distance of two models. 
    \begin{definition}
        Let $\mathfrak{M}_1:=(\Theta_1,\boldsymbol{\eta}_1)$ and $\mathfrak{M}_2:=(\Theta_2,\boldsymbol{\eta}_2)$ be two models. The distance $\mathcal{D}$ between $\mathfrak{M}_1$ and $\mathfrak{M}_2$ is defined by 
        \begin{equation}
            \mathcal{D}(\mathfrak{M}_1,\mathfrak{M}_2) := \max\left\{ \|\Theta_1-\Theta_2\|_\infty, \mathbb{E}[|\boldsymbol{\eta}_1-\boldsymbol{\eta}_2|^2 ]  \right\} ,
        \end{equation}
        where  $\|\Theta\|_\infty$ is the sup-norm of the vector $\Theta${\color{black}, where the natural $L^2$-norm is used for every vector or matrix element in $\boldsymbol{\Theta}$},  and $|\cdot|$ denotes the Euclidean norm.
    \end{definition}
Consider a base model $\mathfrak{M}_0$ such that $J_1^{\mathcal{N},0} > J_1^{\mathcal{P},0}$ (resp. $J_1^{\mathcal{N},0} < J_1^{\mathcal{P},0}$), where $J_1^{\mathcal{N},0}$ and $J_1^{\mathcal{P},0}$ are respectively, the optimal cost functionals of the {\color{black}representative} follower under the Nash and the Pareto games for model $\mathfrak{M}_0$. We show that the ordering of the corresponding cost functional of the {\color{black}representative} follower under the two games will be preserved for any model $\mathfrak{M}$ that only deviates from $\mathfrak{M}_0$ mildly as measured by $\mathcal{D}$. {\color{black} To this end, for $\mathcal{I}=\mathcal{N},\mathcal{P}$, and for any $\delta>0$ and model $\mathfrak{M}$ satisfying $\mathcal{D}(\mathfrak{M},\mathfrak{M}_0)<\delta$, define 
\begin{equation}
\label{eq:E:delta}
\begin{aligned}
    \mathcal{E}_0^\mathcal{I}(\delta) := &\ \frac{1}{2\sqrt{\rho({\bf B})}}\bigg( \rho_T\big({\bf A} - {\bf B}{\bf \Gamma}^\mathcal{I}_0\big)   +   \rho_T\big({\bf C}^\mathcal{I}  - {\bf \Gamma}^\mathcal{I}_0{\bf B}\big)\bigg) \\ &\times \bigg\{   \rho_T\big({\bf \Gamma}^\mathcal{I}_0\tilde{{\bf A}} + \tilde{{\bf C}}^\mathcal{I} {\bf \Gamma}^\mathcal{I}_0 \big)    + \rho(\tilde{{\bf E}}^\mathcal{I}) \Big( \rho_T\big({\bf A} -  {\bf B}{\bf \Gamma}^\mathcal{I}_0\big)   +  \rho_T\big({\bf C}^\mathcal{I}  - {\bf \Gamma}^\mathcal{I}_0 {\bf B}\big)   \Big) \bigg\}^{-\frac{1}{2}},
        \end{aligned}
\end{equation}%
where ${\bf A}_0, {\bf A}$ are the matrices defined in \eqref{bigmatrixA} under $\mathfrak{M}_0$ and $\mathfrak{M}$, respectively, $\tilde{{\bf A}} := {\bf A}-{\bf A}_0$, and similarly for the other matrices; and ${\bf \Gamma}_0^\mathcal{I}$ is the solution of \eqref{dynGamma} under $\mathfrak{M}_0$. Notice that when $\mathcal{D}(\mathfrak{M},\mathfrak{M}_0)<\delta$, $\mathcal{E}_0^\mathcal{I}(\delta) = O(\delta^{-\frac{1}{2}})$. As shown in Theorem \ref{thm:gamma:norm}, when $\mathcal{E}^\mathcal{I}_0(\delta)>1$, along with a size condition on the terminal time $T$, \eqref{eq:E:delta} allows us to bound the differences of the associated Riccati equations and thus the cost functionals of the {\color{black}representative} follower under $\mathfrak{M}$ and $\mathfrak{M}_0$.} The following is the main result in this subsection, whose proof is relegated to   Appendix \ref{sec:app:thm:extend}. 

\begin{theorem}\label{compare_extend}
Consider a base model $\mathfrak{M}_0$ which satisfies Assumptions \ref{fixedptcond} {\color{black}and \ref{ass:2:combined}}. Suppose also that  $J_1^{\mathcal{N},0} > J_1^{\mathcal{P},0}$ (resp. $J_1^{\mathcal{N},0} < J_1^{\mathcal{P},0})$. Let $\Delta J^0_1 :=|J_1^{\mathcal{N},0} -J_1^{\mathcal{P},0}|>0$.  Then there exists a constant $\mathcal{C}_{\mathfrak{M}_0,T}>0$ depending solely on $\mathfrak{M}_0$ and $T$, such that, for any model $\mathfrak{M}$ satisfying 
    \begin{enumerate}[label=\arabic*)]
        \item $\mathcal{D}(\mathfrak{M},\mathfrak{M}_0)<\delta:= \mathcal{C}_{\mathfrak{M}_0,T}\Delta J_1^0$; and
        \item for $\mathcal{I}=\mathcal{N},\mathcal{P}$, 
            $\mathcal{E}_0^\mathcal{I}(\delta)>1$ and
            \begin{equation}
\label{eq:T:extension}
        T \leq    \frac{\log(\mathcal{E}^\mathcal{I}_0(\delta))}{\rho_T\big({\bf A} -  {\bf B} {\bf \Gamma}^\mathcal{I}_0\big)  + \rho_T\big({\bf C}^\mathcal{I}  - {\bf \Gamma}^\mathcal{I}_0 {\bf B}\big)},
    \end{equation}
    \end{enumerate}
we have $J_1^{\mathcal{N}} > J_1^{\mathcal{P}}$ (resp.~$J_1^{\mathcal{N}} < J_1^{\mathcal{P}})$, where $J^{\mathcal{N}}_1$ and $J_1^{\mathcal{P}}$ are respectively, the optimal cost functionals of the {\color{black}representative} follower under the Nash and the Pareto games for this model $\mathfrak{M}$.
\end{theorem}


  By choosing the base model $\mathfrak{M}_0$ to be $\mathfrak{M}_{{\bf b}_0}$, Theorems \ref{thm:compare} and  \ref{compare_extend} allow us to decide which mode of interaction between the leaders is more beneficial to the followers when the model parameters do not deviate significantly from the one specified in Model \ref{degeneratedmodel}{\color{black}, for instance, when the mean field effect on the leaders is small.} We remark that the actual values of $\Delta J^1_0$ and $\mathcal{C}_{\mathfrak{M}_{{\bf b}_0},T}$  can be deduced by following the lines of calculations in  Appendix \ref{sec:app:thm:extend}. However, writing down the explicit expression requires intricate calculations, {\color{black}and thus we leave the details in  Appendix \ref{sec:app:thm:extend}}. In the next section, we shall consider some numerical examples and illustrate the feasible size of the deviation $\delta$. 
  

\section{Numerical Studies}\label{sec:numerical}
Using the Riccati and Lyapunov differential equations, we can compute numerically the optimal ``cost-to-go" functions of the {\color{black}representative} follower, $J^\mathcal{I}_1(t)$, $0\leq t\leq T$, $\mathcal{I}=\mathcal{N}, \mathcal{P}$, using the optimal strategy derived over the horizon $[0,T]$: 
\begin{align*}
        J_1^\mathcal{I}(t) &:= \mathbb{E}\bigg[ \int_t^T \big|x^\mathcal{I}_1(s)-F_1z^\mathcal{I}(s)-G_1x^\mathcal{I}_\alpha(s)   -  H_1x^\mathcal{I}_\beta(s)-M_1\big|^2_{Q_1}   +|u_1^\mathcal{I}(s)|_{R_1}^2 ds \\ & \ \quad + \big|x_1^\mathcal{I}(T)  - \bar{F}_1z^\mathcal{I}(T)-\bar{G}_1x_\alpha^\mathcal{I}(T) -\bar{H}_1x_\beta^\mathcal{I}(T) - \bar{M}_1 \big|^2_{\bar{Q}_1}    \bigg],
    \end{align*}%
{\color{black} with the initial condition ${\bf x}_0$.} The objective of this section is twofold. First, we verify the comparison result in Theorem \ref{thm:compare} under Model \ref{degeneratedmodel} ($\mathfrak{M}_{{\bf b}_0}$) by looking at the cost-to-go functions. Second, we illustrate the feasible deviation from $\mathfrak{M}_{{\bf b}_0}$  such that the comparison result for $\mathfrak{M}_{{\bf b}_0}$ is still preserved for models within that deviation. For simplicity, we assume that $\eta_1$, $\eta_\alpha$ and $\eta_\beta$ are real constants.

\subsection{Cooperation versus Competition under Model \ref{degeneratedmodel}}
\label{sec:num:simple}

We verify the comparison result in Theorem \ref{thm:compare} and examine the size of the threshold $L$. In particular, we fix $T=2$ and consider the following choices of parameters for the {\color{black}representative} follower: {\color{black}$A_1=-1.2,\ B_1=0.2, \ C_1=0.1, \ D_1=0.2, \   E_1=1,\ F_1=0.2 , \ G_1=-0.1,\ H_1=-0.2,\ Q_1=1,\ R_1=1,\ \bar{F}_1=  \bar{G}_1 = \bar{H}_1 =  \bar{Q}_1=0, \ \sigma_1=0.1$; $\alpha$-leader: $ A_{\alpha}=-1, \  C_\alpha = 0,\ D_{\alpha}=1,\ F_{\alpha}=0,\ Q_{\alpha}=1,\ R_{\alpha}=1,\ \sigma_{\alpha}=0.1,\ \bar{G}_\alpha =  \bar{F}_\alpha = \bar{Q}_\alpha = 0,\ \eta_{\alpha}=2$; and $\beta$-leader: $B_{\beta}=-1,\ {\color{black}C_{\beta}=0,\ } D_{\beta}=1,\ F_{\beta}=0, \ M_{\beta}=0,\ Q_{\beta}=4,\ R_{\beta}=4,\  \bar{G}_\beta =  \bar{F}_\beta  =\bar{Q}_\beta = 0,\  \sigma_{\beta}=0.1,\ \eta_{\beta}=2$.}
We also consider different values of $\eta_1$ and the following two individual sets of interaction parameters between the leaders:
 {\color{black}    \begin{equation}
    \label{eq:num:compare:1}
    \tag{Ia}
    B_\alpha = 0.4, \ A_\beta = 0.1, \ G_\alpha = -0.4, \ G_\beta = -0.1
    \end{equation}
    \begin{equation}
    \tag{Ib}
        \label{eq:num:compare:2}
       B_\alpha = -0.4, \ A_\beta = -0.1, \ G_\alpha = 0.4, \ G_\beta = 0.1.
    \end{equation}}%
The above choices of parameters follow Model \ref{degeneratedmodel} with factor $k=2$.

{\color{black} To ensure the well-posedness of the problem, by Theorems \ref{ZMEX} and \ref{wellpose:dominate}, we shall verify that the chosen parameters satisfy Assumptions \ref{fixedptcond} and \ref{ass:2:combined}. By a direct calculation, it is straightforward to verify Assumption \ref{fixedptcond}. On the other hand, by choosing all the matrices $\Pi_1^\mathcal{P}, \Pi_2^\mathcal{P}, \bar{\Pi}_1^\mathcal{N}, \bar{\Pi}_2^\mathcal{N}$ to be identity matrices with suitable dimensions, $\hat{\bf A}^\mathcal{N} = \text{diag}(B_\beta,A_\alpha)$, and 
$$ \hat{\bf B}^\mathcal{N} = \begin{pmatrix}
  0 & 0 & 0 & 0 & 0 &  D_\beta^2/R_\beta& 0\\ 0 & 0 & 0 & 0 & 0 &0 &D_\alpha^2/R_\alpha
\end{pmatrix},$$ it can be verified that Assumption \ref{Fullcond2} (resp. Assumption \ref{Fullcond2:Nash}) is fulfilled for the Pareto (resp. Nash) game. To apply the comparison result in Theorem \ref{thm:compare} for Model \ref{degeneratedmodel}, we also need Assumption \ref{ass:compare}, which  can again be verified by direct calculations.   } 

\begin{figure}[!h]
    \centering
    \includegraphics[scale=0.4]{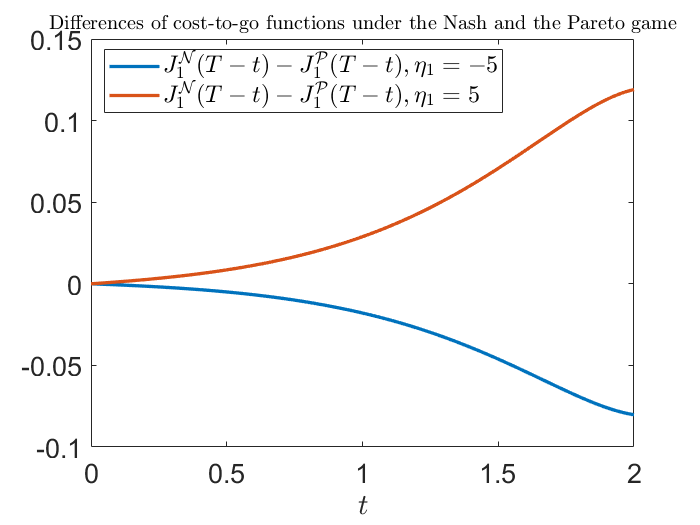}
    \caption{$J_1^\mathcal{N}(T-t)-J^\mathcal{P}_1(T-t)$ with parameters \eqref{eq:num:compare:1} and $\eta_1=\pm5$}
    \label{fig:CTG_D1}
\end{figure}

\begin{figure}[!h]
    \centering
    \includegraphics[scale=0.4]{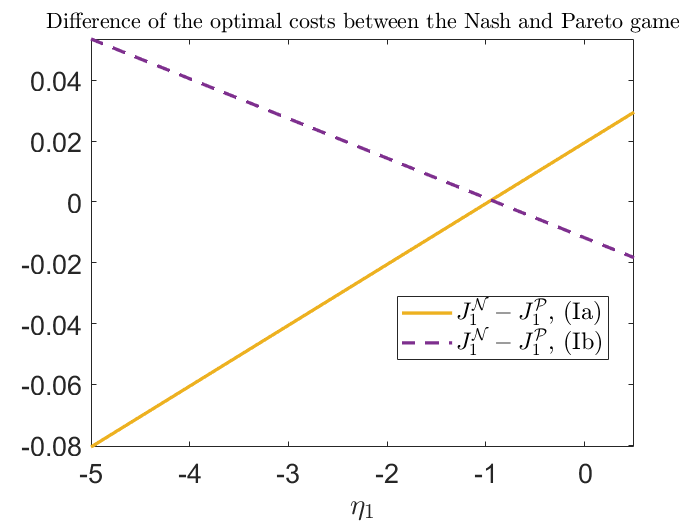}
    \caption{$J_1^\mathcal{N}-J^\mathcal{P}_1$ with parameters \eqref{eq:num:compare:1} (solid) and \eqref{eq:num:compare:2} (dashed) with $\eta_1\in [-5,0.5]$}
    \label{fig:DJ1}
\end{figure}

Figure \ref{fig:CTG_D1} shows the plot of the difference of the cost-to-go functions between  the Nash and the Pareto games with the parameters \eqref{eq:num:compare:1} and $\eta_1 = \pm5$. In particular, we observe that the Pareto game between leaders renders a lower cost-to-go function when $\eta_1=5$, since the interests of both leaders and the followers align with each other. However, this advantage is offset by an excessively negative $\eta_1$, and a smaller cost-to-go function is observed under the Nash game when $\eta_1=-5$. Specifically, the Nash game between leaders are more favourable to leaders when $\eta_1< {\color{black} -0.9731}$; see Figure \ref{fig:DJ1}, which plots the differences of values of the optimal cost functionals with respect to $\eta_1$. This verifies Statement 1) of Theorem \ref{thm:compare}. In the same figure, we observe that under \eqref{eq:num:compare:2}, the Pareto game is favourable to followers only when $\eta_1<{\color{black} -0.9030}$, since the two leaders have conflicting interests. This agrees with {\color{black}our} findings in Statement {\color{black} 3}) of Theorem \ref{thm:compare}.


\subsection{Non-symmetric Example}
By Theorem \ref{thm:compare}, we are able to {\color{black}decide} which mode of interaction between leaders are more favourable to the followers under Model \ref{degeneratedmodel} ($\mathfrak{M}_{{\bf b}_0}$), and by Theorem \ref{compare_extend}, we know that the same conclusion can be carried to any model $\mathfrak{M}$ satisfying $$\mathcal{D}(\mathfrak{M},\mathfrak{M}_{{\bf b}_0})<  O(\min\{ \Delta J^{{\bf b}_0}_1,(\Delta J^{{\bf b}_0}_1)^2\}),$$ where $\Delta J^{{\bf b}_0}_1:=|J^{\mathcal{N},{\bf b_0}}-J^{\mathcal{P},{\bf b_0}}|$. The feasible size of the deviation is examined under a particular example herein. Specifically, {\color{black}our study is based on the same set of parameters in Section \ref{sec:num:simple} with interaction coefficients \eqref{eq:num:compare:1}, except the following changes:}
    \begin{enumerate}
        \item {\color{black}  We include the terminal costs for the follower: $\bar{F}_1 = 0.2, \ \bar{G}_1 = -0.1, \ \bar{H}_1 = -0.2, \  \bar{Q}_1=0.2$; $\alpha$-leader: $\bar{G}_\alpha = -0.4, \ \bar{F}_\alpha = 0, \ \bar{Q}_\alpha = 0.2$; and $\beta$-leader: $\bar{G}_\beta = -0.1, \ \bar{F}_\beta = 0, \ \bar{Q}_\beta = 0.8$. }
        \item  We perturb the following parameters: for $\delta\in \mathbb{R}$, set
    \begin{equation}
        \label{eq:numer:perturb}
        \begin{aligned}
           &  C_1 = 0.1+\delta, \ D_1 = 0.2-\delta, M_1 = \delta, \ B_\alpha = 0.4-0.25\delta, \ C_\alpha = \delta, \ C_\beta = \delta. 
        \end{aligned}
    \end{equation}
        \end{enumerate}
 When $\delta=0$, the model belongs to $\mathfrak{M}_{{\bf b}_0}$ with $k=2$  and satisfies Assumptions {\color{black}\ref{fixedptcond} and \ref{ass:2:combined}  with the same choice of auxiliary matrices as in the previous subsection.} By Theorem \ref{thm:compare}, when $\delta=0$ and $\eta_1>-L$ for some $L>0$, we have $J^{\mathcal{N},0,\eta_1}_1>J^{\mathcal{P},0,\eta_1}_1$, where $J^{\mathcal{I},\delta,\eta_1}_1$, $\mathcal{I}=\mathcal{N},\mathcal{P}$ denotes the cost functional of the {\color{black}representative} follower with the above choices of parameters with the perturbation $\delta$. By solving numerically, it is found that $L = {\color{black} -0.9735}$.
 
\begin{figure}[!h]
    \centering
    \includegraphics[scale=0.4]{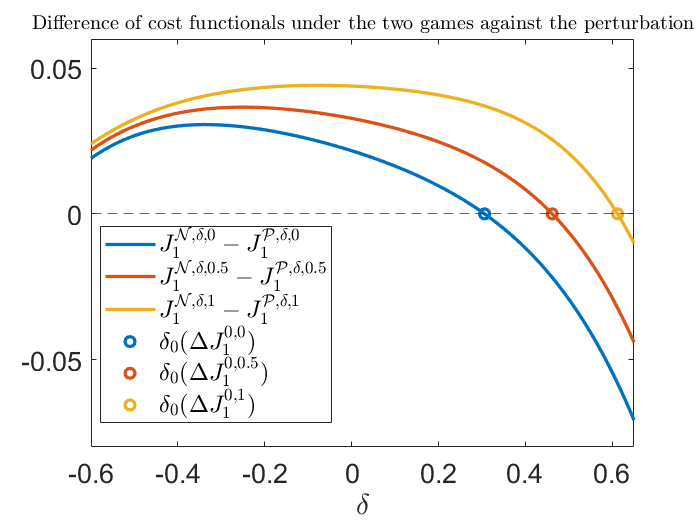}
    \caption{$J_1^{\mathcal{N},\delta,\eta_1}-J_1^{\mathcal{N},\delta,\eta_1}$ as a function of $\delta$ for $\eta_1=0$ (blue), {\color{black}$\eta_1=0.5$} (red) and {\color{black}$\eta_1=1$} (yellow)   }
    \label{fig:DJ_delta}
\end{figure}

 The introduction of $\delta$ herein breaks the symmetry of the underlying model and brings the mean field effect to the leaders. Denote the model with perturbation $\delta$ by $\mathfrak{M}_\delta$.  To illustrate the size of the feasible perturbation $\delta$ and its relationship with $\Delta J_1^{0,\eta_1} = |J_1^{\mathcal{N},0,\eta_1}-J_1^{\mathcal{P},0,\eta_1}|$, we vary the size of $\eta_1>0$ and compute 
    \begin{equation*}
        \delta_0 =  \delta_0(\Delta J_1^{0,\eta_1}) := \inf\left\{ |\delta| : J_1^{\mathcal{N},\delta,\eta_1}-J_1^{\mathcal{P},\delta,\eta_1} <0 \right\} ,
    \end{equation*}
which is the maximum perturbation such that the conclusion for $\mathfrak{M}_{{\bf b}_0}$ is also valid for any model $\mathfrak{M}_\delta$ with $|\delta|<\delta_0$. Figure \ref{fig:DJ_delta} illustrates the difference $J_1^{\mathcal{N},\delta,\eta_1}-J_1^{\mathcal{P},\delta,\eta_1}$ as a function of $\delta$ for $\eta_1=0,{\color{black}0.5,1}$, along with the corresponding maximum feasible perturbation $\delta_0(\Delta J_1^{0,\eta_1})$.\footnote{{\color{black}Using a grid search of 1,000 points over the interval $[-1,1]$, we find that both Assumptions \ref{fixedptcond} and \ref{ass:2:combined} are satisfied when $\delta\in {\color{black}[-0.719,0.721]}$  by some monotonicity observations. Assumption \ref{ass:compare} is only used for the comparison result in Theorem \ref{thm:compare} for Model \ref{degeneratedmodel}, which is not relevant when $\delta\neq 0$. }}

\begin{figure}[!h]
    \centering
    \includegraphics[scale=0.4]{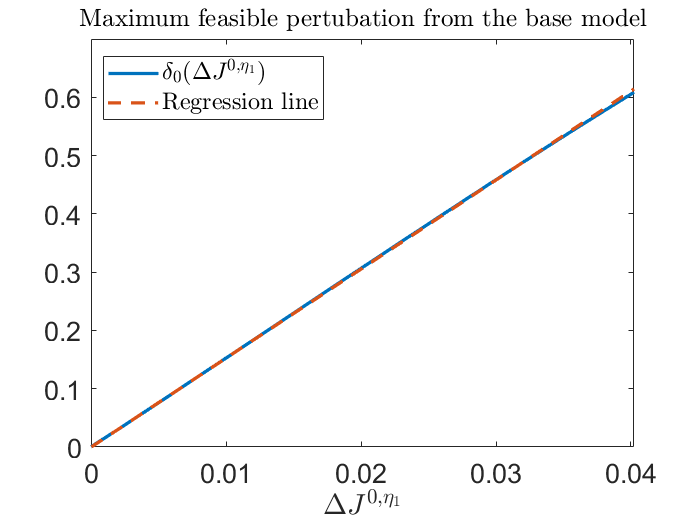}
    \caption{The maximum feasible perturbation $\delta_0$ against $\Delta J_1^{0,\eta_1}$ as $\eta_1$ varies from {\color{black}$-0.9735$ to $1$} (solid), and the best {\color{black}straight} line (dashed) {\color{black}$\Delta J_1^{0,\eta_1}={\color{black}e^{2.7248}\delta_0^{0.9997}}$} computed in \eqref{eq:regression}    }
    \label{fig:Perturb}
\end{figure}

Figure \ref{fig:Perturb} presents the size of $\delta_0$ with respect to $\Delta J^{0,\eta_1}_1$ as $\eta_1$ varies from $-L=-{\color{black}0.9735}$ to {\color{black}$1$}.
We observe that the size of $\delta_0$ increases nearly linearly with respect to $\Delta J^{0,\eta_1}_1$ (and thus with respect to $\eta_1$). By {\color{black}fitting the best straight line on}
    \begin{equation}
    \label{eq:regression}
        \log(\Delta J_1^{0,\eta_1}) = a \log(\delta_0) + b
    \end{equation}
for {\color{black}$\eta_1\in (-L,1]$, we find that {\color{black}$a =0.9997$} and {\color{black}$b=2.7248$}} (see the dashed line in Figure \ref{fig:Perturb}), which affirms our observation and the findings in Theorem \ref{compare_extend}. 
We also observe that the size $\delta_0$ is reasonably large {\color{black}when} comparing to the size of the {\color{black}coefficients of the state dynamics and the cost functionals}. 

        \label{fig:MF:follower}


\section{Conclusion}
\label{sec:conclusion}
In this article, we considered a natural extension of the single-leader linear-quadratic MFG in \cite{BCY2} by introducing multiple leaders over a group of followers. Depending on whether or not the leaders cooperate, we solved respectively the competitive Nash game and the cooperative Pareto game. Under the simplified model proposed in Section \ref{sec:compare}, one can decide which mode of interactions between leaders is more favourable to followers in terms of rendering a smaller optimal cost functional. We also extend the comparison result for models which are sufficiently close to the proposed simplified model by an asymptotic analysis. Due to intricate interactions between the leaders and followers, a general decisive condition is {\color{black}essentially absent before our work}. {\color{black}This article provided the first systematic analysis under a class of preliminary tractable models, and results for more generic models} shall be explored in future works.

\section*{Acknowledgment}

{\color{black}The authors would like to thank the editor and three anonymous referees for their valuable comments and suggestions that have helped to greatly improve the paper.} Dantong Chu acknowledges the support from UGC for his Ph.D. study at the Chinese University of Hong Kong (CUHK). Kenneth Ng acknowledges the financial support from UIUC and CUHK.  Phillip Yam acknowledges the financial supports from HKGRF-14301321 with the project title ``General Theory for Infinite Dimensional Stochastic Control: Mean Field and Some Classical Problems", and  HKGRF-14300123 with project title ``Well-posedness of Some Poisson-driven Mean Field Learning Models and their Applications". Harry Zheng acknowledges the financial support from  EPSRC (UK)  grant (EP/V008331/1) with the project title ``Deep  Learning for Time-Inconsistent Dynamic Optimization".

\nocite{BFY}
\bibliographystyle{acm}

\bibliography{mfg}

\begin{appendices}
 
\section{Comparison of Pareto Game with Two Leaders and One-Leader Game}
\label{sec:app:pareto:one}
{\color{black}This section discusses the relationship of a two-leader Pareto game with an one-leader game.} {\color{black}Consider an $\mathbb{R}^{n_{\alpha}+n_{\beta}}$-valued} process $x_0:= (x_\alpha,x_\beta)$ {\color{black} which represents the state of a combined leader}, and an $\mathbb{R}^{m_{\alpha}+m_{\beta}}$-valued  process $v_0:=(v_\alpha,v_\beta)$ representing the associated control process. Then, the dynamics of the combined leader and the representative follower are respectively given by  
\begin{align*}
&\begin{cases}
dx_0(t)&=\Big(A_0x_0(t)+{\color{black}C_0}z(t)+{\color{black}D_0}v_0(t)\Big)dt+\sigma_0dW^0(t);\\
x_0(0)&=\eta_0 := (\eta_\alpha\ \eta_\beta)^\top ,\\
\end{cases}\\
&\begin{cases}
dx_1(t)&=\Big(A_1x_1(t)+B_1z(t)+\left(\begin{matrix}
    C_1 & D_1
\end{matrix}\right)x_0(t)  +E_1{\color{black}v_{1}}(t)\Big)dt+\sigma_1dW^1(t);\\
{\color{black}x_{1}}(0)&=\eta_1,
\end{cases}
\end{align*}%
and the cost functionals of the {\color{black}combined leader, $J_0=J_0(v_0):=J_\alpha(v_\alpha;v_\beta) + J_\beta(v_\beta;v_\alpha)$} and the representative follower {\color{black}$J_1=J_1(v_1;x_0,z)$} can be respectively written as
 
\begin{align*}
        J_0 =&\ \mathbb{E}\bigg[\int_0^T \big(\big|S_0x_0(t)-E_0z(t) -M_0   \big|_{Q_0}^2  +|v_0(t)|_{R_0}^2\big)dt     +\big|\bar{S}_0x_0(T)
         -\bar{E}_0z(T)-\bar{M}_0\big|_{\bar{Q}_0}^2\bigg],\\
         J_1 =&\  \mathbb{E}\bigg[\int_0^T\big(\big|x_1(t)-F_1z(t)-\left(\begin{matrix}
    G_1 & H_1
\end{matrix}\right)x_0(t) - M_1\big|_{Q_1}^2  +|v_{1}(t)|_{R_1}^2\big)dt  +\big|x_{1}(T)-\bar{F}_1z(T)\nonumber\\
    &\quad -\left(\begin{matrix}
    \bar{G}_1 & \bar{H}_1
\end{matrix}\right)x_0(T)-\bar{M}_1\big|_{\bar{Q}_1}^2\bigg].
\end{align*} %
Here,
  \begin{align*}
 & A_0:=
\left(\begin{matrix}
A_{\alpha} & B_{\alpha} \\
A_{\beta} & B_{\beta}
\end{matrix}\right),\ {\color{black}C_0}:=
\left(\begin{matrix}
C_{\alpha}\\
C_{\beta}
\end{matrix}\right),\ {\color{black}D_0}:=
\left(\begin{matrix}
D_{\alpha} & 0 \\
0 & D_{\beta}
\end{matrix}\right), \\
& {\sigma}_0:=
\left(\begin{matrix}
 \sigma_{\alpha} & 0\\
 0 & \sigma_{\beta}
\end{matrix}\right),
\ W^0(t):=
\left(\begin{matrix}
W^{\alpha}(t)\\W^{\beta}(t)
\end{matrix}\right),\ S_0:=
\left(\begin{matrix}
I & -G_{\alpha} \\
-G_{\beta} & I
\end{matrix}\right),\\
& E_0:=
\left(\begin{matrix}
F_{\alpha}\\
F_{\beta}
\end{matrix}\right),  M_0:= 
\left(\begin{matrix}
M_{\alpha}\\
M_{\beta}
\end{matrix}\right), Q_0:=
\left(\begin{matrix}
Q_{\alpha} & 0 \\
0 & Q_{\beta}
\end{matrix}\right),\\
& R_0:=
\left(\begin{matrix}
R_{\alpha} & 0 \\
0 & R_{\beta}
\end{matrix}\right),  \bar{S}_0:=
\left(\begin{matrix}
I & -\bar{G}_{\alpha} \\
-\bar{G}_{\beta} & I
\end{matrix}\right), \ 
 \bar{E}_0:=
\left(\begin{matrix}
\bar{F}_{\alpha}\\
\bar{F}_{\beta}
\end{matrix}\right),\\
&\bar{M}_0:= 
\left(\begin{matrix}
\bar{M}_{\alpha}\\
\bar{M}_{\beta}
\end{matrix}\right),\
\bar{Q}_0:=
\left(\begin{matrix}
\bar{Q}_{\alpha} & 0 \\
0 & \bar{Q}_{\beta}
\end{matrix}\right) .
 \end{align*}%
By the above formulations, the Pareto game of multiple leaders  can be treated as the mean field game with one combined leader with a higher-dimensional state. In particular, when $\text{det}(S_0),\text{det}(\bar{S}_0)\neq 0$, the formulation can be considered as the linear-quadratic case in \cite{BCY}.

\section{Proofs in Section \ref{sec:sol} }
\label{app:bsde}
This section provides the proofs of the main results in Section \ref{sec:sol}.  



\subsection{Proof of Theorem \ref{ZMEX}}

 We verify that the system  \eqref{ZM}   satisfies Assumptions 2.1 and {\color{black}3.1} 
 in \cite{Peng}, and hence the result follows immediately from Theorem {\color{black}3.2} 
 therein. Given $x_\alpha\in\mathbb{R}^{d_\alpha}$ and $x_\beta\in\mathbb{R}^{d_\beta}$, let $\mathcal{T} : \mathbb{R}^{d_1} \times  \mathbb{R}^{d_1}  \to  \mathbb{R}^{d_1} \times \mathbb{R}^{d_1}$ and  $\mathcal{G}:\mathbb{R}^{d_1}\to \mathbb{R}^{d_1}$ be the operators defined by 
	$$\mathcal{T}\begin{pmatrix}
	z \\ \zeta
	\end{pmatrix} =  \begin{pmatrix} -A_1^{\top} \zeta-Q_1((I-F_1)z-G_1x_{\alpha}-H_1x_{\beta}-M_1)  \\ (A_1+B_1)z+C_1x_{\alpha}+D_1x_{\beta}-E_1R_1^{-1}E_1^{\top} \zeta      \end{pmatrix},$$
and 
\begin{equation*}
\mathcal{G}(z)=\bar{Q}_1\Big((I-\bar{F}_1)z-\bar{G}_1x_{\alpha}-\bar{H}_1x_{\beta}-\bar{M}_1\Big).
\end{equation*} 
Since $\mathcal{T}$ and $\mathcal{G}$ are linear, it is clear that Assumption 2.1 of \cite{Peng} is satisfied. 
To verify Assumption 2.2 of \cite{Peng}, by Young's inequality, we have, for any $\epsilon>0$ and $z_1,z_2,\zeta_1,\zeta_2 \in \mathbb{R}^{d_1}$,  
    \begin{align*}
      & \ \bigg\langle \mathcal{T} \begin{pmatrix}  z_1-z_2\\ \zeta_1-\zeta_2 \end{pmatrix}, \begin{pmatrix}  z_1-z_2\\ \zeta_1-\zeta_2  \end{pmatrix} \bigg\rangle\\
      =&\  \langle B_1 (z_1-z_2), \zeta_1-\zeta_2\rangle  -\langle Q_1(I-F_1)(z_1-z_2), z_1-z_2 \rangle  - \langle E_1R_1^{-1}E_1^{\top} (\zeta_1-\zeta_2),\zeta_1-\zeta_2  \rangle \\
        \leq&\ \frac{\rho(B_1)^2}{2\epsilon }|z_1-z_2|^2 +\frac{\epsilon}{2} |\zeta_1-\zeta_2|^2     -  \lambda_{\min}\left(\frac{Q_1(I-F_1)+(I-F_1)^{\top}Q_1}2\right) |z_1-z_2|^2 \\
        &\ -\lambda_{\min}(E_1R_1^{-1}E_1^{\top}) |\zeta_1-\zeta_2|^2 \\
        =&\  - |z_1-z_2|^2\bigg(\lambda_{\min}\left(\frac{Q_1(I-F_1)+(I-F_1)^{\top}Q_1}2\right)  - \frac{\rho(B_1)^2}{2\epsilon }\bigg)\\
        &\quad - |\zeta_1-\zeta_2|^2 \left(\lambda_{\min}(E_1R_1^{-1}E_1^{\top})-\frac{\epsilon}{2}\right).
    \end{align*}
    By \eqref{fixedptcond1}, we can find $\epsilon>0$ such that $\lambda_{\min}(E_1R_1^{-1}E_1^\top) > \epsilon/2$ and $\lambda_{\min}(Q_1(1-F_1)+(1-F_1)^\top Q_1)>\epsilon^{-1}\rho(B_1)^2$.
Hence, there exists $c_\epsilon>0$ such that  
\begin{equation*}
\bigg\langle \mathcal{T} \begin{pmatrix}  z_1-z_2\\ \zeta_1-\zeta_2 \end{pmatrix}, \begin{pmatrix}  z_1-z_2\\ \zeta_1-\zeta_2  \end{pmatrix} \bigg\rangle \le {\color{red} -} c_\epsilon (|z_1-z_2|^2+|\zeta_1-\zeta_2|^2),
\end{equation*}%
for any $z_1,z_2,\zeta_1,\zeta_2 \in \mathbb{R}^{d_1}$. On the other hand,  
\begin{equation*}
\begin{aligned}
&\ \langle \mathcal{G}(z_1)-\mathcal{G}(z_2),  z_1-z_2 \rangle \geq  \lambda_{\min}\left(\frac{\bar{Q}_1(I-\bar{F}_1)+(I-\bar{F}_1)^{\top}\bar{Q}_1}2\right) |z_1-z_2|^2 ,
\end{aligned}
\end{equation*}%
for any $z_1,z_2  \in \mathbb{R}^{d_1}$. Hence, by Assumption \ref{fixedptcond} herein, Assumption {\color{black}3.1} 
in \cite{Peng} is verified. \hspace*{\fill} \qed
 
\subsection{Proof of Theorem \ref{thm_Nash}}
	\color{black}
    Let $u_\alpha^\mathcal{N}$ and $u_\beta^\mathcal{N}$ be the optimal control for the $\alpha$-leader and $\beta$-leader, respectively. We deduce that $u_\alpha^{\mathcal{N}}$ is given by \eqref{Nashoptimalcontrol}, and the derivation of $u^{\mathcal{N}}_\beta$ can then be done by a parallel argument. {\color{black} Under the Nash game, the $\alpha$-leader assumes that the $\beta$-leader uses his optimal strategy $u^{\mathcal{N}}_\beta$, which is thus considered as an exogenous process.}  For $t\in[0,T]$ and $\tau>0$, consider a perturbation of the optimal control $u^\tau_\alpha(t):=u^{\mathcal{N}}_{\alpha}(t)+\tau \tilde{u}_{\alpha}(t)$, where $\tilde{u}_{\alpha}$ is an arbitrary square-integrable process adapted to $\mathcal{F}^{\alpha}\vee\mathcal{F}^{\beta}$. Let $x^\tau_\alpha$ and $z^\tau$ be respectively the state dynamics of the $\alpha$-leader and the mean field term under the pair of controls $u^\tau_\alpha$ and $u^{\mathcal{N}}_\beta$. By the linearity of the state dynamics, it is easy to see that $x^\tau_\alpha = x^{\mathcal{N}}_\alpha + \tau\tilde{x}_\alpha$, {\color{black}$x^\tau_\beta = x^{\mathcal{N}}_\beta + \tau\tilde{x}_\beta$} and $z^\tau = z^\mathcal{N}+\tau\tilde{z}$, where  ${\bf \Tilde{x}}_\alpha := (\tilde{x}_\alpha,{\color{black}\tilde{x}_\beta},\tilde{z})$ satisfies  
        \begin{equation*}
    \left\{ \begin{aligned}
        d\tilde{{\bf x}}_\alpha(t) =&\ \left(\left(\begin{matrix}
            A_\alpha & {\color{black} B_\alpha }  & C_\alpha \\
            {\color{black} A_\beta} & {\color{black} B_\beta}  &  {\color{black} C_\beta}  \\ 
            {\color{black} C_1} & {\color{black} D_1} {\color{black} }  & {\color{black} A_1+B_1} \\
             & & 
        \end{matrix}\right)\tilde{\bf x}_\alpha(t) + \left(\begin{matrix} D_\alpha \tilde{u}_\alpha(t) \\ {\color{black} 0} \\ {\color{black}-E_1R_1^{-1}E_1^\top\tilde{\zeta}(t) }\end{matrix}\right)\right)dt\\
        -d\tilde{\zeta}(t) = &\ \left(A_1^{\top}\tilde{\zeta}(t) + \left(\begin{matrix}
            -Q_1G_1 & {\color{black} -Q_1H_1} & Q_1({\color{black}I}-F_1)            
        \end{matrix}\right)\tilde{{\bf x}}_\alpha(t)\right)dt  - \left(\begin{matrix}
            Z_{\tilde{\zeta},\alpha} & {\color{black}  Z_{\tilde{\zeta},\beta}} 
        \end{matrix}\right)d{\bf W}(t), \\
        \tilde{{\bf x}}_\alpha(0) =&\ {\bf 0},\    
        {\color{black}\tilde{\zeta}(T):=\bar{Q}_1\left( \begin{matrix}
             -\bar{G}_1 & -\bar{H}_1 & (1-\bar{F}_1) 
        \end{matrix}\right){\bf x}_\alpha(T) .}
        \end{aligned}
        \right.
    \end{equation*}%
	The optimality of $u_{\alpha}^{\mathcal{N}}$ yields the   first-order condition: 
 \begin{align}\label{EulerC}
	0=&\ \frac{d}{d\tau}\bigg|_{\tau=0}J_{\alpha}(u_{\alpha}^{\mathcal{N}}+\tau \tilde{u}_{\alpha};u_{\beta}) \nonumber\\
	=&\ 2\mathbb{E}\bigg[\int_{0}^{T} \Big\{ \Big( \tilde{x}_{\alpha}(t)-{\color{black}G_{\alpha}\tilde{x}_{\beta}(t)}-F_{\alpha}\tilde{z}(t)\Big)^{\top}  Q_{\alpha}\Big(x^{\mathcal{N}}_{\alpha}(t)-F_{\alpha}z^{\mathcal{N}}(t)-G_{\alpha}x^{\mathcal{N}}_{\beta}(t)-M_{\alpha}\Big) \nonumber\\&\quad+\tilde{u}_{\alpha}^\top(t) R_{\alpha}u^{\mathcal{N}}_{\alpha}(t) \Big\} dt \bigg] + 2\bar{J}_\alpha,
\end{align}%
where  
    \begin{align*}
        \bar{J}_\alpha :=&\   \mathbb{E}\bigg[ \Big( \tilde{x}_{\alpha}(T)-{\color{black}\bar{G}_{\alpha}\tilde{x}_{\beta}(T)}-\bar{F}_{\alpha}\tilde{z}(T)\Big)^{\top}  \bar{Q}_{\alpha}\Big( x^{\mathcal{N}}_{\alpha}(T)-\bar{F}_{\alpha}z^{\mathcal{N}}(T)  -\bar{G}_{\alpha}x^{\mathcal{N}}_{\beta}(T)-\bar{M}_{\alpha}\Big)\bigg] .
    \end{align*}%
	Using  \eqref{NumericalNG} and  applying It\^o's lemma on $(p_{\alpha}^{\mathcal{N}})^{\top}\tilde{x}_{\alpha}+(r^{\mathcal{N}}_{\alpha})^{\top}\tilde{z}-(s^{\mathcal{N}}_{\alpha})^{\top}\tilde{\zeta}$, we arrive at  
         \begin{align}
    \label{eq:thm:nash:pf:1}
         \bar{J}_\alpha   =&\ \mathbb{E}\bigg[ \int_0^T \bigg( (p_\alpha^{\mathcal{N}})^\top(t)D_\alpha \tilde{u}_\alpha(t) -(\tilde{x}_\alpha(t)-{\color{black}G_{\alpha}\tilde{x}_{\beta}(t)}-F_\alpha\tilde{z}(t))^\top\nonumber \\
            &\ \quad Q_\alpha(x_\alpha^\mathcal{N}(t)-F_\alpha z^\mathcal{N}(t)-G_\alpha x_\beta^\mathcal{N}(t)-M_\alpha)  \bigg)dt\bigg].
        \end{align}%
   Combining \eqref{EulerC} and \eqref{eq:thm:nash:pf:1}, we obtain  
	\begin{equation}
	    \label{EulerC:final}	0=\mathbb{E}\bigg[\int_{0}^{T} \tilde{u}^{\top}_{\alpha}(t)\Big(R_{\alpha}u_{\alpha}^{\mathcal{N}}(t)+D_{\alpha}^{\top}p_{\alpha}^{\mathcal{N}}(t)\Big) dt\bigg].
	\end{equation}%
Since $\tilde{u}_{\alpha}$ is arbitrary, we conclude that $u_\alpha^\mathcal{N}$ is given by \eqref{Nashoptimalcontrol}. \hspace*{\fill} \qed
	\color{black}

\subsection{Proof of Theorem \ref{thm_Pareto}}
  Under the Pareto game, each leader considers the state dynamics of the other leader endogenously.  Hence, we perturb the pairs of controls of the leaders simultaneously: for $t\in[0,T]$ and $\tau>0$, consider the pair of controls $(u_\alpha^\tau(t),u_\beta^\tau(t)):=(u_{\alpha}^{\mathcal{P}}(t),u_{\beta}^{\mathcal{P}}(t))+\tau(\tilde{u}_{\alpha}(t),\tilde{u}_{\beta}(t))$,  where $\tilde{u}_{\alpha}$ and $\tilde{u}_{\beta}$ are arbitrary square integrable processes adapted to $\mathcal{F}^{\alpha}\vee\mathcal{F}^{\beta}$. Let $x_\alpha^\tau$, $x_\beta^\tau$ and $z^\tau$ be respectively the state dynamics of the $\alpha$-leader, the $\beta$-leader and the mean field term under the control $(u^\tau_\alpha,u^\tau_\beta)$. By the linearity of the state dynamics, it is easy to see that, for $\gamma=\alpha$ and $\beta$, $x_\gamma^\tau= x_{\gamma}^\mathcal{P}+\tau\tilde{x}_{\gamma}$ and $z^\tau=z^{\mathcal{P}}+\tau\tilde{z}$, where   ${\bf \tilde{x}} = (\Tilde{x}_\alpha,\tilde{x}_\beta,\tilde{z})$ satisfies  
    \begin{equation*}
    \left\{ \begin{aligned}
        d\tilde{{\bf x}}(t) =&\ \left(\left(\begin{matrix}
            A_\alpha & B_\alpha & C_\alpha \\
            A_\beta & B_\beta & C_\beta \\
            C_1 & D_1 & A_1+B_1
        \end{matrix}\right)\tilde{\bf x}+\left(\begin{matrix}
            D_\alpha\tilde{u}_\alpha(t) \\ D_\beta\tilde{u}_\beta(t) \\
          -  E_1R_1^{-1}E_1^\top\tilde{\zeta}(t)
        \end{matrix}\right) \right) dt\\
        -d\tilde{\zeta}(t) = &\ \left(A_1^{\top}\tilde{\zeta}(t) + \left(\begin{matrix}
            -Q_1G_1 & -Q_1H_1 & Q_1(1-F_1)            
        \end{matrix}\right)\tilde{{\bf x}}(t)\right)dt  - \left(\begin{matrix}
            Z_{\tilde{\zeta},\alpha} &  Z_{\tilde{\zeta},\beta}  
        \end{matrix}\right)d{\bf W}(t), \\
        \tilde{{\bf x}}(0) =&\ {\bf 0}, \quad   \tilde{\zeta}(T)= \bar{Q}_1\Big((I-\bar{F}_1)\tilde{z}(T)-\bar{G}_1\tilde{x}_{\alpha}(T)  -\bar{H}_1\tilde{x}_{\beta}(T)\Big).
        \end{aligned}
        \right.
    \end{equation*}%
	We consider the first-order condition   
	\begin{align}\label{EulerC_P}
	0=&\ \frac{d}{d\tau}\bigg|_{\tau=0}\left(J_{\alpha}(u^\tau_\alpha;u^\tau_\beta) + J_\beta(u^\tau_\beta;u^\tau_\alpha) \right) \nonumber \\
	=&\ 2\mathbb{E}\bigg[\int_{0}^{T}\Big\{\Big( \tilde{x}_{\alpha}(t)-F_{\alpha}\tilde{z}(t)-   G_\alpha \tilde{x}_\beta(t) \Big)^{\top} Q_{\alpha}\Big(x^{\mathcal{P}}_{\alpha}(t)-F_{\alpha}z^{\mathcal{P}}(t)-G_{\alpha}x^{\mathcal{P}}_{\beta}(t)-M_{\alpha}\Big) \nonumber\\
 &\quad +\tilde{u}_{\alpha}(t)^{\top}R_{\alpha}u_{\alpha}^{\mathcal{P}}(t)  +\Big( \tilde{x}_{\beta}(t)-F_{\beta}\tilde{z}(t)  - G_\beta \tilde{x}_\alpha(t)\Big)^{\top} Q_{\beta}\Big(x^{\mathcal{P}}_{\beta}(t)-F_{\beta}z^{\mathcal{P}}(t) \nonumber\\
 &\quad-G_{\beta}x^{\mathcal{P}}_{\alpha}(t)-M_{\beta}\Big)+ \tilde{u}_{\beta}(t)^{\top} R_{\beta}u_{\beta}^{\mathcal{P}}(t)\Big\} dt\bigg] + 2\bar{J}_{\alpha,\beta},
	\end{align}%
where  
    \begin{align*}
        \bar{J}_{\alpha,\beta} &:= \mathbb{E}\bigg[ \Big( \tilde{x}_{\alpha}(T)-\bar{F}_{\alpha}\tilde{z}(T)-\bar{G}_\alpha \tilde{x}_\beta(T)\Big)^{\top}  \bar{Q}_{\alpha}\Big( x^{\mathcal{P}}_{\alpha}(T)    -\bar{F}_{\alpha}z^{\mathcal{P}}(T)-\bar{G}_{\alpha}x^{\mathcal{P}}_{\beta}(T)-\bar{M}_{\alpha}\Big) \nonumber\\
        &\quad + \Big( \tilde{x}_{\beta}(T) -\bar{F}_{\beta}\tilde{z}(T)    -\bar{G}_\beta \tilde{x}_\alpha(T)\Big)^{\top}  \bar{Q}_{\beta}\Big( x^{\mathcal{P}}_{\beta}(T) -\bar{F}_{\beta}z^{\mathcal{P}}(T)-\bar{G}_{\alpha}x^{\mathcal{P}}_{\beta}(T)-\bar{M}_{\beta}\Big)\bigg].
    \end{align*}%
Using  \eqref{NumericalPG} and  applying Ito's lemma on $(p_{\alpha}^{\mathcal{P}})^{\top}\tilde{x}_{\alpha}+ (p_{\beta}^{\mathcal{P}})^{\top}\tilde{x}_{\beta} +(r^{\mathcal{P}}_{\alpha}+r^\mathcal{P}_\beta)^{\top}\tilde{z}-(s^{\mathcal{P}}_{\alpha}+s^\mathcal{P}_\beta)^{\top}\tilde{\zeta}$, we arrive at 
   \begin{align}
   \label{eq:thm:pareto:pf:1}
          \bar{J}_{\alpha,\beta}  =&\ \mathbb{E}\bigg[ \int_0^T \bigg( (p_\alpha^{\mathcal{P}})^\top(t)D_\alpha \tilde{u}_\alpha(t) + (p_\beta^{\mathcal{P}})^\top(t)D_\beta \tilde{u}_\beta(t) \nonumber\\
            &\ \quad -\Big(\tilde{x}_\alpha(t)-F_\alpha\tilde{z}(t)-G_\alpha\tilde{x}_\beta\Big)^\top Q_\alpha\Big(x_\alpha^\mathcal{P}(t) - G_\alpha x_\beta^\mathcal{P}(t) -M_\alpha\Big)  \nonumber \\
            &\ \quad -\Big( \tilde{x}_{\beta}(t)-F_{\beta}\tilde{z}(t)  - G_\beta \tilde{x}_\alpha(t)\Big)^{\top} Q_{\beta}\Big(x^{\mathcal{P}}_{\beta}(t)   -F_{\beta}z^{\mathcal{P}}(t)  -G_{\beta}x^{\mathcal{P}}_{\alpha}(t)-M_{\beta}\Big)  \bigg)dt\bigg].
        \end{align}%
  Combining \eqref{EulerC_P} and \eqref{eq:thm:pareto:pf:1}, we obtain 
	\begin{multline*}
	0=\mathbb{E}\bigg[\int_{0}^{T}\Big\{\tilde{u}_{\alpha}(t)^{\top}\Big( R_{\alpha}u_{\alpha}^{\mathcal{P}}(t)+D_{\alpha}^{\top}p_{\alpha}^{\mathcal{P}}(t)\Big)+\tilde{u}_{\beta}(t)^{\top}\Big( R_{\beta}u_{\beta}^{\mathcal{P}}(t)+D_{\beta}^{\top}p_{\beta}^{\mathcal{P}}(t)\Big)\Big\} dt \bigg].
	\end{multline*}
	Since $\tilde{u}_{\alpha}$ and $\tilde{u}_{\beta}$ are arbitrary, we deduce that $u^\mathcal{P}_\alpha$ and $u^\mathcal{P}_\beta$ are given by \eqref{Paretooptimalcontrol}. 
 \hspace*{\fill} \qed

     \subsection{Relationship of Linear FBSDEs and Riccati Equations}
    \label{sec:app:FBSDE:riccati}
    We discuss the relationship of well-posedness of the FBSDEs \eqref{NumericalNG}, \eqref{NumericalPG}, and the associated Riccati equations \eqref{dynGamma} and \eqref{dynGamma:agumented}. Since the treatment for both games are identical, we shall drop the superscript $\mathcal{I}$ in the subsequent proofs. 
\begin{lemma}\label{unique:Riccati}
    There exists at most one solution to the Riccati equations \eqref{dynGamma} and \eqref{dynGamma:agumented}. 
\end{lemma}

\begin{proof}
    Assume that \eqref{dynGamma} admits two distinct solutions $\mathbf{\Gamma}_1$ and $\mathbf{\Gamma}_2$. Since $[0,T]$ is compact, we have $\rho_T(\mathbf{\Gamma}_i)<\infty$, $i=1,2$. Let $\tilde{\mathbf{\Gamma}}:=\mathbf{\Gamma}_1-\mathbf{\Gamma}_2$, which satisfies 
	\begin{equation}\left\{
	\begin{aligned}
	-\frac{d \tilde{\mathbf{\Gamma}}(t)}{dt} =&\ \tilde{\mathbf{\Gamma}}(t){\bf  A}+{\bf  C}\tilde{\mathbf{\Gamma}}(t)-{\bf \Gamma}_1(t){\bf  B}\tilde{\mathbf{\Gamma}}(t)-\tilde{\mathbf{\Gamma}}(t){\bf  B}{\bf \Gamma}_2(t);\\
  \tilde{\mathbf{\Gamma}}(T)=&\ {\bf 0}.
	\end{aligned}\right.
        \label{eqn:tildeGamma}
	\end{equation}%
    Consider a small interval $[T-\Delta, T]$, where 
    \begin{equation*}
        \Delta:= \min\left\{\frac1{2(\rho({\bf  A})+\rho({\bf  C})+\rho({\bf  B})(\rho_T(\mathbf{\Gamma}_1)+\rho_T(\mathbf{\Gamma}_2)))} ,T\right\} .
    \end{equation*}%
    By \eqref{eqn:tildeGamma}, for $t\in[T-\Delta,T]$, we have  
    \begin{equation*}
        \tilde{\mathbf{\Gamma}}(t)=\int_t^T \big(\tilde{\mathbf{\Gamma}}(s){\bf  A}+{\bf  C}\tilde{\mathbf{\Gamma}}(s)-{\bf \Gamma}_1(s){\bf  B}\tilde{\mathbf{\Gamma}}(s)-\tilde{\mathbf{\Gamma}}(s){\bf  B}{\bf \Gamma}_2(s)\big) ds. 
    \end{equation*}
    Hence,   
    \begin{eqnarray*}
        \sup\limits_{t\in [T-\Delta , T]} \rho(\tilde{\mathbf{\Gamma}}(t)) &\le& \Delta (\rho({\bf  A})+\rho({\bf  C})+\rho({\bf  B})(\rho_T(\mathbf{\Gamma}_1)+\rho_T(\mathbf{\Gamma}_2)))  \sup\limits_{t\in [T-\Delta , T]} \rho(\tilde{\mathbf{\Gamma}}(t)) \\ 
        &\le& \frac12 \sup\limits_{t\in [T-\Delta , T]} \rho(\tilde{\mathbf{\Gamma}}(t)),
    \end{eqnarray*}%
   which implies $\tilde{\mathbf{\Gamma}} \equiv {\bf 0}$ on $[T-\Delta, T]$. Since $\Delta$ is time-independent, by induction, we can deduce that $\tilde{\mathbf{\Gamma}}  \equiv {\bf 0}$ on $[0, T]$, and thus $\mathbf{\Gamma}_1\equiv \mathbf{\Gamma}_2$.
\end{proof}

\begin{lemma}\label{unique:FBSDE}
    System \eqref{NumericalNG} (resp. \eqref{NumericalPG}) admits a unique solution if and only if \eqref{dynGamma} does for $\mathcal{I}=\mathcal{N}$ (resp. $\mathcal{I}=\mathcal{P}$).
\end{lemma}

\begin{proof}
By Proposition 3.2 of \cite{BSYY}, the Riccati equation \eqref{dynGamma} for $\mathcal{I}=\mathcal{N}$ (resp.~$\mathcal{I}=\mathcal{P}$) admits a unique solution if  System \eqref{NumericalNG} (resp.~\eqref{NumericalPG}) admits a unique solution. Thus, we only show the only if part. Denote by ${\bf \Gamma}$ the solution of \eqref{dynGamma}. Since \eqref{NumericalAll2g} is a linear backward ODE and ${\bf \Gamma}$ is exogeneously given,  \eqref{NumericalAll2g} admits a unique solution ${\bf g}$. It is easy to check $({\bf x}, {\bf p})$ with ${\bf p}(t)={\bf \Gamma}(t){\bf x}(t)+{\bf g}(t)$   is a solution of System \eqref{NumericalNG} (or \eqref{NumericalPG}). Next, we show that there exists no other solution of System \eqref{NumericalNG} (or \eqref{NumericalPG}).

    To proceed, we first show that ${\bf p}$ must be affine in ${\bf x}$. Assume the contrary that System \eqref{NumericalNG} (or \eqref{NumericalPG}) admits a non-affine solution $({\bf x}^{\prime}, {\bf p}^{\prime})$. Denote by $\tilde{{\bf x}}={\bf x}-{\bf x}^{\prime}$, $\tilde{{\bf p}}={\bf p}-{\bf p}^{\prime}$. Then $(\tilde{{\bf x}},\tilde{{\bf p}})$ satisfies  
        \begin{equation}\label{eqn:tildexp}
        \left\{\begin{aligned}
        d\tilde{{\bf x}}(t)&=\Big({\bf  A}\tilde{{\bf x}}(t)-{\bf  B}\tilde{{\bf p}}(t)\Big)dt;\\
        -d\tilde{{\bf p}}(t)&=\Big({\bf  C}\tilde{{\bf p}}(t)+{\bf  D}\tilde{{\bf x}}(t)\Big)dt;\\
        \tilde{{\bf x}}(0)&={\bf  0},\quad \tilde{{\bf x}}(T)={\bf  E}\tilde{{\bf x}}(T).
        \end{aligned}\right.
        \end{equation}%
        By \eqref{eqn:tildexp}, $\tilde{{\bf p}}$ only depends on $\tilde{{\bf x}}$. However, since ${\bf p}'$ is not affine in ${\bf x}'$, $\tilde{{\bf p}}(t)={\bf \Gamma}(t){\bf x}(t)+{\bf g}(t)-{\bf p}'(t)$ depends on both ${\bf x}$ and ${\bf x}^{\prime}$, which leads to a contradiction. 

        Next, we show that affine solution of System \eqref{NumericalNG} (or \eqref{NumericalPG}) is unique. Assume the contrary that the FBSDE system admits another affine solution $({\bf x}^{\prime}, {\bf p}^{\prime})$, where ${\bf p}^{\prime}(t)={\bf \Gamma}^{\prime}(t){\bf x}^{\prime}(t)+{\bf g}^{\prime}(t)$. If ${\bf \Gamma}^{\prime} \not\equiv {\bf \Gamma}$, then ${\bf \Gamma}^{\prime}$ is another solution of \eqref{dynGamma}, which contradicts the uniqueness of \eqref{dynGamma}. If ${\bf \Gamma}^{\prime} \equiv {\bf \Gamma}$, we must have ${\bf g}^{\prime} \equiv {\bf g}$ by the uniqueness of \eqref{NumericalAll2g}. In other words, both  ${\bf x}$ and ${\bf x}^{\prime}$ are solutions of the following SDE
    \begin{equation}\label{eqn:SDEx}
    \left\{\begin{aligned}
    d{\bf  x}(t)&=\Big(({\bf  A}-{\bf  B}(t){\bf \Gamma}(t)){\bf x}(t)-{\bf  B}(t){\bf g}(t)\Big)dt+{\bf  \Sigma} d{\bf  W}(t);\\
    {\bf  x}(0)&={\bf  x}_0.
    \end{aligned}\right.
    \end{equation}
    Since \eqref{eqn:SDEx} is linear, the solution of it must be unique. Hence, we have $({\bf x},{\bf p})\equiv ({\bf x}^{\prime},{\bf p}')$  as desired. 

    {\color{black}Further discussions on the relationship between the linear FBSDE system and its associated Riccati equation can be found in  Section 5 in \cite{yong1999linear}.}

\end{proof}

\section{Proofs in Section \ref{sec:compare} }
\label{sec:app:compare}
This section is devoted to providing the proofs in Section \ref{sec:compare}. In the sequel, we shall write ${\bf 0}_{m,n}$ and ${\bf e}_{m,n}$ to denote the $m\times n$ matrices with all entries being 0 and 1, respectively.


\subsection{Proof of Theorem \ref{thm:cost}}
\label{app:cost}
     Since the proofs for the Nash and the Pareto game are identical, for notational convenience, we shall drop the superscript $\mathcal{I}$ for distinguishing the two games. 
    \begin{enumerate}
        \item By Theorem \ref{wellpose:dominate}, we have  $\zeta(t)=\mathbf{i}_3^{\top}({\bf \Gamma}(t) \mathbf{x}(t)+{\bf g}(t))$.  Let $\hat{x}_1(t):=x_1(t)-z(t)$ and $\hat{\xi}(t):=\xi(t)-\zeta(t)$. Then $\hat{x}_1$ and $\hat{\xi}$ solve the following system:  
               \begin{equation}\label{hatxform}
            \left\{\begin{aligned}
            d\hat{x}_1(t)&=\Big(A_1\hat{x}_1(t)-E_1R_1^{-1}E_1^{\top}\hat{\xi}(t)\Big)dt+\sigma_1dW^1(t),\\
            -d\hat{\xi}(t)&=\Big(A_1^{\top}\hat{\xi}(t)+Q_1 \hat{x}_1(t)\Big)dt-Z_{\hat{\xi},1}(t)dW^1(t) - Z_{\hat{\xi},\alpha}(t)dW^\alpha(t) - Z_{\hat{\xi},\beta}(t)dW^\beta(t),\\
             \hat{x}_1(0)&=\eta_1-\mathbb{E}[\eta_1], \quad  \hat{\xi}(T)=\bar{Q}_1 \hat{x}_1(T).
            	\end{aligned}\right.
            \end{equation}%
        By applying It\^o's lemma, it is easy to see that $\hat{\xi}(t)=\Gamma_0(t) \hat{x}_1(t)$, where $\Gamma_0$ solves the Riccati equation \eqref{Gamma0}. Using this and Theorem \ref{C3_follower}, we see that  
            \begin{equation*}
                u_1(t) = -R_1^{-1}E_1^\top \xi(t) = -R_1^{-1}E_1^\top( \Gamma_0(t)\hat{x}_1(t) + \zeta(t)),
            \end{equation*}%
        which leads to   \eqref{optimalfollowerform}. 

        \item By Fubini's theorem and the tower property of conditional expectations, we have $J_1 = J_{11} + J_{12} + J_{13}$, where 
            \begin{align*}
                J_{11} :=&\  \int_0^T  \mathbb{E}\big[|\hat{x}_1(t)|_{Q_1}^2\big] dt +  \mathbb{E}\left[ | \hat{x}(T)|^2_{\bar{Q}_1} \right], \\
                J_{12} :=&\  \int_0^T \mathbb{E}\left[ \left|u_1(t) - \mathbb{E}^{\mathcal{F}^{\alpha}_t\vee\mathcal{F}^{\beta}_t}\left[u_1(t)\right] \right|^2_{R_1}  \right] dt ,\\
                J_{13} :=&\  \int_0^T \bigg\{  \mathbb{E}\bigg[ \Big|(I-F_1)z(t)    - G_1x_\alpha(t)-M_1  -H_1x_\beta(t) \Big|^2_{Q_1} \bigg] +\Big|\mathbb{E}^{\mathcal{F}^\alpha_t\vee \mathcal{F}^\beta_t}[u_1(t)]\Big|^2_{R_1}   \bigg\}dt  \\ &\quad   + \mathbb{E}\bigg[\Big|(1-\bar{F}_1)z(T) -\bar{G}_1x_\alpha(T)-\bar{H}_1x_\beta(T)-\bar{M}_1 \Big|^2_{\bar{Q}_1} \bigg].
            \end{align*}%
        Using the fact that $\hat{\xi}(t)=\Gamma_0(t)\hat{x}_1(t)$ and solving \eqref{hatxform}, we obtain  
            \begin{equation}
            \label{eq:xhat:integral}
                \begin{aligned}
                    \hat{x}_1(t)=&\ \int_0^t e^{\int_s^t \bar{A}_1(u)du}\sigma_1 dW^1(s) + e^{\int_0^t\bar{A}_1(s)ds}(\eta_1-\mathbb{E}[\eta_1]),
                \end{aligned}
            \end{equation}%
        where $\bar{A}_1$   is defined in  \eqref{Gamma0Lambda0A}. Hence, we have  
          \begin{equation}
                \label{eq:lyapunov:pf:2}
            \int_0^T \mathbb{E} \left[ \left| \hat{x}(t) \right|^2_{Q_1}     \right] dt  = \int_0^T \mathrm{tr}\left(Q_1 (\Sigma_0(t)+P_0(t))\right)dt,
        \end{equation}%
    and  
    \begin{equation}
    \label{eq:lyaounov:pf:2.5}
    \mathbb{E}\left[ | \hat{x}(T)|^2_{\bar{Q}_1} \right]=\mathrm{tr}\left(\bar{Q}_1 (\Sigma_0(T)+P_0(T))\right),
    \end{equation}%
    where $\Sigma_0,P_0$ are defined in \eqref{Gamma0Lambda0A}.  Using   \eqref{optimalfollowerform} and \eqref{eq:xhat:integral}, we also have   
        \begin{align}
            \label{eq:lyapunov:pf:3}
        J_{12} =&\  \int_0^T \mathrm{tr}\left(\Gamma^{\top}_0(t) E_1R^{-1}_1E_1^{\top} \Gamma_0(t) (\Sigma_0+P_0)(t)\right)dt \\
           J_{13}  =&\ \mathbb{E}\bigg[\int_0^T \bigg\{ \Big\langle  {\bf x}(t),({\bf Q}+{\bf R}(t)){\bf x}(t) - 2{\bf q}Q_1M_1   + 2{\bf \Gamma}^\top(t){\bf i}_3 E_1R_1^{-1}E_1^\top {\bf i}_3^\top {\bf g}(t) \Big\rangle  \nonumber \\
            &\quad  + M_1^\top Q_1M_1+ {\bf g}^\top(t){\bf i}_3E_1R_1^{-1}E_1^\top {\bf i}_3^\top {\bf g}(t) \bigg\}dt + \Big\langle {\bf x}(T), {\bf \bar{Q}} {\bf x}(T) -2{\bf \bar{q}}\bar{Q}_1\bar{M}_1 \Big\rangle   \bigg]\nonumber \\
            &\quad + \bar{M}_1^\top \bar{Q}_1\bar{M}_1,
                          \label{eq:lyapunov:pf:4}
        \end{align}%
    By applying It\^o's lemma on ${\bf x}^\top({\bf K}{\bf x}+2{\bf k}) +l$ and using  \eqref{eq:K} - \eqref{eq:l},  we see that \eqref{eq:lyapunov:pf:4} coincides with  
        \begin{equation}
                    \label{eq:lyapunov:pf:5}
            \begin{aligned}
               &\ \mathbb{E}[{\bf x}^\top(0){\bf K}(0){\bf x}(0)] + \int_0^T \text{tr}(\langle {\bf \Sigma}, {\bf K}(t){\bf \Sigma}\rangle dt + 2\mathbb{E}[{\bf x}^\top(0){\bf k}(0)] +l(0).
            \end{aligned}
        \end{equation}%
   By summing \eqref{eq:lyapunov:pf:2}-\eqref{eq:lyapunov:pf:3} and \eqref{eq:lyapunov:pf:5}, we arrive at \eqref{eq:lem:cost}. 
    
    \end{enumerate}

 \subsection{Proof of Theorem \ref{thm:compare}}

    Under Model \ref{degeneratedmodel}, {\color{black}by defining}  $\tilde{x}_\beta := k x_\beta$, $\tilde{v}_\beta := kv_\beta$ and $\tilde{\eta}_\beta := k\eta_\beta$, we see that the states $(x_\alpha,\tilde{x}_\beta,x_1)$ and the associated cost functionals belong to Model \ref{degeneratedmodel} with $k=1$. Hence, without loss of generality, it suffices to consider $k=1$.
With a slight abuse of notations, we continue to use $x_{\beta}$ to denote the state dynamics of the $\beta$-leader instead of $\tilde{x}_{\beta}$. In the sequel, for any matrix ${\bf M}$, we denote by ${\bf M}^{{\bf b}_0}$   the corresponding matrix under Model \ref{degeneratedmodel}. To proceed, we need the following auxiliary comparison results for ODEs.

 \begin{lemma}
 \label{lem:strict:ineq}
    Let $w(t)$, $t\in[0,T]$, be the solution of the differential equation $w'(t) = f(t,w(t))$ with $w(0)=0$,
    where $f : [0,T]\times \mathbb{R}\to \mathbb{R}$ is not identical to zero, and  is locally Lipschitz  in the spatial variable  uniformly in $[0,T]$: for any $n>0$, there exists $L_n>0$ such that for any $t\in[0,t]$ and $w_1,w_2\in \mathbb{R}$ with $|w_1|,|w_2|\leq n$, it holds that
        \begin{equation*}
            |f(t,w_1)-f(t,w_2)| \leq L_n|w_1-w_2|.
        \end{equation*}
    Let $q(t)$ be a differentiable function in $[0,T]$ such that $q(0)\geq 0$,   $q'(t)\geq f(t,q(t))$ for $t\in[0,T]$, {\color{black}and $\max\{q'(t)-f(t,q(t)),q(0)-w(0)\}>0$ for $t\in(0,T]$}. Then $q(t)>w(t)$ for $t\in(0,T]$.
    \end{lemma}
\begin{proof}
    Let $h(t):= q(t)-w(t)$. {\color{black} If $h(0)=q(0)-w(0)>0$ and $h'(t)=q'(t)-w'(t)\geq 0$, it is clear that $h(t)>0$ for all $t\in[0,T]$. Next, we consider the case $h(0)=0$, which implies that $h'(t)>0$ for all $t\in(0,T]$.  } Assume the contrary that there exists $t_0\in(0,T]$ such that $h(t_0)=0$. Notice that $h'(t_0) > f(t_0,q(t_0))-f(t_0,w(t_0)) =0$, which   implies that $h$ cannot be identically zero in $[0,t_0]$. In addition, there exists $t_1\in(0,t_0)$ with $h(t_1)<0$. By continuity of $h$, there exists $t_2\in [0,t_1)$ such that $h(t_2)=0$ and $h(t)\leq 0$ for $t\in [t_2,t_1]$. Hence, for $t\in[t_2,t_1]$,
        \begin{equation*}
            h'(t) > f(t,q(t))-f(t,w(t))  = L_mh(t),
        \end{equation*}
    where $m:= \sup_{t\in[0,T]} \{|w(t)|,|q(t)| \} <\infty$. Therefore $h(t_1) > h(t_2)e^{L_m(t_1-t_2)}=0$, which contradicts with the fact that $h(t_1)<0$.
\end{proof}

\begin{lemma}
    \label{lem:ode:compare}
   For $i=1,2$, let $w_i(t)$, $t\in[0,T]$  be the solution of the differential equation
         \begin{equation*}
            -w_i'(t) = -w_i^2(t) + a_i(t) w_i(t) + b_i(t)   , \ w_i(T)=w_{iT},
        \end{equation*}
    where $a_i,b_i$ are continuous in $[0,T]$. 
    {\color{black}Assume that $w_{2T}\geq w_{1T}\geq 0$, $a_2(t)\geq  a_1(t)$ and $b_2(t)\geq b_1(t)\geq 0$ for $t\in[0,T]$, and $\max\{a_2(t)-a_1(t),b_2(t)-b_1(t),w_{2T}-w_{1T} \}>0$ for $t\in[0,T)$.}
    Then  $w_2(t)> w_1(t)> 0$ for $t\in[0,T)$.
\end{lemma}
    \begin{proof}
        To see that $w_1(t)>0$ for $t\in[0,T)$, consider the solution of the equation 
        \begin{equation*}
            -w_0'(t) = -w_0^2(t) + a_1(t)w_0(t), \ w_0(T)=0.
        \end{equation*}
        Since $b_1(t)>0$ for $t\in[0,T)$ and $b_1(T)\geq 0$, we infer by Lemma \ref{lem:strict:ineq} that $w_1(t)>w_0(t)=0$ for $t\in[0,T)$. Now consider $\tilde{w}_1:=w_1-w_2$, which satisfies  
       \begin{equation*}
       \left\{\begin{aligned}
            -\tilde{w}_1'(t) =&\ -\tilde{w}_1^2(t)-2w_2(t)\tilde{w}_1(t) + a_1(t) \tilde{w}_1(t)  \\ &\ + (a_1(t)-a_2(t))w_2(t)+ b_1(t)- b_2(t)  , \\
            \tilde{w}_1(T)=&\ \tilde{w}_{1T} := w_{1T}-w_{2T}\geq 0. 
        \end{aligned}\right.
        \end{equation*}%
    Let  $\tilde{w}_2$ be the solution of the following equation  
    \begin{equation}\label{z2}
            -\tilde{w}_2'(t) = -\tilde{w}_2^2(t)-2w_2(t)\tilde{w}_2(t) + a_1(t) \tilde{w}_2(t), \ \tilde{w}_2(T)=0.
        \end{equation}%
    By the unique existence of solution of \eqref{z2}, we infer that $\tilde{w}_2(t)=0$ for all $t\in[0,T]$. Since $w_2(t)> 0$ for all $t\in[0,T)$, we have $(a_1(t)-a_2(t))w_2(t)+ b_1(t)- b_2(t){\color{black} \leq }0$ for all $t\in[0,T)$. Hence, by Lemma \ref{lem:strict:ineq}, we conclude that $0=\tilde{w}_2(t)> \tilde{w}_1(t)=w_1(t)-w_2(t)$ for $t\in[0,T)$.
    \end{proof}

\subsubsection{Analysis of the Riccati Equations}
We first consider the Riccati equations associated with the leader-game, which is decoupled from the {\color{black}representative} follower due to the absence of the mean field term in the leaders' dynamics and cost functionals. For $\mathcal{I}=\mathcal{N}$ and $\mathcal{P}$,
let  ${\bf x}^{\alpha\beta,\mathcal{I}}(t) := (x^\mathcal{I}_\alpha(t),x^\mathcal{I}_\beta(t))$, which satisfies \vspace{-0.6cm} {\small 
    \begin{align*}
        d{\bf x}^{\alpha\beta,\mathcal{I}}(t) =&\ ({\bf A}^{\alpha\beta} - {\color{black}{\bf B}^{\alpha\beta,\mathcal{I}}} {\bf \Gamma}^{\alpha\beta,\mathcal{I}}(t)){\bf x}^{\alpha\beta,\mathcal{I}}(t) dt + {\bf \Sigma}^{\alpha\beta}d{\bf W}(t),
    \end{align*}}%
where ${\bf \Sigma}^{\alpha\beta} := \text{diag}(\sigma_\alpha,\sigma_\beta)$, ${\bf B}^{\alpha\beta,\mathcal{P}} =D_{\alpha\beta}^2R_{\alpha\beta}^{-1} \text{diag}(1,1)$,   
    \begin{align*}
       & {\bf A}^{\alpha\beta} := \begin{pmatrix}
                A_{\alpha\beta}& B_{\alpha\beta} \\ B_{\alpha\beta}& A_{\alpha\beta}
        \end{pmatrix}, \  {\color{black}{\bf B}^{\alpha\beta,\mathcal{N}}:= \begin{pmatrix}
            {\bf B}^{\alpha\beta,\mathcal{P}} & {\bf 0}_{2,2}
        \end{pmatrix}},
    \end{align*}%
and ${\bf \Gamma}^{\alpha\beta,\mathcal{I}}(t)$, $t\in[0,T]$, is the solution of the Riccati equation 
    \begin{equation}\label{Gamma}
       \left\{ \begin{aligned}
        -\frac{d{\bf \Gamma}^{\alpha\beta,\mathcal{I}}(t)}{dt} =&\  {\bf \Gamma}^{\alpha\beta,\mathcal{I}}(t){\bf A}^{\alpha\beta} + {\bf C}^{\alpha\beta,\mathcal{I}} {\bf \Gamma}^{\alpha\beta,\mathcal{I}}(t)  - {\bf \Gamma}^{\alpha\beta,\mathcal{I}}(t){\bf B}^{\alpha\beta{\color{black},\mathcal{I}}} {\bf \Gamma}^{\alpha\beta,\mathcal{I}}(t) + {\bf D}^{\alpha\beta,\mathcal{I}},
        \\ {\bf \Gamma}^{\alpha\beta,\mathcal{I}}(T)=&\ {\bf E}^{\alpha\beta,\mathcal{I}},
        \end{aligned}\right.
    \end{equation}%
where 
    \begin{align*}
    &   {\color{black}  {\bf C}^{\alpha\beta,\mathcal{N}} := \left(\begin{matrix}
            A_{\alpha\beta} & 0 & B_{\alpha\beta} & 0 \\
            0 & A_{\alpha\beta} & 0 & B_{\alpha\beta} \\
            B_{\alpha\beta} & 0 & A_{\alpha\beta} & 0 \\
            0 & B_{\alpha\beta} & 0 & A_{\alpha\beta} \\
        \end{matrix}\right) }, \ {\bf C}^{\alpha\beta,{\color{black}\mathcal{P}}} :=  ({\bf A}^{\alpha\beta})^\top, \\
        & {\color{black} 
         {\bf D}^{\alpha\beta,\mathcal{N}} := Q_{\alpha\beta}\left(\begin{matrix}
         1 & -G_{\alpha\beta} \\
-G_{\alpha\beta} & 1 \\
{\color{black}-G_{\alpha\beta} } & {\color{black}  G_{\alpha\beta}^2 } \\  
G_{\alpha\beta}^2 &  -G_{\alpha\beta}   
         \end{matrix}
         \right)},\ {\bf D}^{\alpha\beta,\mathcal{P}} = {\bf I}^{\alpha\beta}{\bf D}^{\alpha\beta,\mathcal{N}},  \\     
       &  {\color{black} 
         {\bf E}^{\alpha\beta,\mathcal{N}} := \bar{Q}_{\alpha\beta}\left(\begin{matrix}
         1 & -\bar{G}_{\alpha\beta} \\
-\bar{G}_{\alpha\beta} & 1 \\
{\color{black}-\bar{G}_{\alpha\beta} } & {\color{black}  \bar{G}_{\alpha\beta}^2 } \\  
\bar{G}_{\alpha\beta}^2 &  -\bar{G}_{\alpha\beta}   
         \end{matrix}
         \right)}, \ {\bf E}^{\alpha\beta,\mathcal{P}} = {\bf I}^{\alpha\beta}{\bf E}^{\alpha\beta,\mathcal{N}}, \\
        & {\color{black}{\bf I}^{\alpha\beta} := \left(\begin{matrix}
             1 & 0 & 0 & 1 \\
             0 & 1 & 1 & 0
         \end{matrix}\right)}.
    \end{align*}%
Since Assumption \ref{ass:2:combined}  holds, we infer that the solution of the Riccati equation \eqref{Gamma} uniquely exists for both $\mathcal{I}=\mathcal{N}$,  $\mathcal{P}$. Write ${\bf \Gamma}^{\alpha\beta,\mathcal{N}}(t) = (\Gamma^{\alpha\beta,\mathcal{N}}_{ij}(t))_{{\color{black}1\leq i\leq 4, 1\leq j\leq 2}}$ {\color{black}and ${\bf \Gamma}^{\alpha\beta,\mathcal{P}}(t) = (\Gamma^{\alpha\beta,\mathcal{P}}_{ij}(t))_{{\color{black}1\leq i,j\leq 2}}$}. It is clear that ${\bf \Gamma}^{\alpha\beta,\mathcal{P}}(t)$ is symmetric, since ${\bf C}^{\alpha\beta,\mathcal{P}} = ({\bf A}^{\alpha\beta})^\top$, and ${\bf D}^{\alpha\beta,\mathcal{P}}$,  ${\bf E}^{\alpha\beta,\mathcal{P}}$ are symmetric. {\color{black} For both games, using a simple uniqueness argument, we have the following:}
    \begin{proposition}
    \label{symmetry:nash}
       For $\mathcal{I}=\mathcal{N}$ and $\mathcal{P}$, the solution ${\bf \Gamma}^{\alpha\beta,\mathcal{I}}(t)$, {\color{black}$t\in[0,T]$, of} the Riccati equation \eqref{Gamma}   satisfies:
       \begin{enumerate}
           \item  $\Gamma_{11}^{\alpha\beta,\mathcal{I}}(t) = \Gamma_{22}^{\alpha\beta,\mathcal{I}}(t)$ and  $\Gamma_{12}^{\alpha\beta,\mathcal{I}}(t) = \Gamma_{21}^{\alpha\beta,\mathcal{I}}(t)$;
           \item {\color{black}  $\Gamma_{31}^{\mathcal{N}}(t) =\Gamma_{42}^{\mathcal{N}}(t)$ and $\Gamma_{32}^{\mathcal{N}}(t) =\Gamma_{41}^{\mathcal{N}}(t)$. }
       \end{enumerate}
       
    \end{proposition}
\begin{proof}
{\color{black} It suffices to consider $\mathcal{I}=\mathcal{N}$.} For $t\in[0,T]$, let 
$$ {\color{black}  \tilde{{\bf \Gamma}}^{\alpha\beta,\mathcal{N}}(t) := \left( \begin{matrix}
     \Gamma^{\alpha\beta,\mathcal{N}}_{11}(t) &  \Gamma^{\alpha\beta,\mathcal{N}}_{12}(t) \\
      \Gamma^{\alpha\beta,\mathcal{N}}_{12}(t) &  \Gamma^{\alpha\beta,\mathcal{N}}_{11}(t) \\
       \Gamma^{\alpha\beta,\mathcal{N}}_{31}(t) &  \Gamma^{\alpha\beta,\mathcal{N}}_{32}(t) \\ 
        \Gamma^{\alpha\beta,\mathcal{N}}_{32}(t) &  \Gamma^{\alpha\beta,\mathcal{N}}_{31}(t)
\end{matrix} \right)} .$$
 It {\color{black} is easy to verify} that  $\Tilde{{\bf \Gamma}}^{\alpha\beta,\mathcal{N}}$ satisfies the Riccati equation \eqref{Gamma} {\color{black} for $\mathcal{I}=\mathcal{N}$}. Hence, by the uniqueness of solution, we have $\Tilde{{\bf \Gamma}}^{\alpha\beta,\mathcal{N}}\equiv {\bf \Gamma}^{\alpha\beta,\mathcal{N}}$ and the result follows. %
\end{proof}

\begin{proposition}
\label{pp:leader:compare}
    For $\mathcal{I}=\mathcal{N}$ and $\mathcal{P}$, let   $\Gamma^{\alpha\beta,\mathcal{I}}(t):=\Gamma^{\alpha\beta,\mathcal{I}}_{11}(t)+\Gamma^{\alpha\beta,\mathcal{I}}_{12}(t)$. {\color{black}Suppose $\max\{\vert B_{\alpha\beta} \vert, \vert G_{\alpha\beta} \vert,  \vert \bar{Q}_{\alpha\beta}\bar{G}_{\alpha\beta}  \vert\} $ $> 0$.} Then for any $t\in[0,T)$ {\color{black} and $G_{\alpha\beta},\bar{G}_{\alpha\beta}\in (-1,1)$}, we have $\Gamma^{\alpha\beta,\mathcal{I}}(t)> 0$, and in addition, 
\begin{enumerate}
    \item  $\Gamma^{\alpha\beta,\mathcal{N}}(t) < \Gamma^{\alpha\beta,\mathcal{P}}(t)$ if $B_{\alpha\beta}\geq 0$ {\color{black} and $G_{\alpha\beta},\bar{G}_{\alpha\beta} \in (-1,0]$};
    \item $\Gamma^{\alpha\beta,\mathcal{N}}(t) > \Gamma^{\alpha\beta,\mathcal{P}}(t)$ if $B_{\alpha\beta}\leq  0$ {\color{black} and $G_{\alpha\beta},\bar{G}_{\alpha\beta} \in [0,1)$}.
\end{enumerate}
\end{proposition}
\begin{proof}
   Notice that for $\mathcal{I}=\mathcal{N},\mathcal{P}$,  $\Gamma^{\alpha\beta,\mathcal{I}}(t)$ satisfies the following equation: 
    \begin{equation}
    \label{eq:Gamma:alpha:beta}
        \left\{\begin{aligned}
             -\frac{d\Gamma^{\alpha\beta,\mathcal{I}}(t)}{dt}=&\ 2(A_{\alpha\beta}+ B_{\alpha\beta})\Gamma^{\alpha\beta,\mathcal{I}}(t) - {\color{black}  B_{\alpha\beta} \tilde{\Gamma}^{\alpha\beta,\mathcal{N}}(t) \mathbbm{1}_{\mathcal{I}=\mathcal{N}}} 
               -\frac{D_{\alpha\beta}^2(\Gamma^{\alpha\beta,\mathcal{I}}(t))^2}{  R_{\alpha\beta}} \\
               &\ +  Q_{\alpha\beta}(1-G_{\alpha\beta})   
             {\color{black}-\mathbbm{1}_{\mathcal{I}=\mathcal{P}}Q_{\alpha\beta}G_{\alpha\beta}(1-G_{\alpha\beta})}, \\  \Gamma^{\alpha\beta,\mathcal{I}}(T)=&\ \bar{Q}_{\alpha\beta}(1-\bar{G}_{\alpha\beta}){\color{black}-\mathbbm{1}_{\mathcal{I}=\mathcal{P}}\bar{Q}_{\alpha\beta}\bar{G}_{\alpha\beta}(1-\bar{G}_{\alpha\beta})}.
       \end{aligned}  \right. 
    \end{equation}%
     {\color{black}
    Here, $\tilde{\Gamma}^{\alpha\beta,\mathcal{N}}(t):= \Gamma^{\alpha\beta,\mathcal{N}}(t) - \hat{\Gamma}^{\alpha\beta,\mathcal{N}}(t)$, and $\hat{\Gamma}^{\alpha\beta,\mathcal{N}}(t):=\Gamma^{\alpha\beta,\mathcal{N}}_{31}(t)+\Gamma^{\alpha\beta,\mathcal{N}}_{32}(t)$, which respectively satisfies 
\begin{align}
 \label{eq:Gamma:tilde:alpha:beta:N}
      &  \left\{\begin{aligned}
             -\frac{d\tilde{\Gamma}^{\alpha\beta,\mathcal{N}}(t)}{dt}=&\ \left(2A_{\alpha\beta} -\frac{D^2_{\alpha\beta} \Gamma^{\alpha\beta,\mathcal{N}}(t) }{R_{\alpha\beta}}\right)\tilde{\Gamma}^{\alpha\beta,\mathcal{N}}(t)  + Q_{\alpha\beta}(1-G_{\alpha\beta}^2),\\ \tilde{\Gamma}^{\alpha\beta,\mathcal{N}}(T)=&\  \bar{Q}_{\alpha\beta} (1-\bar{G}_{\alpha\beta}^2), 
       \end{aligned}  \right. \\
\label{eq:Gamma:hat:alpha:beta:N}
    &  \left\{\begin{aligned}
             -\frac{d\hat{\Gamma}^{\alpha\beta,\mathcal{N}}(t)}{dt}=&\ \left(2A_{\alpha\beta}+ 2B_{\alpha\beta} -\frac{D_{\alpha\beta}^2\Gamma^{\alpha\beta,\mathcal{N}}(t) }{R_{\alpha\beta}} \right)\hat{\Gamma}^{\alpha\beta,\mathcal{N}}(t)  \\
             &\ + B_{\alpha\beta} \tilde{\Gamma}^{\alpha\beta,\mathcal{N}}(t)     - Q_{\alpha\beta}G_{\alpha\beta}(1-G_{\alpha\beta}),\\ \hat{\Gamma}^{\alpha\beta,\mathcal{N}}(T)=&\ - \bar{Q}_{\alpha\beta}\bar{G}_{\alpha\beta}(1-\bar{G}_{\alpha\beta}). 
       \end{aligned}  \right.
    \end{align}%
 For $\mathcal{I}=\mathcal{P}$, using the fact that $G_{\alpha\beta},\bar{G}_{\alpha\beta}\in(-1,1)$, we have $\Gamma^{\alpha\beta,\mathcal{P}}(t)>0$ for $t\in[0,T)$, thanks to Lemma \ref{lem:ode:compare}. For $\mathcal{I}=\mathcal{N}$, \eqref{eq:Gamma:tilde:alpha:beta:N} and Lemma \ref{lem:ode:compare} imply that  $\tilde{\Gamma}^{\alpha\beta,\mathcal{N}}(t)>0$ for $t\in[0,T]$. We now consider the two cases specified in the proposition. 
    \begin{enumerate}[label = \arabic*)]
        \item If $B_{\alpha\beta} \geq 0$,  $G_{\alpha\beta},\bar{G}_{\alpha\beta}\in (-1,0]$, and at least one of these coefficients are non-zero, then for any $t\in[0,T]$, $B_{\alpha\beta}\tilde{\Gamma}^{\alpha,\beta,\mathcal{N}}(t)-Q_{\alpha\beta}G_{\alpha\beta}(1-G_{\alpha\beta}) \geq 0$, and 
            \begin{equation}
           B_{\alpha\beta}\tilde{\Gamma}^{\alpha\beta,\mathcal{N}}(t) - Q_{\alpha\beta}G_{\alpha\beta}(1-G_{\alpha\beta})    +\bar{Q}_{\alpha\beta}\bar{G}_{\alpha\beta}(1-\bar{G}_{\alpha\beta})  
                  \label{eq:pp13:proof:strict}
            \end{equation}
        is strictly positive for any $t\in[0,T)$. Hence, we conclude from \eqref{eq:Gamma:alpha:beta} and Lemma \ref{lem:ode:compare} that $\Gamma^{\alpha\beta,\mathcal{P}}(t) >\Gamma^{\alpha\beta,\mathcal{N}}(t)$ for $t\in[0,T)$. On the other hand, from \eqref{eq:Gamma:hat:alpha:beta:N} and Lemma \ref{lem:ode:compare}, we also have $0<\hat{\Gamma}^{\alpha\beta,\mathcal{N}}(t) <  \Gamma^{\alpha\beta,\mathcal{N}}(t)$.

       
        \item If $B_{\alpha\beta}\leq 0$, $G_{\alpha\beta},\bar{G}_{\alpha\beta}\in [0,1)$ and at least one of these coefficients are non-zero, then \eqref{eq:pp13:proof:strict} is strictly negative. Hence, we conclude from \eqref{eq:Gamma:alpha:beta} and Lemma \ref{lem:ode:compare} that $0<\Gamma^{\alpha\beta,\mathcal{P}}(t) <\Gamma^{\alpha\beta,\mathcal{N}}(t)$ for $t\in[0,T)$.
    \end{enumerate}}
\end{proof}

We now investigate the effect of the leaders on the follower. Since the system is symmetric and the leaders do not depend on the mean field term, we can reduce  \eqref{dynGamma} to the following system of Riccati equations: 
\begin{equation}
\label{eq:riccati:reduced}
\left\{\begin{aligned}
-\frac{d {\bf {\bf \Gamma}}^{\mathcal{I},{\bf b}_0}(t)}{dt}=&\ {\bf {\bf \Gamma}}^{\mathcal{I},{\bf b}_0}(t){\bf A}^{{\bf b}_0}+{\bf C}^{\mathcal{I},{\bf b}_0}{\bf \Gamma}^{\mathcal{I},{\bf b}_0}(t)  - {\bf \Gamma}^{\mathcal{I},{\bf b}_0}(t){\bf B}^{{\color{black}\mathcal{I}},{\bf b}_0}{\bf \Gamma}^{\mathcal{I},{\bf b}_0}(t)+{\bf D}^{\mathcal{I},{\bf b}_0},\\ 
{\bf \Gamma}^{\mathcal{I},{\bf b}_0}(T)=&\ {\bf E}^{\mathcal{I},{\bf b}_0},
\end{aligned}\right.
\end{equation}%
where {\color{black}${\bf \Gamma}^{\mathcal{N},{\bf b}_0}(t) \in \mathbb{R}^{5\times 3}$, ${\bf \Gamma}^{\mathcal{P},{\bf b}_0}(t) \in \mathbb{R}^{3\times 3}$}, and  for $\mathcal{I}=\mathcal{N},\mathcal{P}$,  
 \begin{align}
 \label{eq:riccati:reduced:matrix}
      & {\bf A}^{{\bf b}_0}:= \begin{pmatrix}
                {\bf A}^{\alpha\beta}&{\bf 0}_{2,1}\\
                C_1{\bf e}_{1,2}&A_1+B_1
        \end{pmatrix}, \ {\bf B}^{{\color{black}\mathcal{I}},{\bf b}_0}:= \begin{pmatrix}
             {\bf B}^{\alpha\beta,{\color{black}\mathcal{I}}}&{\bf 0}_{2,1} \nonumber\\
               {\color{black}({\bf 0}^\mathcal{I})^\top}&  E_{1}^2 R_{1}^{-1}
        \end{pmatrix},\\
        &{\bf C}^{\mathcal{I},{\bf b}_0} := \begin{pmatrix}
                {\bf C}^{\alpha\beta,\mathcal{I}}&{\color{black}{\bf 0}^\mathcal{I}}\\
                {\color{black}({\bf 0}^\mathcal{I})^\top}& A_1
        \end{pmatrix}\\     
        &{\bf D}^{\mathcal{I},{\bf b}_0}: = \begin{pmatrix}
               {\bf D}^{\alpha\beta,\mathcal{I}}&{\color{black} {\bf 0}^\mathcal{I} }\\
                -Q_1G_1{\bf e}_{1,2}&Q_1(1-F_1)
        \end{pmatrix} \nonumber \\  
        &{\bf E}^{\mathcal{I},{\bf b}_0}: = \begin{pmatrix}
               {\bf E}^{\alpha\beta,\mathcal{I}}&{\color{black} {\bf 0}^\mathcal{I}} \\
                -\bar{Q}_1\bar{G}_1{\bf e}_{1,2}&\bar{Q}_1(1-\bar{F}_1)
        \end{pmatrix}, 
    \end{align}%
{\color{black} with ${\bf 0}^\mathcal{N}={\bf 0}_{4,1}$ and ${\bf 0}^\mathcal{P} = {\bf 0}_{2,1}$}. For $\mathcal{I}=\mathcal{N}$ and $\mathcal{P}$, write the solution of  \eqref{eq:riccati:reduced}   as    \begin{equation}
    \label{Gamma0form}
        {\bf \Gamma}^{\mathcal{I},{\bf b}_0}(t)= \begin{pmatrix}
        {\bf \Gamma}^{\mathcal{I},{\bf b}_0}_{11}(t) & {\bf \Gamma}^{\mathcal{I},{\bf b}_0}_{12}(t) \\
        {\bf \Gamma}_{21}^{\mathcal{I},{\bf b}_0}(t) & \Gamma^{\mathcal{I},{\bf b}_0}_{22}(t)
        \end{pmatrix}
    \end{equation}%
where {\color{black}${\bf \Gamma}_{11}^{\mathcal{N},{\bf b}_0}(t) \in \mathbb{R}^{4\times2}$, ${\bf \Gamma}_{11}^{\mathcal{P},{\bf b}_0}(t) \in \mathbb{R}^{2\times2}$,  ${\bf \Gamma}_{12}^{\mathcal{N},{\bf b}_0}(t) \in \mathbb{R}^{4\times1}$, ${\bf \Gamma}_{12}^{\mathcal{P},{\bf b}_0}(t) \in \mathbb{R}^{2\times1}$, and for $\mathcal{I}=\mathcal{N},\mathcal{P}$,} ${\bf \Gamma}_{21}^{\mathcal{I},{\bf b}_0}(t) \in \mathbb{R}^{1\times2}$, and $\Gamma_{22}^{\mathcal{I},{\bf b}_0}(t) \in \mathbb{R}$. We have the following component-wise comparisons.
\begin{proposition}
    \label{pp:riccati:compare}
    {\color{black}Suppose that $\max\{\vert C_1 \vert,  \vert G_1 \vert,  \vert \bar{Q}_1\bar{G}_1 \vert\}, \max\{\vert B_{\alpha\beta} \vert, \vert G_{\alpha\beta} \vert,  \vert \bar{Q}_{\alpha\beta}\bar{G}_{\alpha\beta}\vert \}>0$.} For $t\in[0,T]$, the components of ${\bf \Gamma}^{\mathcal{N},{\bf b}_0}(t), {\bf \Gamma}^{\mathcal{P},{\bf b}_0}(t), \ t\in[0,T]$ have following properties.
    \begin{enumerate}
    \item ${\bf \Gamma}^{\mathcal{N},{\bf b}_0}_{12}(t)={\color{black}{\bf 0}_{4, 1}}$ and ${\bf \Gamma}^{\mathcal{P},{\bf b}_0}_{12}(t) = {\bf 0}_{2,1}$. 

    \item ${\bf \Gamma}^{\mathcal{N},{\bf b}_0}_{11}(t)={\bf \Gamma}^{\alpha\beta,\mathcal{N}}(t)$, ${\bf \Gamma}^{\mathcal{P},{\bf b}_0}_{11}(t)={\bf \Gamma}^{\alpha\beta,\mathcal{P}}(t)$.

    \item ${\bf \Gamma}^{\mathcal{N},{\bf b}_0}_{21}(t)=\Gamma^{\mathcal{N},{\bf b}_0}_{21}(t){\bf e}_{1, 2}$, ${\bf \Gamma}^{\mathcal{P},{\bf b}_0}_{21}(t)=\Gamma^{\mathcal{P},{\bf b}_0}_{21}(t){\bf e}_{1, 2}$ for some $\mathbb{R}$-valued function $\Gamma^{\mathcal{N},{\bf b}_0}_{21}(t)$, $\Gamma^{\mathcal{P},{\bf b}_0}_{21}(t)$.

    \item $\Gamma^{\mathcal{N},{\bf b}_0}_{22}(t)=\Gamma^{\mathcal{P},{\bf b}_0}_{22}(t) =: \Gamma^{{\bf b}_0}_{22}(t) $, where $\Gamma^{{\bf b}_0}_{22}(t)> 0$ for $t\in[0,T)$.



    \item Suppose that $C_1\geq 0$, $G_1, {\color{black}\bar{G}_1\leq 0}$. Then for any $t\in[0,T)$, we have $\Gamma^{\mathcal{N},{\bf b}_0}_{21}(t)$, $\Gamma^{\mathcal{P},{\bf b}_0}_{21}(t) > 0$, and in addition,
        \begin{enumerate}
            \item $\Gamma^{\mathcal{N},{\bf b}_0}_{21}(t) > \Gamma^{\mathcal{P},{\bf b}_0}_{21}(t)$ if $B_{\alpha\beta}\geq  0$ {\color{black} and $-1<G_{\alpha\beta} , \bar{G}_{\alpha\beta}\leq 0$};
            \item $\Gamma^{\mathcal{N},{\bf b}_0}_{21}(t) <  \Gamma^{\mathcal{P},{\bf b}_0}_{21}(t)$ if $B_{\alpha\beta} \leq  0$ and $0 \leq G_{\alpha\beta},{\color{black}\bar{G}_{\alpha\beta}}< 1$.
        \end{enumerate}
    \item Suppose that $C_1 \leq 0$, $G_1, {\color{black}\bar{G}_1 \ge} 0$. Then for any $t\in[0,T)$, we have $\Gamma^{\mathcal{N},{\bf b}_0}_{21}(t)$, $\Gamma^{\mathcal{P},{\bf b}_0}_{21}(t) < 0$, and in addition,
        \begin{enumerate}
            \item $\Gamma^{\mathcal{N},{\bf b}_0}_{21}(t) <  \Gamma^{\mathcal{P},{\bf b}_0}_{21}(t)$ if $B_{\alpha\beta} \geq  0$ {\color{black} and $-1<G_{\alpha\beta}, \bar{G}_{\alpha\beta}\leq 0$};
            \item $\Gamma^{\mathcal{N},{\bf b}_0}_{21}(t) >  \Gamma^{\mathcal{P},{\bf b}_0}_{21}(t)$ if $B_{\alpha\beta} \leq  0$ and $0 \leq G_{\alpha\beta},{\color{black}\bar{G}_{\alpha\beta}}< 1$.
        \end{enumerate}

    \end{enumerate}
\end{proposition}
\begin{proof}
    We shall omit the proof of  Statements 1) - 4), which can be shown by writing down the equations satisfied by the components and applying Lemma \ref{lem:ode:compare}, along with a uniqueness argument as in the proof of Proposition \ref{symmetry:nash}.  
To prove Statements 5) - 6), by expanding the system and using Statement 3), we see that for $\mathcal{I}=\mathcal{N},\mathcal{P}$, ${ \Gamma}^{\mathcal{I},{\bf b}_0}_{21}$ satisfies \vspace{-0.6cm} {\small 
  \begin{equation}
    \label{eq:G12}
       \left\{ \begin{aligned}
              -\frac{d\Gamma_{21}^{\mathcal{I},{\bf b}_0}(t)}{dt} = &\ \bigg(\tilde{A}- \frac{E_1^2\Gamma^{{\bf b}_0}_{22}(t)}{R_1}  - \frac{D_{\alpha\beta}^2\Gamma^{\alpha\beta,\mathcal{I}}(t)}{R_{\alpha\beta}} \bigg) \Gamma_{21}^{\mathcal{N},{\bf b}_0}(t)  +C_1\Gamma^{{\bf b}_0}_{22}(t)  
              -Q_1G_1 , \\ 
              \Gamma^{\mathcal{I},{\bf b}_0}_{21}(T)=&\ -\bar{Q}_1\bar{G}_1,
        \end{aligned}\right.
    \end{equation}}%
    where $\tilde{A} := A_1+A_{\alpha\beta}+B_{\alpha\beta}$. By Lemma \ref{lem:ode:compare} and the fact that $\Gamma^{{\bf b}_0}_{22}(t)> 0$, for any $t\in[0,T)$, we have $\Gamma^{\mathcal{N},{\bf b}_0}_{21}(t)$, $\Gamma^{\mathcal{P},{\bf b}_0}_{21}(t) > 0$ when $C_1\geq 0$, $G_1 \leq  0$ and $\bar{G}_1\leq 0$; $\Gamma^{\mathcal{N},{\bf b}_0}_{21}(t)$, $\Gamma^{\mathcal{P},{\bf b}_0}_{21}(t) < 0$ when $C_1 \leq  0$, $G_1 \geq  0$ and $\bar{G}_1 \geq 0$. The ordering of $\Gamma^{\mathcal{N},{\bf b}_0}_{21}$ and  $\Gamma^{\mathcal{P},{\bf b}_0}_{21}$ now follows from Lemma \ref{lem:ode:compare} and Proposition \ref{pp:leader:compare}.
\end{proof}



\subsubsection{Analysis of Lyapunov Equations}
In Theorem \ref{thm:cost}, we have expressed the cost functional of the follower in terms of the solution of a Lyapunov equation. Under Model \ref{degeneratedmodel}, we can express the the optimal control for the follower using the reduced Riccati equation \eqref{eq:riccati:reduced}, whereas the associated cost functional can be expressed in terms of the solution of a reduced Lyapunov equation.  The following result is a simple consequence of Theorem \ref{thm:cost}, whose proof is omitted. 
\begin{lemma}
\label{lem:cost}
  Consider Model \ref{degeneratedmodel}. For $\mathcal{I}=\mathcal{N}$ and $\mathcal{P}$, the optimal control $u^{\mathcal{I},{\bf b}_0}$ for the follower can be written as  
         \begin{equation}
       \label{eq:control:follower:base}
       \begin{aligned}
        u^{\mathcal{I},{\bf b}_0}_1(t) =&\ -\frac{E_1}{R_1} \Big( \Gamma_{0}(t) x^{\mathcal{I},{\bf b}_0}(t) + \Gamma_{21}^{\mathcal{I},{\bf b}_0}(t)(x_\alpha^{\mathcal{I},{\bf b}_0}(t) + x^{\mathcal{I},{\bf b}_0}_\beta(t))  + (\Gamma^{{\bf b}_0}_{22}(t)-\Gamma_0(t))z^{\mathcal{I},{\bf b}_0}(t)  \Big),
        \end{aligned}
    \end{equation}%
    where $x_\alpha^{\mathcal{I},{\bf b}_0}, x^{\mathcal{I},{\bf b}_0}_\beta$ and $z^{\mathcal{I},{\bf b}_0}$ are respectively the state dynamics of the $\alpha$-leader, $\beta$-leader and the mean field term when the optimal controls are employed. In addition, the cost functional of the follower $J_1^{\mathcal{I},{\bf b}_0}$ is given by     
        \begin{equation}
        \label{eq:lem:cost:2}
        \begin{aligned}
            J_1^{\mathcal{I},{\bf b}_0} =&\  \mathbb{E}[ \langle \tilde{{\bf x}} ,  {\bf K}^{\mathcal{I},{\bf b}_0}\tilde{{\bf x}} \rangle(0)]   + \int_0^T \textup{tr}(\tilde{ {\bf \Sigma} }^\top {\bf K}^{\mathcal{I},{\bf b}_0}(t) \tilde{ {\bf \Sigma} })  dt \\ &\  + \int_0^T \left(\frac{E^2_1\Gamma^2_{0}(t)}{R_1}+Q_1\right)(\Sigma_0+P_0)(t) dt  + \bar{Q}_1(\Sigma_0(T) + P_0(T)),
            \end{aligned}
        \end{equation}%
    where   $\tilde{\bf x}(0) := (\eta_\alpha,\eta_\beta,\mathbb{E}[\eta_1])$, $\tilde{ {\bf \Sigma} } := \textup{diag}(\sigma_\alpha,\sigma_\beta,0)$, $\Gamma_0$ is the solution of \eqref{Gamma0}, 
     \begin{equation}\label{Sigma0Lambda0}
     \begin{aligned}
    \Sigma_0(t)&=\int_0^te^{\int_s^t 2[A_1-\frac{E_1}{R_1}\Gamma_{0}(u)]du}\sigma^2_1ds,  \ 
    P_0(t)&=e^{\int_0^t 2[A_1-\frac{E_1}{R_1}\Gamma_{0}(s)] ds} \mathbb{E}[|\eta_1-\mathbb{E}[\eta_1]|^2],
    \end{aligned}
    \end{equation}%
     and  ${\bf K}^{\mathcal{I},{\bf b}_0}(t)$ is the solution of the Lyapunov equation  
     \begin{equation}
     \label{Lyapunov} \left\{
         \begin{aligned}
             -\frac{d{\bf K}^{\mathcal{I},{\bf b}_0}(t)}{dt} =&\  {\bf K}^{\mathcal{I},{\bf b}_0}(t) ({\bf A}^{{\bf b}_0} - {\bf B}^{{\color{black}\mathcal{I},}{\bf b}_0} {\bf \Gamma}^{\mathcal{I},{\bf b}_0}(t))  + ({\bf A}^{{\bf b}_0} - {\bf B}^{{\color{black}\mathcal{I},}{\bf b}_0}{\bf \Gamma}^{\mathcal{I},{\bf b}_0}(t))^{\top}{\bf K}^{\mathcal{I},{\bf b}_0}(t) \\
             &\ +  {\bf Q}^{{\bf b}_0} + {\bf R}^{\mathcal{I},{\bf b}_0}(t) , \\
             {\bf K}^{\mathcal{I},{\bf b}_0}(T) =&\ {\bf \bar{Q}}^{{\bf b}_0}.
         \end{aligned}\right.
     \end{equation}%
   Here,  ${\bf \Gamma}^{\mathcal{I},{\bf b}_0}(t)$  is the solution of   \eqref{eq:riccati:reduced}, ${\bf A}^{{\bf b}_0}, {\bf B}^{{\color{black}\mathcal{I},}{\bf b}_0}$ are given by  \eqref{eq:riccati:reduced:matrix}
,  and 
        \begin{align*}
        &{\bf R}^{\mathcal{I},{\bf b}_0}(t) := \frac{E_1^2}{R_1} \begin{pmatrix} ({\bf \Gamma}_{21}^{\mathcal{I},{\bf b}_0} )^{\top}(t) \\ \Gamma_{22}^{{\bf b}_0}(t) \end{pmatrix} \begin{pmatrix} {\bf \Gamma}^{\mathcal{I},{\bf b}_0}_{21}(t)& \Gamma^{{\bf b}_0}_{22}(t) \end{pmatrix},\\
        &     {\bf Q}^{{\bf b}_0} := Q_1 \begin{pmatrix}
                G^2_1{\bf e}_{2, 2} &-(1-F_1)G_1{\bf e}_{2,1}\\
                -(1-F_1)G_1{\bf e}_{1,2}& (1-F_1)^2
        \end{pmatrix},\\
        &     {\bf \bar{Q}}^{{\bf b}_0} := \bar{Q}_1 \begin{pmatrix}
                \bar{G}^2_1{\bf e}_{2, 2} &-(1-\bar{F}_1)\bar{G}_1{\bf e}_{2,1}\\
                -(1-\bar{F}_1)\bar{G}_1{\bf e}_{1,2}& (1-\bar{F}_1)^2
        \end{pmatrix}.
        \end{align*}%
\end{lemma}
We proceed to compare the solutions ${\bf K}^{\mathcal{I},{\bf b}_0}$ for $\mathcal{I}=\mathcal{N},\mathcal{P}$ in a component-wise manner, which eventually leads to  Theorem \ref{thm:compare}. Notice that it is clear from \eqref{Lyapunov} that ${\bf K}^{\mathcal{I},{\bf b}_0}$ is symmetric. In what follows, we denote the components of ${\bf K}^{\mathcal{I},{\bf b}_0}$ by  
    \begin{equation}
    \label{eq:K:b}
        {\bf K}^{\mathcal{I},{\bf b}_0}(t) = \begin{pmatrix}
                {\bf K}^{\mathcal{I},{\bf b}_0}_{11}(t) & {\bf K}^{\mathcal{I},{\bf b}_0}_{12}(t) \\
                ({\bf K}^{\mathcal{I},{\bf b}_0}_{12}(t))^{\top} & K^{\mathcal{I},{\bf b}_0}_{22}(t)
        \end{pmatrix},
    \end{equation}%
where ${\bf K}^{\mathcal{I},{\bf b}_0}_{11}(t) \in \mathbb{R}^{2\times 2}$, ${\bf K}^{\mathcal{I},{\bf b}_0}_{12}(t) \in \mathbb{R}^{2\times 1}$, $K^{\mathcal{I},{\bf b}_0}_{22}(t) \in \mathbb{R}$.

\begin{proposition}
\label{poK}
 {\color{black}Under Assumption \ref{ass:compare}, f}or any $t\in[0,T]$ and $\mathcal{I}=\mathcal{N},\mathcal{P}$, the solution ${\bf K}^{\mathcal{I},{\bf b}_0}(t)$  of \eqref{Lyapunov}  satisfies the following:
\begin{enumerate}
    \item $K^{\mathcal{N},{\bf b}_0}_{22}(t)=K^{\mathcal{P},{\bf b}_0}_{22}(t)=:K^{{\bf b}_0}_{22}(t)$, where $K^{{\bf b}_0}_{22}(t)>0$ for $t\in[0,T)$. 

    \item ${\bf K}_{12}^{\mathcal{I},{\bf b}_0}(t) = K^{\mathcal{I},{\bf b}_0}_{12}(t){\bf e}_{2, 1}$ for some $\mathbb{R}$-valued function $K^{\mathcal{I},{\bf b}_0}(t)$.

    \item ${\bf K}^{\mathcal{I},{\bf b}_0}_{11}(t)=K^{\mathcal{I},{\bf b}_0}_{11}(t){\bf e}_{2,2}$ for some $\mathbb{R}$-valued function $K^{\mathcal{I},{\bf b}_0}_{11}(t)$;
    
    \item Suppose that Assumption \ref{ass:compare} holds. Then for $t\in [0,T)$, we have
        \begin{enumerate}
            \item $K^{\mathcal{I},{\bf b}_0}_{12}(t) >0$   if $C_1 \geq 0$ ,  $G_1 \leq 0$  and  $\bar{G}_1\leq 0$; 
            \item $K^{\mathcal{I},{\bf b}_0}_{12}(t) < 0$ if $C_1 \leq 0$ ,  $G_1 \geq  0$ and  $\bar{G}_1 \geq 0$;
           \item {\color{black} $K^{{\bf b}_0}_{11}(t)>0$.}
        \end{enumerate}
  {\color{black}   \item Suppose that Conditions 1-4) of Assumption \ref{ass:compare} holds. Then, for $t\in[0,T)$. we have
        \begin{enumerate}
            \item $K^{\mathcal{N},{\bf b}_0}_{11}(t) > K^{\mathcal{P},{\bf b}_0}_{11}(t)$,  $K^{\mathcal{N},{\bf b}_0}_{12}(t) > K^{\mathcal{P},{\bf b}_0}_{12}(t)$ if  $B_\alpha,A_\beta \geq 0$, $C_1\geq 0$, $-1 < G_\alpha,G_\beta,\bar{G}_\alpha, \bar{G}_\beta\leq 0$ and $G_1,\bar{G_1}\leq 0$;

            \item $K^{\mathcal{N},{\bf b}_0}_{11}(t) > K^{\mathcal{P},{\bf b}_0}_{11}(t)$,  $K^{\mathcal{N},{\bf b}_0}_{12}(t) < K^{\mathcal{P},{\bf b}_0}_{12}(t)$ if  $B_\alpha,A_\beta \ge 0$, $C_1 \le 0$, $-1 < G_\alpha,G_\beta,\bar{G}_\alpha, \bar{G}_\beta\leq 0$ and $G_1,\bar{G_1}\geq 0$;

            \item    $K^{\mathcal{N},{\bf b}_0}_{11}(t)<  K^{\mathcal{P},{\bf b}_0}_{11}(t)$,  $K^{\mathcal{N},{\bf b}_0}_{12}(t) < K^{\mathcal{P},{\bf b}_0}_{12}(t)$ if  $B_\alpha,A_\beta \le 0$, $C_1 \ge 0$, $0\le G_\alpha,G_\beta,\bar{G}_\alpha, \bar{G}_\beta< 1$ and $G_1,\bar{G_1}\leq 0$;

             \item $K^{\mathcal{N},{\bf b}_0}_{11}(t) < K^{\mathcal{P},{\bf b}_0}_{11}(t)$,  $K^{\mathcal{N},{\bf b}_0}_{12}(t) > K^{\mathcal{P},{\bf b}_0}_{12}(t)$ if  $B_\alpha,A_\beta \le 0$, $C_1 \le 0$, $0\le G_\alpha,G_\beta,\bar{G}_\alpha, \bar{G}_\beta<1$ and $G_1,\bar{G_1}\geq 0$.
           
        \end{enumerate}}

\end{enumerate}
\end{proposition}

\begin{proof}
    We shall omit the proof of Statements 1)-3), 
which can be shown by using the same argument as in the proof of Proposition \ref{pp:riccati:compare} and utilizing Lemma \ref{lem:ode:compare}.
 \begin{enumerate}
 \setcounter{enumi}{3}
 
 
     \item {\color{black} We shall consider prove Statement 4a), the proof of Statement 4b) can be shown by a parallel argument. By expanding the system \eqref{Lyapunov} and using Statement 2), we see that for $\mathcal{I}=\mathcal{N}, \mathcal{P}$, $K^{\mathcal{I},{\bf b}_0}_{12}$ satisfies
  
     \begin{equation}
     \label{eq:K12}
    \left\{ \begin{aligned}
          -\frac{dK^{\mathcal{I},{\bf b}_0}_{12}(t)}{dt}&= K^{\mathcal{I},{\bf b}_0}_{12}(t)\bigg(\tilde{A}+B_1  -\frac{E_1^2\Gamma^{{\bf b}_0}_{22}(t)}{R_1}   -\frac{D_{\alpha\beta}^2\Gamma^{\alpha\beta,\mathcal{I}}(t)}{R_{\alpha\beta}}\bigg) + \varrho^{\mathcal{I},{\bf b}_0}(t), \\
          K^{\mathcal{I},{\bf b}_0}_{12}(t) &= -\bar{Q}_1\bar{G}_1(1-\bar{F}_1), 
         \end{aligned}\right.
     \end{equation}%
where  
    \begin{align*}
        \varrho^{\mathcal{I},{\bf b}_0}(t) &=   C_1 K^{{\bf b}_0}_{22}(t) +\frac{E_1^2\left(\phi^{{\bf b}_0}(t) + K^{{\bf b}_0}_{22}(t) \right)\Gamma^{\mathcal{I},{\bf b}_0}_{21}(t)}{R_1} \  -(1-F_1)G_1Q_1,
    \end{align*}
and $\phi^{{\bf b}_0}(t):=\Gamma^{{\bf b}_0}_{22}(t)-2K^{{\bf b}_0}_{22}(t)$.  We shall prove that
    \begin{equation}
        \label{eq:pp16:compare:rho}
        \varrho^{\mathcal{I},{\bf b}_0}(t)\geq 0, \ \max\{\varrho^{\mathcal{I},{\bf b}_0}(t), - \bar{Q}_1\bar{G}_1(1-\bar{F}_1)\} >0,
    \end{equation}
 for any $t\in[0,T)$. This allows us to deduce $K_{12}^{\mathcal{I},{\bf b}_0}(t)>0$. 

Suppose that Conditions 1-4) of Assumption \ref{ass:compare} hold. Consider the equation satisfied by $\phi^{{\bf b}_0}$: 
 
     \begin{equation}
    \left\{ \begin{aligned}
          -\frac{d\phi^{{\bf b}_0}(t)}{dt}&=2
          \phi^{{\bf b}_0}(t)\bigg(A_1+B_1  -\frac{E_1^2\Gamma^{{\bf b}_0}_{22}(t)}{R_1}\bigg)  -
      \left(B_1+\frac{E_1^2\Gamma^{{\bf b}_0}_{22}(t)}{R_1}\right) \Gamma^{{\bf b}_0}_{22}(t)\\
      &\ +Q_1(2F_1-1)(1-F_1),\\ \phi^{\mathcal{I},{\bf b}_0}(T)&= \bar{Q}_1(2\bar{F}_1-1)(1-\bar{F}_1).
         \end{aligned}\right.
     \end{equation}%
By a direct expansion, one can verify that $0\leq \Gamma_{22}^{{\bf b}_0}(t) \leq U$, where $U$ is given in Assumption \ref{ass:compare}. Hence, under Condition 4) of Assumption \ref{ass:compare}, an application of Lemma \ref{lem:ode:compare} implies that $\phi^{{\bf b}_0}(t)> 0$ for any $t\in[0,T)$. By Condition 3) of Assumption \ref{ass:compare} and Statement 1) of the proposition, we obtain \eqref{eq:pp16:compare:rho} as desired. 
  
Now, suppose that Conditions 1-3), and 5) of Assumption \ref{ass:compare} hold. In this case, we write 
    \begin{align}
    \label{eqn:varrho}
        \varrho^{\mathcal{I},{\bf b}_0}(t) &=  \psi^{\mathcal{I},{\bf b}_0}(t)K^{{\bf b}_0}_{22}(t)  +\frac{E_1^2\Gamma^{\mathcal{I},{\bf b}_0}_{21}(t)\Gamma^{{\bf b}_0}_{22}(t)}{R_1}  -(1-F_1)G_1Q_1,
    \end{align}%
}%
     where $\psi^{\mathcal{I},{\bf b}_0}(t) := C_1-E_1^2\Gamma^{\mathcal{I},{\bf b}_0}_{21}(t)/R_1$. We claim that for any $t\in[0,T)$, $\psi^{\mathcal{I},{\bf b}_0}(t){\color{black}>}0$ when $C_1 {\color{black} >} 0$, $G_1{\color{black}\le} 0$, $\bar{G}_1\leq 0$, {\color{black} which again implies \eqref{eq:pp16:compare:rho}.  }
     
     {\color{black}To establish the claim, b}y \eqref{eq:G12}, we see that $\psi^{\mathcal{I},{\bf b}_0}$ satisfies the equation
     \begin{equation}
        \label{eq:theta}
        \left\{\begin{aligned}
            -\frac{d\psi^{\mathcal{I},{\bf b}_0}(t)}{dt} = &\  \psi^{\mathcal{I},{\bf b}_0}(t) \bigg(\tilde{A}  -\frac{E^2_1 \Gamma^{{\bf b}_0}_{22}(t)}{R_1}\bigg) +{\color{black}V}  -   C_1\tilde{A}  +\frac{D_{\alpha\beta}^2E_1^2\Gamma^{\alpha\beta,\mathcal{I}}(t)\Gamma^{\mathcal{I},{\bf b}_0}_{21}(t)}{R_1R_{\alpha\beta}}, \\
            \psi^{\mathcal{I},{\bf b}_0}(T) =&\ C_1+{\color{black}\bar{V}}. 
        \end{aligned}\right.
        \end{equation}%
     By Propositions \ref{pp:leader:compare}-\ref{pp:riccati:compare},  $\Gamma^{\alpha\beta,\mathcal{I}}(t)\Gamma^{\mathcal{I},{\bf b}_0}_{21}(t) \geq 0$ when $C_1 {\color{black} \geq }0$, $G_1{\color{black}\le}0$ and $\bar{G}_1\leq0$. Hence, solving \eqref{eq:theta} gives
     \begin{align}
             \label{eq:pp:lyapunov:s4:pf}
      &\ e^{-\int_t^T \left(\tilde{A} - \frac{E_1^2\Gamma^{{\bf b}_0}_{22}(s)}{R_1}\right) ds} \psi^{\mathcal{I},{\bf b}_0}(t) \geq  C_1+{\color{black}\bar{V}}
        - \bigg( C_1\tilde{A}  - {\color{black}V}
      \bigg)   \int_t^T e^{-\int_s^T \left(\tilde{A} - \frac{E_1^2\Gamma^{{\bf b}_0}_{22}(u)}{R_1}\right) du } ds.
       \end{align}%
     Suppose that {\color{black}Condition 5a) of Assumption \ref{ass:compare}}  holds. Then it is clear from \eqref{eq:pp:lyapunov:s4:pf} that $\psi^{\mathcal{I},{\bf b}_0}(t)> 0$ for $t\in[0,T)$. Now suppose alternatively that,    {\color{black}Conditions 1,2, and 5b)} of Assumption \ref{ass:compare}  hold. 

{\color{black} Finally, for $t\in[0,T)$, the fact that $K_{11}^{\mathcal{I},{\bf b}_0}(t)$ can be shown by the positive definiteness of ${\bf R}^{\mathcal{I},{\bf b}_0}(t)$, and the positive semi-definiteness of ${\bf Q}^{{\bf b}_0}$ and $\bar{{\bf Q}}^{{\bf b}_0}$. }

\item 
{\color{black} We only show (a), the rest of the statements can be shown by a parallel argument. We first prove that $K_{12}^{\mathcal{N},{\bf b}_0}(t)>K_{12}^{\mathcal{P},{\bf b}_0}(t)$ for $t\in[0,T)$. Notice that, by Propositions \ref{pp:leader:compare} and \ref{pp:riccati:compare}, we have $\Gamma^{\alpha\beta,\mathcal{N}}(t)<\Gamma^{\alpha\beta,\mathcal{P}}(t)$ and $\Gamma_{21}^{\mathcal{N},{\bf b}_0}(t)>\Gamma_{21}^{\mathcal{P},{\bf b}_0}(t)$ for $t\in[0,T)$. Hence, the claim follows from \eqref{eq:K12} and Lemma \ref{lem:ode:compare}.

We proceed to show $K^{\mathcal{N},{\bf b}_0}_{11}(t) > K^{\mathcal{P},{\bf b}_0}_{11}(t)$. By expanding the system \eqref{Lyapunov}  using  Statement 3) of the proposition, we know that $K_{11}^{\mathcal{I},{\bf b}_0}$ satisfies  
 \begin{equation}
 \label{eq:K11}
   \left\{ \begin{aligned}
         -\frac{dK_{11}^{\mathcal{I},{\bf b}_0}(t)}{dt}=&\ 2K^{\mathcal{I},{\bf b}_0}_{11}(t)\bigg(A_{\alpha\beta}+B_{\alpha\beta}  -\frac{D_{\alpha\beta}^2\Gamma^{\alpha\beta,\mathcal{I}}(t)}{R_{\alpha\beta}}\bigg)   +\varsigma^{\mathcal{I},{\bf b}_0}(t) +G_1^2 Q_1,\\ 
          K^{\mathcal{I},{\bf b}_0}_{11}(T)=&\  \bar{G}_1^2 \bar{Q}_1,
         \end{aligned}\right.
     \end{equation}%
where     
\begin{align}
   \varsigma^{\mathcal{I},{\bf b}_0}(t)&=2C_1K_{12}^{\mathcal{I},{\bf b}_0}(t) + \frac{E_1^2\Gamma_{21}^{\mathcal{I},{\bf b}_0}(t)}{R_1}\varphi^{\mathcal{I},{\bf b}_0}(t) \nonumber \\
   &=\frac{E_1^2(\Gamma_{21}^{\mathcal{I},{\bf b}_0}(t))^2}{R_1}+2\psi^{\mathcal{I},{\bf b}_0}(t)K_{12}^{\mathcal{I},{\bf b}_0}(t),
   \label{eqn:varsigma}
\end{align}
and $\varphi^{\mathcal{I},{\bf b}_0}(t):=\Gamma^{\mathcal{I},{\bf b}_0}_{21}(t)-2K^{\mathcal{I},{\bf b}_0}_{12}(t)$ satisfies 

       \begin{equation}
       \label{eq:varphi}
    \left\{ \begin{aligned}
          -\frac{d\varphi^{\mathcal{I},{\bf b}_0}(t)}{dt}=&\ \varphi^{\mathcal{I},{\bf b}_0}(t)\bigg(\tilde{A}+B_1  -\frac{E_1^2\Gamma^{{\bf b}_0}_{22}(t)}{R_1}    -\frac{D_{\alpha\beta}^2\Gamma^{\alpha\beta,\mathcal{I}}(t)}{R_{\alpha\beta}}\bigg) +
      C_1 \phi^{{\bf b}_0}(t)\\ &\  -\left(  
  B_1+\frac{2E_1^2(\Gamma^{{\bf b}_0}_{22}(t)-K^{{\bf b}_0}_{22}(t))}{R_1}\right)\Gamma^{\mathcal{I},{\bf b}_0}_{21}(t)   -(2F_1-1)G_1Q_1,\\  \varphi^{\mathcal{I},{\bf b}_0}(T)=&\ -\bar{Q}_1(2\bar{F}_1-1)\bar{G}_1.
         \end{aligned}\right.
     \end{equation}%
     By Lemma \ref{lem:ode:compare} and the proven fact that  $K^{\mathcal{N},{\bf b}_0}_{12}(t) > K^{\mathcal{P},{\bf b}_0}_{12}(t)$, it suffices to show that $\varphi^{\mathcal{N},{\bf b}_0}(t) > \varphi^{\mathcal{P},{\bf b}_0}(t)$ for $t\in[0,T)$. By noticing that $C_1\phi^{{\bf b}_0}(t)\geq 0$, $B_1+ \frac{2E_1^2\Gamma_{22}^{{\bf b}_0}(t)}{R_1} \leq B_1 + \frac{2E_1^2U}{R_1} \leq 0$ and $\Gamma^{\alpha\beta,\mathcal{N}}(t)<\Gamma^{\alpha\beta,\mathcal{P}}(t)$, the claim follows from Lemma \ref{lem:ode:compare}. }
   \end{enumerate}

\end{proof}

\subsubsection{Proof of Theorem \ref{thm:compare}}
\label{app:proof:compare}

     We shall only prove  Statement 1) of Theorem \ref{thm:compare}.  Statements 2) - 4) can be proven by a parallel argument. Using \eqref{eq:lem:cost}, we observe that the difference 
      \begin{equation}
      \label{eq:cost:difference}
         \begin{aligned}
              J^{\mathcal{N},{\bf b}_0}_1- J^{\mathcal{P},{\bf b}_0}  
           &= \left(K^{\mathcal{N},{\bf b}_0}_{11}(0)-K^{\mathcal{P},{\bf b}_0}_{11}(0)\right)\mathbb{E}[(\eta_\alpha+\eta_\beta)^2]\\ &\hspace{-0.2cm} +\left( \sigma_\alpha^2+\sigma_\beta^2  \right)\int_0^T \left(K^{\mathcal{N},{\bf b}_0}_{11}-K^{\mathcal{P},{\bf b}_0}_{11}\right)(t) dt\\
            &\hspace{-0.2cm} + 2\left(K^{\mathcal{N},{\bf b}_0}_{12}(0)-K^{\mathcal{P},{\bf b}_0}_{12}(0)\right) \mathbb{E}[\eta_1]\mathbb{E}[\eta_\alpha+\eta_\beta] 
        \end{aligned}
    \end{equation}%
     is a linear function of $\mathbb{E}[\eta_1]$ {\color{black} given fixed $\eta_\alpha$ and $\eta_\beta$.} {\color{black} We first prove Statement 1) under Conditions 1-4) of Assumption \ref{ass:compare}. Indeed, the result is an immediate consequence of Statement 5) of Proposition \ref{poK}. }
     
    {\color{black} Next, we assume that Conditions 1-3) and 5) are satisfied.} We begin by considering the case where $\mathbb{E}[\eta_\alpha+\eta_\beta]\neq 0$, and show that there exist constants $\underline{L},\bar{L}>0$ such that $J^{\mathcal{N},{\bf b}_0}_1- J^{\mathcal{P},{\bf b}_0}>0$ if $\mathbb{E}[\eta_1]\mathbb{E}[\eta_\alpha+\eta_\beta] >-\underline{L}$; and $J^{\mathcal{N},{\bf b}_0}_1- J^{\mathcal{P},{\bf b}_0}<0$ if $\mathbb{E}[\eta_1]\mathbb{E}[\eta_\alpha+\eta_\beta] <-\bar{L}$. This shows that the coefficient of $\mathbb{E}[\eta_1]$ is non-zero, and the statement follows from linearity.

    To proceed, we consider the sub-optimal control for the follower under the Pareto game: 
    \begin{equation*}
    \begin{aligned}
   \tilde{v}^{\mathcal{P}}_1(t)=&\ -\frac{E_1}{R_1}\big(\Gamma_{0}(t)\tilde{x}_1^{\mathcal{P},{\bf b}_0}(t) + \Gamma^{\mathcal{N},{\bf b}_0}_{21}(t)(x^{\mathcal{P},{\bf b}_0}_{\alpha}(t) +x^{\mathcal{P},{\bf b}_0}_{\beta}(t)) +(\Gamma^{{\bf b}_0}_{22}(t)-\Gamma_{0}(t))\tilde{z}^{\mathcal{P},{\bf b}_0}(t)\big),
   \end{aligned}
    \end{equation*}%
 where  $\tilde{x}_1^{\mathcal{P},{\bf b}_0}$, $\tilde{z}^{\mathcal{P},{\bf b}_0}$ are the corresponding state dynamics of the follower and the mean field term when the control  $\tilde{v}_1^\mathcal{P}$ is implemented. Comparing with the optimal control \eqref{eq:control:follower:base}, the function  $\Gamma^{\mathcal{P},{\bf b}_0}_{12}$ therein is replaced by $\Gamma^{\mathcal{N},{\bf b}_0}_{12}$, which leads to the sub-optimality of $\tilde{v}^\mathcal{P}$. Following the proof of Lemma \ref{lem:cost}, one can easily show that the sub-optimal cost functional of the follower $\tilde{J}_1^{\mathcal{P},{\bf b}_0}$ under the Pareto game is given by  
      \begin{equation}
        \label{eq:tilde:J}
        \begin{aligned}
            \tilde{J}_1^{\mathcal{P},{\bf b}_0} =&\  \mathbb{E}[ \langle \tilde{{\bf x}},  \tilde{{\bf K}}^{\mathcal{P},{\bf b}_0}\tilde{{\bf x}} \rangle(0)]+\bar{Q}_1(\Sigma_0(T) + P_0(T))  \\ &\  + \int_0^T \left(\textup{tr}( \tilde{ {\bf \Sigma} }^\top  \tilde{{\bf K}}^{\mathcal{P},{\bf b}_0}(t) \tilde{ {\bf \Sigma} })  +\left(\frac{E^2_1\Gamma^2_{0}(t)}{R_1}+Q_1\right)(\Sigma_0(t)+P_0(t))\right)dt,
            \end{aligned}
        \end{equation}%
where ${\bf \tilde{K}}^{\mathcal{P},{\bf b}_0}$ is the solution of \eqref{Lyapunov} for $\mathcal{I}=\mathcal{P}$, except that all the occurrences of $\Gamma_{12}^{\mathcal{P},{\bf b}_0}$ in ${\bf \tilde{\Gamma}}^{\mathcal{P},{\bf b}_0}$ and ${\bf R}^\mathcal{P}_0$ of the equation are replaced by $\Gamma_{12}^{\mathcal{N},{\bf b}_0}$. To be exact, by following the proof of Proposition  \ref{poK}, we can write  
    \begin{equation}
        {\bf \tilde{K}}^{\mathcal{P},{\bf b}_0}(t) =  \begin{pmatrix}
                 \tilde{K}^{\mathcal{P},{\bf b}_0}_{11}(t){\bf e}_{2,2}  &  \tilde{K}^{\mathcal{P},{\bf b}_0}_{12}(t){\bf e}_{2,1} \\
             \tilde{K}^{\mathcal{P},{\bf b}_0}_{12}(t){\bf e}_{1,2} & K^{{\bf b}_0}_{22}(t)
        \end{pmatrix},
    \end{equation}%
where $\tilde{K}^{\mathcal{P},{\bf b}_0}_{11}$ and  
$\tilde{K}^{\mathcal{P},{\bf b}_0}_{12}$ satisfies the equations \eqref{eq:K11} and 
\eqref{eq:K12} for $\mathcal{I}=\mathcal{P}$, except the occurrences of the terms 
{\color{black}$\Gamma_{21}^{\mathcal{P},{\bf b}_0}$
and $\psi^{\mathcal{P},{\bf b}_0}$ in \eqref{eqn:varsigma}  and \eqref{eqn:varrho} are replaced by $\Gamma_{21}^{\mathcal{N},{\bf b}_0}$
and $\psi^{\mathcal{N},{\bf b}_0}$, respectively.}
Since $B_{\alpha\beta}> 0$, 
 ${\color{black}-1<} G_{\alpha\beta}, \bar{G}_{\alpha\beta} {\color{black} \leq } 0$, $C_1> 0$, $G_1{\color{black}\le} 0$ and $\bar{G}_1\leq0$, by Proposition \ref{pp:riccati:compare} and the proof of Proposition \ref{poK}, we have, for any $t\in[0,T)$, $\Gamma_{21}^{\mathcal{N},{\bf b}_0}(t), \psi^{\mathcal{N},{\bf b}_0}(t)> 0$. Hence, by Lemma \ref{lem:ode:compare}, $K_{12}^{\mathcal{N},{\bf b}_0}(t), \tilde{K}_{12}^{\mathcal{P},{\bf b}_0}(t),K_{11}^{\mathcal{N},{\bf b}_0}(t),$ $\tilde{K}_{11}^{\mathcal{P},{\bf b}_0}(t)> 0$. {\color{black}Next, by} Proposition \ref{pp:leader:compare}, $\Gamma^{\alpha\beta,\mathcal{N}}(t) <  \Gamma^{\alpha\beta,\mathcal{P}}(t)$ when $B_{\alpha\beta}\ge 0$, ${\color{black}-1<} G_{\alpha\beta}, \bar{G}_{\alpha\beta} {\color{black} \leq } 0${\color{black}. Thus,} {\color{black}by} \eqref{eq:K12}, \eqref{eq:K11}, Lemma \ref{lem:ode:compare} and the equation satisfied by $\tilde{K}^{\mathcal{P},{\bf b}_0}_{11},\tilde{K}^{\mathcal{P},{\bf b}_0}_{12}${\color{black}, we infer} that {\color{black}$K_{11}^{\mathcal{N},{\bf b}_0}(t)> \tilde{K}_{11}^{\mathcal{P},{\bf b}_0}(t)$} and  $K_{12}^{\mathcal{N},{\bf b}_0}(t)> \tilde{K}_{12}^{\mathcal{P},{\bf b}_0}(t)$ for $t\in[0,T)$. 
  By \eqref{eq:lem:cost} and \eqref{eq:tilde:J}, we have 
          \begin{align}
                  \label{eq:JN-tildeJP}
             J^{\mathcal{N},{\bf b}_0}_1- \tilde{J}^{\mathcal{P},{\bf b}_0} 
           &= \left(K^{\mathcal{N},{\bf b}_0}_{11}(0)-\tilde{K}^{\mathcal{P},{\bf b}_0}_{11}(0)\right)\mathbb{E}[(\eta_\alpha+\eta_\beta)^2] \nonumber\\ &\hspace{-0.2cm} +\left( \sigma_\alpha^2+\sigma_\beta^2  \right)\int_0^T \left(K^{\mathcal{N},{\bf b}_0}_{11}-\tilde{K}^{\mathcal{P},{\bf b}_0}_{11}\right)(t) dt \nonumber\\
            &\hspace{-0.2cm} + 2\left(K^{\mathcal{N},{\bf b}_0}_{12}(0)-\tilde{K}^{\mathcal{P},{\bf b}_0}_{12}(0)\right) \mathbb{E}[\eta_1]\mathbb{E}[\eta_\alpha+\eta_\beta] .
        \end{align}%
     Hence, by \eqref{eq:JN-tildeJP} and the above discussions, we infer that there exists $\underline{L}>0$ {\color{black}independent of $\eta_1$} such that  $J_1^{\mathcal{N},{\bf b}_0} > \tilde{J}^{\mathcal{P},{\bf b}_0}_1$ whenever $\mathbb{E}[\eta_1]\mathbb{E}[\eta_\alpha+\eta_\beta] > -\underline{L}$. 

     Next, we consider a sub-optimal strategy $\tilde{v}^\mathcal{N}_1$ for the follower under the Nash game, where   
     \begin{equation*}
 \begin{aligned}
      \tilde{v}^{\mathcal{N}}_1(t):=&\ -\frac{E_1}{R_1}\big(\Gamma_{0}(t)\tilde{x}^{\mathcal{N},{\bf b}_0}(t) + \Gamma^{\mathcal{P},{\bf b}_0}_{21}(t)(x^{\mathcal{N},{\bf b}_0}_{\alpha}(t) +x^{\mathcal{N},{\bf b}_0}_{\beta}(t))+(\Gamma^{{\bf b}_0}_{22}(t)-\Gamma_{0}(t))\tilde{z}^{\mathcal{N},{\bf b}_0}(t) \big).
      \end{aligned}
    \end{equation*}%
Let $\tilde{J}_1^{\mathcal{N},{\bf b}_0}$ be the sub-optimal cost functional of the follower using the control $\tilde{v}_1^\mathcal{N}$. Using the same argument as the above, one can show that there exists $\bar{L}>0$ such that $\tilde{J}^{\mathcal{N},{\bf b}_0}_1- J^{\mathcal{P},{\bf b}_0}<0$ whenever $\mathbb{E}[\eta_1]\mathbb{E}[\eta_\alpha+\eta_\beta]<-\bar{L}$.  The sub-optimality of $\tilde{v}^{\mathcal{N}}_1$ immediately yields
$J_1^{\mathcal{N},{\bf b}_0}\leq \tilde{J}^{\mathcal{N},{\bf b}_0}_1< J^{\mathcal{P},{\bf b}_0}$. Therefore, Statement 1 holds for the case  $\mathbb{E}[\eta_\alpha+\eta_\beta]\neq 0$, along with the linear dependence of the difference  $J^{\mathcal{N},{\bf b}_0}_1- J^{\mathcal{P},{\bf b}_0}$ on $\mathbb{E}[\eta_1]$. Finally, for the case $\mathbb{E}[\eta_\alpha+\eta_\beta]=0>-L$ for any $L>0$, we observe from \eqref{eq:JN-tildeJP} that $J^{\mathcal{N},{\bf b}_0}_1> \tilde{J}^{\mathcal{P},{\bf b}_0}>J^{\mathcal{P},{\bf b}_0}_1$.   \hspace*{\fill} \qed

\subsection{Proof of Theorem \ref{compare_extend}}
\label{sec:app:thm:extend}
 The proof of Theorem \ref{compare_extend} stems from the continuity of the cost functional of the follower with respect to the model parameters. We proceed in  4 steps to verify this claim. In Steps 1 and 2, we consider a general model $\mathfrak{M}$ and provide explicit bounds to the functions ${\bf \Gamma}^\mathcal{I}$, ${\bf g}^\mathcal{I}$, ${\bf K}^\mathcal{I}$, ${\bf k}^\mathcal{I}$ and $l^\mathcal{I}$, which are respectively the solutions of \eqref{dynGamma}, \eqref{NumericalAll2g}, \eqref{eq:K}, \eqref{eq:k} and \eqref{eq:l}, with respect to the norm $\rho_T(\cdot)$ defined in Section \ref{sec:problem:formulation}. In Step 3, we consider a fixed base model $\mathfrak{M}_0$ and show that, for any model $\mathfrak{M}_1$ that is sufficiently close to $\mathfrak{M}_0$, the associated auxiliary functions and hence the optimal cost functional of the follower are also close to those of $\mathfrak{M}_0$. Using this, in Step 4,  we show that the comparison result for $\mathfrak{M}_0$ is preserved for any model $\mathfrak{M}$ that is close to $\mathfrak{M}_0$. 
In the sequel, we shall let $\mathcal{C}_{\Theta}>0$ be a generic constant depending only on the parameter $\Theta$ of the model $\mathfrak{M}=(\Theta,{\bf \eta})$, which may change from line to line. Similarly, we let $\mathcal{C}_{\Theta,T}>0$ be a generic constant that depends only on $\Theta$ and  $T$. 

\noindent\underline{\textit{Step 1}:}\\
 We first show that $\rho_T({\bf \Gamma}^\mathcal{I})<\infty$. Indeed, by Theorem {\color{black}\ref{wellpose:dominate}}, under Assumption \ref{ass:2:combined}, ${\bf \Gamma}^\mathcal{I}$ admits a unique solution. Since the interval $[0,T]$ is compact, ${\bf \Gamma}^\mathcal{I}$ admits a uniform bound which depends on the model parameter $\Theta$ and the duration $T$. In particular, under some additional conditions, we have the following explicit bound.
\begin{theorem}
\label{thm:gamma:norm}
   For $\mathcal{I}=\mathcal{N}$ and $\mathcal{P}$, suppose that 
   \begin{equation}
   \label{eq:gamma:norm:cond:1}
       \frac{\rho({\bf A}) + \rho({\bf C}^\mathcal{I})}{2\sqrt{\rho({\bf B}{\color{black}^\mathcal{I}})(\rho({\bf D}^\mathcal{I})+ \rho({\bf E}^\mathcal{I})(\rho({\bf A}) + \rho({\bf C}^\mathcal{I})))}} > 1
   \end{equation}%
   and 
    \begin{equation}
            \label{eq:T:cond:spectral}
            T \leq \frac{\log\left( \frac{\rho({\bf A}) + \rho({\bf C}^\mathcal{I})}{2\sqrt{\rho({\bf B}{\color{black}^\mathcal{I}})(\rho({\bf D}^\mathcal{I})+ \rho({\bf E}^\mathcal{I})(\rho({\bf A}) + \rho({\bf C}^\mathcal{I})))}} \right)}{\rho({\bf A}) + \rho({\bf C}^\mathcal{I})}.
        \end{equation}%
    Then  the equation \eqref{Gamma} admits a unique solution ${\bf \Gamma}^\mathcal{I}$, which satisfies 
        \begin{equation}
        \label{eq:Lambda}
            \rho_T({\bf \Gamma}^\mathcal{I}) \leq \Lambda^\mathcal{I} :=\sqrt{\frac{\rho({\bf D}^\mathcal{I})+ \rho({\bf E}^\mathcal{I})(\rho({\bf A}) + \rho({\bf C}^\mathcal{I}))}{ \rho({\bf B}{\color{black}^\mathcal{I}})} }  .
        \end{equation}%
\end{theorem}

\begin{proof}
    Denote the dimension of ${\bf \Gamma}^\mathcal{I}$ by {\color{black}$m_\mathcal{I}$} and define $\mathcal{L}{\color{black}^\mathcal{I}}: C([0,T];\mathbb{R}^{{\color{black} m_\mathcal{I}}}) \to C([0,T];\mathbb{R}^{{\color{black} m_\mathcal{I}}})$
    by $\mathcal{L}{\color{black}^\mathcal{I}}({\bf \Phi}^\mathcal{I})={\bf \Psi}^\mathcal{I}$, where ${\bf \Psi}^\mathcal{I}$ is the solution of the equation 
    \begin{equation}
     \label{eq:Gtilde:T}
       \left\{ \begin{aligned}
         -\frac{d{\bf \Psi}^\mathcal{I}(t)}{dt} =&\  {\bf \Psi}^\mathcal{I}(t){\bf A} +    {\bf C}^\mathcal{I}{\bf \Psi}^\mathcal{I}(t) - {\bf \Phi}^\mathcal{I}(t) {\bf B}{\color{black}^\mathcal{I}}{\bf \Phi}^\mathcal{I}(t)  +  {\bf D}^\mathcal{I}, \\
          {\bf \Psi}(T)=&\ {\bf E}^\mathcal{I}.
        \end{aligned}\right.
    \end{equation}%
   With $\Lambda^\mathcal{I}$ as in \eqref{eq:Lambda}, consider the closed ball 
     \begin{equation*}
            \mathcal{B}_{\Lambda^\mathcal{I}}  := \left\{ {\bf \Phi} \in C([0,T];\mathbb{R}^{{\color{black} m_\mathcal{I}}}) : \sup_{0\leq t\leq T} \rho({\bf \Phi(t)}) \leq  \Lambda^\mathcal{I} \right\}.
        \end{equation*}%
    We shall show that $\mathcal{L}{\color{black}^\mathcal{I}}(\mathcal{B}_{\Lambda^\mathcal{I}}) \subseteq \mathcal{B}_{\Lambda^\mathcal{I}}$. Let ${\bf \Phi}^\mathcal{I} \in\mathcal{B}_{\Lambda^\mathcal{I}}$ and ${\bf \Psi}^\mathcal{I} = \mathcal{L}{\color{black}^\mathcal{I}}({\bf \Phi}^\mathcal{I})$, which can be represented as 
    \begin{align*}
        {\bf \Psi}^\mathcal{I}(t) &= \int_t^T e^{{\bf C}^\mathcal{I}(s-t)} ({\bf D}^\mathcal{I} - {\bf \Phi}^\mathcal{I}(s){\bf B}{\color{black}^\mathcal{I}}{\bf \Phi}^\mathcal{I}(s)) e^{{\bf A}(s-t)}ds+e^{{\bf C}^\mathcal{I}(T-t)} {\bf E}^\mathcal{I}  e^{{\bf A}(T-t)}.
    \end{align*}%
   Hence, by  \eqref{eq:T:cond:spectral} and the fact that $ {\bf \Phi}^\mathcal{I} \in  \mathcal{B}_{\Lambda^\mathcal{I}}$, we have 
    \begin{align*}
        \sup_{0\leq t\leq T}\rho( {\bf \Psi}^\mathcal{I}(t) ) 
       &\leq   \left(\rho({\bf D}^\mathcal{I}) + \rho({\bf B}{\color{black}^\mathcal{I}}) \sup_{0\leq t\leq T}\rho^2({\bf \Phi}^\mathcal{I}(t))\right) \times  \int_0^T e^{(\rho({\bf A}) + \rho({\bf C}^\mathcal{I}) )s} ds \\ &\ +\rho({\bf E}^\mathcal{I})e^{(\rho({\bf A}) + \rho({\bf C}^\mathcal{I}))T} 
      \leq  \Lambda^\mathcal{I}.
    \end{align*}%
   Therefore, ${\bf \Psi}{\color{black}^\mathcal{I}}\in \mathcal{B}_{\Lambda{\color{black}^\mathcal{I}}}$.  Next, we shall show that $\mathcal{L}{\color{black}^\mathcal{I}}(\mathcal{B}_{\Lambda^\mathcal{I}})$ is relatively compact with respect to the norm $\rho_T(\cdot)$.  Consider a sequence $\{ {\color{black}{\bf \Psi}^\mathcal{I}_n} \}_{n=1}^\infty \subset \mathcal{L}{\color{black}^\mathcal{I}}( \mathcal{B}_{\Lambda^\mathcal{I}})$. Since $\mathcal{L}{\color{black}^\mathcal{I}}( \mathcal{B}_{\Lambda^\mathcal{I}}) \subseteq  \mathcal{B}_{\Lambda^\mathcal{I}}$, we have $\sup_{n\in\mathbb{N}} \rho_T({\color{black}{\bf \Psi}^\mathcal{I}_n}) <\Lambda^\mathcal{I}$. By \eqref{eq:Gtilde:T}, we also have $\sup_{n\in\mathbb{N}} \rho_T(\frac{d{\color{black}{\bf \Psi}^\mathcal{I}_n}}{dt}) < \infty $, whence the sequence is uniformly equicontinuous. The relatively compactness now follows from the Arzel\'a-Ascoli theorem. It is easy to check that $ \mathcal{B}_{\Lambda^\mathcal{I}}$ is closed and convex. Therefore,  by Schauder's fixed point theorem, the map $\mathcal{L}{\color{black}^\mathcal{I}}$ admits a fixed point ${\bf \Gamma}{\color{black}^\mathcal{I}}\in \mathcal{B}_{\Lambda^\mathcal{I}}$, which is thus the solution of \eqref{dynGamma}. The uniqueness of solution, on the other hand, follows from the local Lipschitz property of \eqref{dynGamma}.
 \end{proof}

 \noindent\underline{\textit{Step 2}:}\\
Next, we provide explicit bounds for ${\bf g}^\mathcal{I}$, ${\bf K}^\mathcal{I}$, ${\bf k}^\mathcal{I}$ and $l^\mathcal{I}$ associated with $\mathfrak{M}=(\Theta,\boldsymbol{\eta})$, in terms of $\rho_T({\bf \Gamma}^\mathcal{I})$ obtained in Step 1. By applying a Gr\"onwall-type inequality to the equations \eqref{NumericalAll2g}, \eqref{eq:K}, \eqref{eq:k} and \eqref{eq:l}, it is easy to see that 
\begin{equation}
\label{eq:asym:step2}
\begin{aligned}
\rho_T({\bf g}^\mathcal{I}) \le&\ \mathcal{C}_\Theta e^{\mathcal{C}_\Theta(1+\rho_T({\bf \Gamma}^\mathcal{I}))T},\\
\rho_T({\bf K}^\mathcal{I}) \le&\ \mathcal{C}_\Theta(1+\rho_T({\bf \Gamma}^\mathcal{I}))e^{\mathcal{C}_\Theta(1+\rho_T({\bf \Gamma}^\mathcal{I}))T},\\
\rho_T({\bf k}^\mathcal{I}) \le&\ \mathcal{C}_\Theta(1+\rho_T({\bf g}^\mathcal{I})(\rho_T({\bf \Gamma}^\mathcal{I})+\rho_T({\bf K}^\mathcal{I}))) e^{\mathcal{C}_\Theta(1+\rho_T({\bf \Gamma}^\mathcal{I}))T}, \\ 
\rho_T(l^\mathcal{I}) \le&\ \mathcal{C}_\Theta T\left(1+\rho_T({\bf g}^\mathcal{I})\left(1+\rho_T({\bf k}^\mathcal{I})+\rho_T({\bf g}^\mathcal{I})\right) \right).
\end{aligned}
\end{equation}%
Thus, with a bound for ${\bf \Gamma}^\mathcal{I}$, the bounds for  ${\bf g}^\mathcal{I}$, ${\bf K}^\mathcal{I}$, ${\bf k}^\mathcal{I}$ and $l^\mathcal{I}$ can be obtained accordingly. 

 \noindent\underline{\textit{Step 3}:}\\
We consider two models $\mathfrak{M}_0=(\Theta_0,\boldsymbol{\eta}_0)$ and $\mathfrak{M}_1=(\Theta_1,\boldsymbol{\eta}_1)$, such that $\mathcal{D}(\mathfrak{M}_0,\mathfrak{M}_1)<\delta$ for some $\delta>0$,  and $\mathfrak{M}_0$ satisfies Assumption \ref{ass:2:combined}. In the sequel, we denote by ${\bf M}_i$ the associated matrix or function ${\bf M}$ for the model $\mathfrak{M}_i$, $i=0,1$.
In this step, we shall provide uniform bounds for the functions ${\bf \tilde{\Gamma}}^\mathcal{I}$, ${\bf \tilde{g}}^\mathcal{I}$, ${\bf \tilde{K}}^\mathcal{I}$, $\tilde{{\bf k}}^\mathcal{I}$ and $\tilde{l}^\mathcal{I}$ in terms of $\delta$, $\Theta_0$ and $T$. 
First, using \eqref{dynGamma}, it is easy to see that 
\begin{equation}
\label{eq:tilde:gamma:compare}
   \left\{ \begin{aligned}
        -\frac{d{\bf \tilde{\Gamma}^\mathcal{I}}(t)}{dt} =&\ {\bf \tilde{D}}^\mathcal{I} + \tilde{{\bf \Gamma}}^\mathcal{I}(t)\left({\bf A}_1 -{\bf B}_1^{{\color{black}\mathcal{I}}} {\bf \Gamma}_0^\mathcal{I}(t) \right) +  \left({\bf C}_1^\mathcal{I}  -{\bf \Gamma}^\mathcal{I}_0(t) {\bf B}_1^{{\color{black}\mathcal{I}}}  \right){\bf \tilde{\Gamma}}^\mathcal{I}(t)\\
        &\ - \tilde{{\bf \Gamma}}^\mathcal{I}(t){\bf B}_1^{{\color{black}\mathcal{I}}}\tilde{{\bf \Gamma}}^\mathcal{I}(t)  + {\bf \Gamma}_0^\mathcal{I}(t){\bf \tilde{A}}  + {\bf \tilde{C}}^\mathcal{I}{\bf \Gamma}_0^\mathcal{I}(t), \\ 
        \tilde{{\bf \Gamma} }(T)= &\ \tilde{{\bf E}}^\mathcal{I}.
    \end{aligned}\right.
\end{equation}%
We proceed to show the well-posedness of \eqref{eq:tilde:gamma:compare} and that $\rho_T({\bf \tilde{\Gamma}}^\mathcal{I})\leq O(\min\{\delta,\sqrt{\delta}\})$ when $\delta>0$ is sufficiently small. In particular, we can pick $\delta>0$ such that, for $\mathcal{I}=\mathcal{N}$ and $\mathcal{P}$,  $\mathcal{E}_0^\mathcal{I}(\delta)>1$,   $\rho({\bf B}^{{\color{black}\mathcal{I}}}_0+{\bf \tilde{B}}^{{\color{black}\mathcal{I}}}) > 0$ and \eqref{eq:T:extension} holds, where $\mathcal{E}^\mathcal{I}_0(\delta)$ is defined in \eqref{eq:E:delta} with ${\bf A}_1, {\bf B}^{{\color{black}\mathcal{I}}}_1, {\bf C}_1^\mathcal{I}, {\bf E}^\mathcal{I}_1$ in place of ${\bf A}, {\bf B}^{{\color{black}\mathcal{I}}}, {\bf C}^\mathcal{I}, {\bf E}^\mathcal{I}$.  Notice that the above is possible, since $\mathcal{E}^\mathcal{I}_0(\delta) = O(\delta^{-\frac{1}{2}})$. By following the proof of Theorem \ref{thm:gamma:norm}, we can readily show that 
    \begin{align}
         \label{eq:tilde:gamma}
        \rho({\bf \tilde{\Gamma}}^\mathcal{I}) \leq &\ \sqrt{ \rho^{-1}({\bf B}^{{\color{black}\mathcal{I}}}_1) \Big( \rho_T\big({\bf \Gamma}^\mathcal{I}_0\tilde{{\bf A}} + \tilde{{\bf C}}^\mathcal{I}{\bf \Gamma}^\mathcal{I}_0 \big) + \rho({\bf \tilde{E}})   \times \Big( \rho_T\big({\bf A}_1-  {\bf B}^{{\color{black}\mathcal{I}}}_1{\bf \Gamma}^\mathcal{I}_0\big)  +   \rho_T\big({\bf C}^\mathcal{I}_1 - {\bf \Gamma}^\mathcal{I}_0 {\bf B}^{{\color{black}\mathcal{I}}}_1\big)   \Big) \Big) }\nonumber\\
    =&\ O(\sqrt{\delta}).
    \end{align}%
 By \eqref{eq:tilde:gamma} and Theorem {\color{black}\ref{wellpose:dominate}}, we have $\rho({\bf \Gamma}^\mathcal{I}_1) \leq \mathcal{C}_{\Theta_0,T}(1+\sqrt{\delta})$. Using this  and rewriting  \eqref{eq:tilde:gamma:compare} as 
\begin{equation}
   \left\{ \begin{aligned}
        -\frac{d{\bf \tilde{\Gamma}^\mathcal{I}}(t)}{dt} =&\ {\bf \tilde{D}}^\mathcal{I} + \tilde{{\bf \Gamma}}^\mathcal{I}(t)\left({\bf A}_0  -{\bf B}^{{\color{black}\mathcal{I}}}_0 {\bf \Gamma}_0^\mathcal{I}(t) \right) +  \left({\bf C}_0^\mathcal{I} -{\bf \Gamma}^\mathcal{I}_1(t){\bf B}^{{\color{black}\mathcal{I}}}_1\right){\bf \tilde{\Gamma}}^\mathcal{I}(t)\\
        &\ +{\bf \Gamma}_1^\mathcal{I}(t) {\bf \tilde{B}}^{{\color{black}\mathcal{I}}}{\bf \Gamma}^\mathcal{I}_0(t)  + {\bf \Gamma}_1^\mathcal{I}(t){\bf \tilde{A}}  + {\bf \tilde{C}}^\mathcal{I}{\bf \Gamma}_1^\mathcal{I}(t),  \\ 
        \tilde{{\bf \Gamma} }(T)= &\ \tilde{{\bf E}}^\mathcal{I},
    \end{aligned}\right.
\end{equation}%
we can arrive, by a Gr\"onwall-type inequality, that 
    \begin{equation}
        \label{eq:tilde:gamma:2}
        \begin{aligned}
         \rho({\bf \tilde{\Gamma}}^\mathcal{I}) \leq &\ \mathcal{C}_{\Theta_0,T}\delta(1+\sqrt{\delta})e^{\mathcal{C}_{\Theta_0,T}(1+\sqrt{\delta})T} =O(\delta). 
         \end{aligned}
    \end{equation}%
Next, using the equation satisfied by $\tilde{{\bf K}}^\mathcal{I}$: 
    \begin{equation}
    \left\{\label{eq:K:tilde}
        \begin{aligned}
            -\frac{d\tilde{{\bf K}}^\mathcal{I}(t)}{dt} =&\ {\bf \tilde{Q}} + {\bf \tilde{R}}^\mathcal{I}(t) + \tilde{{\bf K}}^\mathcal{I}(t)\left({\bf A}_1-{\bf B}^{{\color{black}\mathcal{I}}}_1 {\bf \Gamma}^\mathcal{I}_1 (t)\right) + \left({\bf A}_1- {\bf B}^{{\color{black}\mathcal{I}}}_1{\bf \Gamma}^\mathcal{I}_1(t)\right)^\top {\bf \tilde{K}}^\mathcal{I}(t)  \\
          &\   + {\bf K}_0^\mathcal{I}(t)\left(\tilde{\bf A} - \tilde{\bf B}^{{\color{black}\mathcal{I}}}{\bf \Gamma}^\mathcal{I}_0(t) - {\bf B}^{{\color{black}\mathcal{I}}}_1 \tilde{{\bf \Gamma}}^\mathcal{I}(t) \right)   + \left(\tilde{\bf A} - \tilde{\bf B}^{{\color{black}\mathcal{I}}}{\bf \Gamma}^\mathcal{I}_0(t) - {\bf B}^{{\color{black}\mathcal{I}}}_1\tilde{\Gamma}^\mathcal{I}(t) \right)^\top {\bf K}_0^\mathcal{I}(t)  , \\
            \tilde{{\bf K}}^\mathcal{I}(T)=&\ \bar{{\bf Q}}_1-\bar{{\bf Q}}_0,
        \end{aligned}\right.
    \end{equation} %
along with the estimate of   $\rho_T({\bf K}_0^\mathcal{I})$ in \eqref{eq:asym:step2}, we have 
     \begin{align*}
    \rho_T({\bf \tilde{K}}^\mathcal{I}) \le&\  \mathcal{C}_{\Theta_0}\bigg( \rho_T({\bf K}_0^\mathcal{I})\left(\delta+\delta\rho_T({\bf \Gamma}^\mathcal{I}_0)+\rho_T({\bf \tilde{\Gamma}}^\mathcal{I})\right)  + \delta\left(1 + \rho^2_T({\bf \Gamma}^\mathcal{I}_0)\right) + \rho_T({\bf \Gamma}^\mathcal{I}_0)\rho_T({\bf \tilde{\Gamma}}^\mathcal{I})\bigg)  \\
    &\ \times 
    e^{\mathcal{C}_{\Theta_0}(1+\delta + \rho_T({\bf \Gamma}^\mathcal{I}_0) +  \rho_T({\bf \tilde{\Gamma}}^\mathcal{I}) )T}  \leq  \mathcal{C}_{\Theta_0,T}\delta. 
\end{align*}%
 Similarly, by  a Gr\"onwall-type inequality and \eqref{eq:asym:step2}, we have 
 \begin{equation}
    \label{eq:tilde:estimates}
        \max\left\{ \rho_T({\bf \tilde{g}}^\mathcal{I}),\rho_T({\bf \tilde{k}}^\mathcal{I}), \rho_T(\tilde{l}^\mathcal{I})  \right\} \leq  \mathcal{C}_{\Theta_0,T}\delta. 
    \end{equation}%
\noindent\underline{\textit{Step 4}:}\\
Let $J^{\mathcal{I},i}_1$, $\mathcal{I}=\mathcal{N},\mathcal{P}$ and $i=0,1$, be the cost functional of the follower under game $\mathcal{I}$ and model $\mathfrak{M}_i$. Using the estimates $\rho(\tilde{{\bf \Gamma}}^\mathcal{I}) \leq O(\delta)$,  \eqref{eq:tilde:estimates}, along with the expression \eqref{eq:lem:cost}, one can readily show, by a standard telescoping approach, that  $|J_1^{\mathcal{I},0}-J_1^{\mathcal{I},1}|\leq \mathcal{C}_{\Theta_0,\boldsymbol{\eta}_0,T}\delta$,
where the constant $\mathcal{C}_{\Theta_0,\boldsymbol{\eta}_0,T}  >0$ depends solely on $\Theta_0,\boldsymbol{\eta}_0$ and $T$. Now, suppose that $\Delta J^0_1:=J^{\mathcal{N},0}_1-J^{\mathcal{P},0}_1>0$. Then for any $\mathfrak{M}$ with
$\delta:=\mathcal{D}(\mathfrak{M}_0,\mathfrak{M}) <  \frac{|\Delta J^0_1|}{2\mathcal{C}_{\Theta_0,\boldsymbol{\eta}_0,T}}$, 
we have $|J^{\mathcal{N},0}_1-J^{\mathcal{N},1}_1|<\frac{\Delta J^0_1}{2}$ and   $|J^{\mathcal{P},0}_1-J^{\mathcal{P},1}_1|<\frac{\Delta J^0_1}{2}$, whence $J_1^{\mathcal{N},1} - J_1^{\mathcal{P},1}>0$.
Similarly, we can prove that $J_1^{\mathcal{N},1} - J_1^{\mathcal{P},1}<0$ whenever $J^{\mathcal{N},0}_1-J^{\mathcal{P},0}_1<0$ and  $\mathcal{D}(\mathfrak{M}_0,\mathfrak{M}) <\delta$.
\begin{flushright}
    \qed
\end{flushright}

\end{appendices}

 \end{document}